\documentclass{article}
\DeclareFontFamily{U}{mathb}{\hyphenchar\font45}
\DeclareFontShape{U}{mathb}{m}{n}{
      <5> <6> <7> <8> <9> <10> gen * mathb
      <10.95> mathb10 <12> <14.4> <17.28> <20.74> <24.88> mathb10
      }{}
\DeclareSymbolFont{mathb}{U}{mathb}{m}{n}
\DeclareFontSubstitution{U}{mathb}{m}{n}
\DeclareMathSymbol{\sqsubsetneq}  {0}{mathb}{"88}
\DeclareMathSymbol{\sqsubseteq} {0}{mathb}{"84}
\newcommand{\notsqsubseteq}{\mathrel{\mathrlap{\not}\sqsubseteq}}
\usepackage{amsmath, amssymb}
\usepackage[ style=numeric]{biblatex}
\let\cite\cites
\usepackage{bm}
\usepackage{mathtools}
\usepackage{fancyhdr}
\usepackage[margin=1in]{geometry}
\usepackage{bbm}
\usepackage{circuitikz}
\usepackage{setspace}
\usepackage{pdfpages}
\usepackage{graphicx}
\usepackage{enumitem}
\usepackage{amsthm}
\usepackage{tikz}
\usepackage{tikz-cd}
\newtheoremstyle{dotless}{}{}{}{}{\bfseries}{}{ }{}
\theoremstyle{dotless}
\renewenvironment{proof}{{\bfseries Proof:}}{\qed\vspace{1mm}}
\newtheorem{thm}{Theorem}[section]
\newtheorem{defin}[thm]{Definition}
\newtheorem{lemma}[thm]{Lemma}
\newtheorem{prop}[thm]{Proposition}
\newtheorem{obs}[thm]{Observation}
\newtheorem{cor}[thm]{Corollary}
\newtheorem{q}[thm]{Question}
\newtheorem{C}[thm]{Construction}
\usepackage{hyperref}
\hypersetup{
    colorlinks=true,
    linkcolor=blue,
    citecolor=blue,
    filecolor=blue, 
    urlcolor=blue,       
    pdftitle={Overleaf Example},
    pdfpagemode=FullScreen,
    }
\usepackage{subfiles}
\title{Reverse mathematics of regular countable second countable spaces}
\author{Giorgio G. Genovesi\\
School of mathematics, University of Leeds }

\date{}
\begin{document}
\def\null\vskip 2em{}\maketitle
\begin{abstract}
\noindent
We study the reverse mathematics of characterization theorems of regular countable second countable (CSC) spaces. We prove that arithmetic comprehension is equivalent over $\mathbf{RCA}_0$ to every $T_3$ CSC space being metrizable, and we characterize the $T_3$ spaces which are metrizable over $\mathbf{RCA}_0$.  We show that Lynn's theorem for CSC spaces can be carried out in $\mathbf{ACA}_0$, namely that every zero dimensional separable space is homeomorphic to the order topology of a linear order. We also show that arithmetic comprehension is equivalent to every $T_2$ compact CSC space being well-orderable.
From general topology, we know that the locally compact $T_2$ CSC spaces are the well-orderable CSC spaces, and that the $T_3$ scattered CSC spaces are the completely metrizable CSC spaces. We show that these characterizations and a few others are equivalent to arithmetic transfinite recursion over $\mathbf{RCA}_0$.  We also find a few statements that are equivalent to $\Pi^1_1$ comprehension. In particular we show that every $T_3$ CSC space has a Cantor-Bendixson rank and that every $T_3$ CSC space is the disjoint union of a scattered space and dense in itself space are equivalent to $\Pi^1_1$ comprehension over $\mathbf{RCA}_0$.
\end{abstract}
\setcounter{section}{-1}
\begin{spacing}{0.1}
\tableofcontents
\end{spacing}
\newpage

\section{Introduction}
In this work, we study the reverse mathematics of theorems characterizing regular countable second countable spaces. Reverse mathematics, in a very broad sense, is the study of the necessary and sufficient axioms needed to prove a given theorem. Most of the time, the given theorem is one that can be expressed in the seemingly poor language of second order arithmetic. Over a suitable base theory, one can use a series of codes to talk about pairs, $n$-tuples, finite sets and sequences, functions, $k$-ary relations, countable rings and groups, complete separable metric spaces and Banach spaces, and many other mathematical objects. This allows us to express, in the language of second order arithmetic, many theorems of ordinary mathematics.\\
\\
However, the language of second order arithmetic is not suited to consider objects of arbitrary large cardinality. One can therefore only hope to formalize a modest part of general topology. One approach is the study of MF spaces done by Mummert \cite{Mummert}, in which the points of a space are identified with second order objects. This approach is similar to how Polish spaces are formalized in second order arithmetic. Another approach is to consider Countable Second Countable (CSC) spaces. Friedman and Hirst studied CSC spaces in the form of countable subsets of complete separable metric spaces \cite{First}{Hirst2}.  Friedman also considered CSC spaces in the form of countable metric spaces,  countable linear orders with the order topology, and subspaces of the rationals \cite{FriedC}. Greenberg and Montalb{\'a}n looked at $T_2$ compact CSC spaces mostly over $\textbf{ACA}_0$, which they referred to as \textit{very countable spaces}. The study of CSC spaces in general over $\textbf{RCA}_0$ started with Dorais \cite{Dorais}, who showed that arithmetic comprehension is needed for many basic topological facts. For example, in the absence of arithmetic comprehension, there can be compact infinite discrete spaces and products of compact spaces that are not compact.\\
\\
Dorais' formalization of CSC spaces has been looked at from different perspectives.
Frittaion, Hendtlass, Marcone, Shafer, and Van der Meeren looked at Noetherian CSC spaces \cite{Frittaion}, which arise naturally as the Alexandroff topology of well quasi orders. Countable linear orders with their order topology  were studied by Shafer \cite{Shafer}. More recently, the Ginsburg-Sands theorem for CSC spaces, a Ramsey type principle, has been studied by Benham, Delapo, Dzhafarov, Solomon, and Villano \cite{Benham}. In their work they showed that every infinite $T_1$ CSC space has an infinite subspace which is either discrete or has the cofinite topology is strictly between the strength of the Ramsey theorem for pairs and $\textbf{ACA}_0$.\\
\\
It is natural to ask what are the conditions for metrizability. Understanding which CSC spaces are metrizable or completely metrizable allows us to frame Friedman and Hirst's results about countable subspaces of complete metric spaces in terms of Dorais' formalization of CSC spaces. A sketch of a metrization proof is given in Greenberg and Montalb{\'a}n \cite{Greenberg}, which follows the classic proof of the Urysohn metrization theorem. However, the proof seems to work only for compact $T_2$ CSC spaces. We prove that over $\textbf{RCA}_0$, the metrizability of effectively $T_2$ regular CSC spaces is equivalent to arithmetic comprehension. 
Lynn \cite{Lynn} proved that the zero dimensional subspaces of $\mathbb{R}^n$ are linearly orderable. We show that this theorem for $T_3$ CSC spaces can be carried out over $\textbf{ACA}_0$, and we get a partial reversal. We will also prove that $T_2$ compact CSC spaces are well-orderable over $\textbf{ACA}_0$ and give a partial characterization of the effectively $T_2$ effectively compact CSC spaces over $\textbf{RCA}_0$ in terms of well-orders. \\
\\
In the second part of this work, we will consider scattered regular CSC spaces and locally compact $T_2$ CSC spaces. We have that scattered $T_2$ regular CSC spaces are precisely the countable completely metrizable spaces, and locally compact $T_2$ CSC spaces are the well-orderable spaces. Both of these characterizations turn out to be equivalent to arithmetic transfinite recursion. The proofs of both these characterizations rely on the fact that arithmetic transfinite recursion is equivalent to every locally compact $T_2$ CSC space has a choice of compact neighborhoods. We will also be able to give alternative proofs to some results by Friedman and Hirst \cite{First}{Hirst2}.
Finally, we show that every $T_2$ regular  CSC space has a rank and every $T_2$ regular CSC space is the disjoint union of a dense in itself space and a scattered space are both equivalent to $\Pi^1_1$ comprehension over $\textbf{RCA}_0$. We will also show that being compact, being connected, being scattered, and being a complete metric space are all $\Pi^1_1$ universal formulas over $\textbf{ACA}_0$.

\part{Metrization and orderability of CSC spaces}
In this part we will study the conditions needed for a space to be metrizable or orderable over $\mathbf{RCA}_0$. We will also study the conditions for orderability. A large portion of this part will be dedicated to introducing basic known facts and definitions as well as presenting lemmas which, while being technical and not of immediate interest, will be needed for the next part. 
\section{Notation}

For a poset $(P,\leq_P)$ and $S\subseteq P$ we write:
\begin{equation*}
{\uparrow} S =\{p\in P:\exists s\in S \,(s\leq_P p)\}
\end{equation*}
\begin{equation*}
{\downarrow} S =\{p\in P:\exists s\in S\,( p\leq_P s)\}
\end{equation*}
We call ${\uparrow} S$ the upwards closure of $S$ in $P$ and ${\downarrow} S$ the downwards closure of $S$ in $P$.\\
\\
We use the French notation for intervals. That is, for a linear order $(L,<_L)$ and $a,b\in L$ we write:

\begin{equation*}
\mathopen[a,b\mathclose]_{<_L} =\{l\in L: a\leq_L l \leq_L b\}
\end{equation*}
\begin{equation*}
\mathopen]a,b\mathclose]_{<_L} =\{l\in L: a<_L l \leq_L b\}
\end{equation*}
\begin{equation*}
\mathopen[a,b\mathclose[_{<_L} =\{l\in L: a\leq_L l <_L b\}
\end{equation*}
\begin{equation*}
\mathopen]a,b\mathclose[_{<_L} =\{l\in L: a<_L l <_L b\}
\end{equation*}
If $L$ has a minimal element, we will denote it usually by $0$. We use the convention that $-\infty$ denotes a new element strictly smaller than any element of $L$ and $+\infty$ or $\infty$ to denote an element that is strictly greater than any element of $L$. By $L+1$, we mean  the linear order with field $L\cup\{\infty\}$ and the order is extended so that $\infty$ is greater than any other element of $L$. \\
\\
By $\mathbb{N}^{<\mathbb{N}}$ we mean the set of all finite sequences of numbers and by $2^{<\mathbb{N}}$ we mean the set of all finite $0$-$1$ valued sequences. For $\sigma,\tau\in \mathbb{N}^{<\mathbb{N}}$ we write $\sigma\sqsubseteq \tau$ if $|\sigma|\leq |\tau|\wedge \forall i<|\sigma|\;(\sigma(i)=\tau(i))$. We will also use the following conventions:
\begin{enumerate}[label={(\arabic*)}]

\item $2^{=n}=\{\sigma\in 2^{<\mathbb{N}}:|\sigma|=n\}$ and $\mathbb{N}^{=n}=\{\sigma\in \mathbb{N}^{<\mathbb{N}}:|\sigma|=n\}$
\item $2^{\leq n}=\{\sigma\in 2^{<\mathbb{N}}:|\sigma|\leq n\}$ and $\mathbb{N}^{\leq n}=\{\sigma\in \mathbb{N}^{<\mathbb{N}}:|\sigma|\leq n\}$
\end{enumerate} Given a sequence $\sigma\in \mathbb{N}^{<\mathbb{N}}$ we write
\begin{equation*} 
[[\sigma]]=\{f:\mathbb{N}\rightarrow \mathbb{N}: \forall i<|\sigma|\, f(i)=\sigma(i)\}
\end{equation*}
A tree is a subset $T$ of $\mathbb{N}^{<\mathbb{N}}$ such that $T={\downarrow} T$. We say that $f:\mathbb{N}\rightarrow \mathbb{N}$ is a branch in $T$ if for all $n\in \mathbb{N}$ we have $f|_{<n}\in T$. By $[T]$, we mean the set of branches of $T$. We say that a tree is well-founded if $[T]\neq \emptyset$.\\
\\
Given $A\subseteq \mathbb{N}$ and $e,x\in \mathbb{N}$ we write $\Phi^A_e(x)$ to mean the Turing machine with index $e$, oracle $A$, and input $x$. We write $\Phi^A_e(x){\downarrow}_{\leq t}$ to mean that the Turing machine with index $e$ and oracle $A$ halts at input $x$ in less than $t$ steps and we write $\Phi^A_e(x){\downarrow}$ to mean that there exists a $t$ such that $\Phi^A_e(x){\downarrow}_{\leq t}$ and $\Phi^A_e(x){\uparrow}$ if such $t$ doesn't exist. We write $B\leq_T A$ if there exists an $e\in\mathbb{N}$ such that $\Phi^A_e$ is the characteristic function of $B$.
\\
\\
We will also assume that there is a fixed $\Delta^0_0$ definable pairing function for numbers and we use the Ackermann interpretation to encode finite sets. That is, we have that $x\in_{Ack} y$ if the $x$-th digit of the binary representation of $y$ is $1$. We will not distinguish between finite sets and codes for finite sets. We will use the fact that membership in the Ackermann interpretation is $\Delta^0_1$ definable. More details about coding in arithmetic can be found in \cite{Mofa}.
\section{Reverse mathematics}
We briefly introduce the Big Five systems of reverse mathematics and some classic theorems. More details can be found in \cite{Simp},\cite{Hirschfeldt}, and \cite{Rmprp}.\\
\\
$\textbf{RCA}_0$ is the system consisting of the axioms of a discretely ordered commutative semiring with identity, induction for $\Sigma^0_1$ formulas, and comprehension for $\Delta^0_1$ predicates.\\
\\
$\textbf{WKL}_0$ is the system $\textbf{RCA}_0$ plus the statement that any infinite subtree of the binary tree $2^{<\mathbb{N}}$ has an infinite branch.\\
\\
$\textbf{ACA}_0$ is the system $\textbf{RCA}_0$ plus comprehension for all arithmetical formulas.\\
\\
$\textbf{ATR}_0$ is the system $\textbf{RCA}_0$ plus arithmetic transfinite recursion, which is the statement that for any well-order $(L,<_L)$ with least element $0$, any set $X$, and any arithmetical formula $\varphi$ there exists a sequence of sets $(X_j)_{j\in L}$ such that $X_0=X$ and for all $j\in L\setminus\{0\}$
\begin{equation*}
X_j=\{n\in\mathbb{N}:\varphi(n, (X_i)_{i<_L \, j})\}.
\end{equation*}
Informally, we can view $X_j$ as being the set obtained by applying some arithmetic procedure given by $\varphi$ to $X$ $ j$ many times.\\
\\
$\Pi^1_1 \textbf{-CA}_0$ is the system $\textbf{RCA}_0$ plus comprehension for $\Pi^1_1$ formulas.
\begin{prop} $\textbf{RCA}_0$ proves the following:
\begin{enumerate}[label={(\arabic*)}]
\item For any strictly increasing function $f:\mathbb{N}\rightarrow \mathbb{N}$, $\text{rng}(f)$ exists.
\item $\textbf{B}\Sigma^0_1$: for every $\Sigma^0_1$ formula $\varphi$ we have
\begin{equation*}
\forall a\;(\forall x<a\; \exists y\;\varphi(x,y))\rightarrow (\exists b\;\forall x<a\;\exists y<b \;\varphi(x,y))
\end{equation*}
(See \cite[Exercise II.3.14]{Simp} or \cite[Theorem 1.44]{Mofa}).
\end{enumerate}
\end{prop}

\begin{prop} \label{ACAnot} Over $\textbf{RCA}_0$ the following are equivalent:
\begin{enumerate}[label={(\arabic*)}]
\item Arithmetic comprehension.
\item For every injective function $f:\mathbb{N}\rightarrow \mathbb{N}$, the range of $f$ exists (See \cite[Lemma III.1.3]{Simp}).
\item For every set $X$, the Turing jump $X'=\{e\in\mathbb{N}:\Phi^X_e(e){\downarrow}\}$ exists (See \cite[Corollary 5.6.3]{Rmprp}).
\item K{\H o}nig's lemma: every infinite, finitely branching tree $T\subseteq \mathbb{N}^{<\mathbb{N}}$ has a branch (See \cite[Theorem III.7.2]{Simp}).
\end{enumerate}
\end{prop}

\begin{thm} \label{ATRnot} Over $\textbf{RCA}_0$ the following are equivalent:
\begin{enumerate}[label={(\arabic*)}]
\item Arithmetic transfinite recursion.
\item For every sequence of trees $(T_i)_{i\in\mathbb{N}}$ such that $\forall i\in\mathbb{N}\;|[T_i]|\leq 1$, the set $\{i\in\mathbb{N}: [T_i]\neq \emptyset\}$ exists (See \cite[Theorem V.5.2]{Simp}).
\item Weak comparability of well-orders: for any two well-orders $L$ and $W$, there either exists an increasing function from $L$ to $W$ or an increasing function from $W$ to $L$ (See \cite[Section 12.1]{Rmprp}).
\item Strong comparability of well-orders: for any two well-orders $L$ and $W$, either $L$ is isomorphic to an initial segment of $W$ or $W$ is isomorphic to an initial segment of $L$ (See \cite[Section V.6]{Simp}).
\end{enumerate}
\end{thm}

\begin{thm}(See \cite[Theorem V.8.3]{Simp}, \cite [Corollary 12.1.14]{Rmprp}) $\textbf{ATR}_0$ proves $\Sigma^1_1$ choice, or rather for every $\Sigma^1_1$ formula $\varphi$ we have
\begin{equation*}
\forall n\; \exists X\; \varphi(n,X)\rightarrow \exists (X_n)_{n\in\mathbb{N}}\;\forall n\;\varphi(n,X_n)
\end{equation*}
\end{thm}

\begin{thm}(See \cite[Section V.1]{Simp}) For any $\Sigma_1^1$ formula $\varphi(X)$ there exists a $\Delta^0_0$ formula $\theta(\sigma,\tau)$ such that $\textbf{ACA}_0$ proves that
\begin{equation*}
\forall X\; \varphi(X)\leftrightarrow \exists f\;\forall m\;\theta(X|_{\leq m},f|_{\leq m})
\end{equation*}
the formula to the righthand side of the biconditional is called the Kleene normal form of $\varphi$.
\end{thm}

\begin{obs} The formula $\theta$ in the previous theorem can be seen as defining a tree with respect to $X$. So over $\textbf{ACA}_0$ being a ill founded tree is a universal $\Sigma^1_1$ formula and being a well-founded tree is a universal $\Pi^1_1$ formula.
\end{obs}

\begin{defin} Let $T\subseteq \mathbb{N}^{<\mathbb{N}}$ be a tree. The Kleene Brouwer order on $T$ is given by
\begin{equation*}
\sigma<_{\text{KB}}\tau\leftrightarrow \tau\sqsubseteq \sigma\vee \exists j\;( \forall i<j\;\sigma(i)=\tau(i)\wedge \sigma(j)<\tau(j))
\end{equation*}
$\textbf{ACA}_0$ proves that every $T$ is well-founded if and only if it is well-ordered with respect to the Kleene Brouwer order (See \cite[Section V.1]{Simp}).
\end{defin}

\begin{cor}\label{Pioneone} (See \cite[Lemma VI.1.1]{Simp}) over $\textbf{RCA}_0$ the following are equivalent:
\begin{enumerate}[label={(\arabic*)}]
\item $\Pi^1_1$ comprehension.
\item For any sequence of trees $(T_i)_{i\in\mathbb{N}}$ the set $\{i\in\mathbb{N}: T_i \text{ is well-founded}\}$ exists.
\item For any sequence of linear orders $(L_i,<_i)_{i\in \mathbb{N}}$ the set $\{i\in\mathbb{N}: (L_i,<_i) \text{ is a well order}\}$ exists.
\end{enumerate}
\end{cor}

\section{Countable second countable spaces}
In this section we will go over some basic definitions and facts about CSC spaces. We will also go over how countable metric spaces may be treated as CSC spaces over $\mathbf{RCA}_0$.
\begin{defin} A topological space is said to be second countable if it has a countable basis of open sets. We say that a topological space is first countable if the neighborhood filter of every point has a countable cofinal sequence. It is easy to check that for a countable space, being first countable is equivalent to being second countable.
\end{defin}

\begin{obs} One may wonder if all countable spaces are second countable. This is not the case. A simple procedure to construct countable but not second countable spaces is to take any downwards directed family $\mathcal{F}\subseteq\mathcal{P}(\mathbb{N}_{>0})$ that does not have a countable cofinal sequence, and define on $\mathbb{N}$ the topology in which every $n\in \mathbb{N}_{>0}$ is isolated and the basic open neighborhoods of $0$ are the ones of the form $\{0\}\cup S$ where $S\in\mathcal{F}$.\\
\\
To find such a directed subset of $\mathcal{P}(\mathbb{N}_{>0})$ one can use the fact that there exists a family $\mathcal{B}$ of size $2^{\aleph_0}$ of infinite subsets of $\mathbb{N}_{>0}$ such that the intersection of any two members is finite. Let $\mathcal{F}$ be the family of all finite intersections of the complements of elements in $\mathcal{B}$. We have that $\mathcal{F}$ is downwards directed, but it does not have a cofinal sequence.
\end{obs}

\begin{defin}\textbf{(Dorais \cite[Definition 2.2]{Dorais})} In the context of second order arithmetic, we have that a Countable Second Countable (CSC) space  is a tuple $(X,(U_i)_{i\in \mathbb{N}},k)$ where $X$ is a subset of $\mathbb{N}$, for each $i \in \mathbb{N}$ $U_i\subseteq X$ and $k: X\times\mathbb{N}\times\mathbb{N}\rightarrow \mathbb{N}$ is a partial function such that
\begin{equation*}(x\in U_i\cap U_j) \rightarrow ( x\in U_{k(x,i,j)}\subseteq U_i\cap U_j)
\end{equation*}
A set $A\subseteq X$ is open if for all $x\in A\;\exists i\in \mathbb{N}\,(x\in U_i\subseteq A)$, and a set $C\subseteq X$ is closed if it is the complement of an open set.
\end{defin}

\begin{defin} \textbf{(Dorais \cite[Definition 2.4]{Dorais})}
Given a CSC space $(X,(U_i)_{i\in\mathbb{N}},k)$, an open code is a partial function $f:\mathbb{N}\rightarrow \mathbb{N}$. The open collection coded by $f$ is the collection $A=\{x\in X:\exists n\in\mathbb{N}\, x\in U_{f(n)}\}$.
Similarly, a closed code is a partial function $g:\mathbb{N}\rightarrow \mathbb{N}$. The closed collection coded by $g$ is the collection $C=\{x\in X:\forall n\in\mathbb{N}\, x\notin U_{g(n)}\}$.\\
\\
In general, $\textbf{RCA}_0$ may not be sufficient to show that a coded open or closed  collection exists as a set. Similarly, we have that $\mathbf{RCA}_0$ is not sufficient to prove that every open or closed set has a code. A coded open collection that exists as a set is called an effectively open set, and a coded closed collection which exists as a set is called an effectively closed set. We leave it to the reader to show that every coded open collection exists as a set and every open set is coded are both equivalent to arithmetic comprehension over $\mathbf{RCA}_0$.
\end{defin}

\begin{obs}
The coded open collections will be closed under finite intersections and countable unions. Given $A_0$ and $A_1$ open collections with codes $f_0$ and $f_1$, the partial function $h:\mathbb{N}\times\mathbb{N}\times X\rightarrow \mathbb{N}$ given by $h(x,i,j)=k(x,f_0(i),f_1(j))$ will code $A_0\cap A_1$. Given a sequence of coded open collections $(A_i)_{i\in\mathbb{N}}$ with codes $(f_i)_{i\in\mathbb{N}}$ the partial function $F:\mathbb{N}\times \mathbb{N}\rightarrow \mathbb{N}$ given by $F(i,n)=f_i(n)$ is the code for the open collection $\bigcup_{i\in\mathbb{N}} A_i$.  Similarly, coded closed collections are closed under finite unions and countable intersections.
\end{obs}

\begin{defin} Let $(X,(U_i)_{i\in\mathbb{N}},k)$ be a CSC space, we say that a sequence of open sets $(A_n)_{n\in\mathbb{N}}$ is a sequence of uniformly effectively open sets if there exists a sequence $(f_n)_{n\in\mathbb{N}}$, such that for all $n\in\mathbb{N}$ $f_n$ is an open code for $A_n$. Similarly, $(C_n)_{n\in\mathbb{N}}$ is a sequence of uniformly effectively closed sets if there exists a sequence $(g_n)_{n\in\mathbb{N}}$ such that for all $n\in\mathbb{N}$ $g_n$ is a closed code for $C_n$.
\end{defin}

\begin{defin}\textbf{(Dorais \cite[Definition 2.5]{Dorais})} \label{effhom} For a pair $(X,(U_i)_{i\in\mathbb{N}},k)$ and $(Y,(V_i)_{i\in\mathbb{N}},k^*)$ of CSC spaces, we say that a function $f:X\rightarrow Y$ is effectively continuous if there exists a function $v:X\times\mathbb{N}\rightarrow \mathbb{N}$ such that for every $x\in X$ and $j\in \mathbb{N}$ if $f(x)\in V_j$ then $x\in U_{v(x,j)}$ and $f(U_{v(x,j)})\subseteq f(V_j)$. We say that the function $v$ verifies that $f$ is continuous. We say that $f:X\rightarrow Y$ is effectively open if there exists a function $v:X\times\mathbb{N}\rightarrow \mathbb{N}$ such that for every $x\in X$ and $j\in \mathbb{N}$ such that $x\in U_j$ then $f(x)\in V_{v(x,j)}$ and $V_{v(x,j)}\subseteq f(U_j)$. If $f$ is effectively continuous and has an effectively continuous inverse, or equivalently, $f$ is an effectively continuous and effectively open bijection, then we say $f$ is an effective homeomorphism. We will write $X\cong Y$ to say that $X$ is effectively homeomorphic to $Y$. We say that $f:X\rightarrow Y$ is an effective embedding if $f$ is effectively continuous and there is a function $v:X\times \mathbb{N}\rightarrow \mathbb{N}$ such that for all $x\in X$ and $i\in \mathbb{N}$ such that $x\in U_i$ then $x\in f^{-1}(V_{v(x,i)})\subseteq U_i$. We will say that the function $v$ verifies that $f$ is effectively open with its range. We note, however, that in general, over $\textbf{RCA}_0$, the range of an effective embedding may not exist as a set.
\end{defin}

\begin{obs} Since over $\textbf{ACA}_0$ every CSC space has a function $k$ we will often omit it from the definition of CSC space. Similarly, we won't distinguish homeomorphisms from effective homeomorphisms over $\textbf{ACA}_0$.
\end{obs}
\begin{defin}(See \cite[Definition II.4.4]{Simp}) A real number is a sequence $(q_k)_{k\in\mathbb{N}}$ of rational numbers such that $\forall k\;\forall i\;(|q_{k}-q_{k+i}|\leq 2^{-k})$. We denote the class of real numbers by $\mathbb{R}$.  Let  $(p_k)_{k\in\mathbb{N}}$ and $(q_k)_{k\in\mathbb{N}}$ be real numbers, we define:
\begin{enumerate}[label={(\arabic*)}]
\item $(p_k)_{k\in\mathbb{N}}=_{\mathbb{R}}(q_k)_{k\in\mathbb{N}}$ if and only if $\forall k\;|p_k-q_k|\leq 2^{-k+1}$.
\item $(p_k)_{k\in\mathbb{N}}\leq_{\mathbb{R}}(q_k)_{k\in\mathbb{N}}$ if and only if $\forall k\;p_k-q_k\leq 2^{-k+1}$.
\item $(p_k)_{k\in\mathbb{N}}<_{\mathbb{R}}(q_k)_{k\in\mathbb{N}}$ if and only if $\exists k\;q_k-p_k>2^{k+1}$.
\end{enumerate} 
Addition on the reals is defined pointwise. Given a rational $q$ and a real $r=(r_n)_{n\in\mathbb{N}}$ we define the product $q\cdot r$ to be the sequence $(q\cdot r_{n+ k})_{n\in\mathbb{N}}$ where $k=\min\{j\in\mathbb{N}:2^j\geq q\}$.
\end{defin}
\noindent
We would like to be able to give a CSC structure to countable metric spaces. We will not be considering metric spaces the way they are usually defined in the context of reverse math (see \cite[Defintion I.8.1]{Simp} or \cite[Definition 10.1.5]{Rmprp}) which is that of coded Polish spaces. Since we are working with CSC spaces, for our purposes we will only be considering countable metric spaces which will be a pair $(X,d)$ where $X\subseteq \mathbb{N}$ is the collection of points and $d:X\times X\rightarrow \mathbb{R}$ is a metric.\\
\\
It is not immediate how a CSC structure may be given effectively to a countable metric space. In general, for a countable metric space $(X,d)$, the open balls $B(x,r)$ are $\Sigma^0_1$ definable but may not be $\Pi^0_1$ definable relative $(X,d)$. So $\Delta^0_1$ comprehension does not suffice to show that the open balls exist as sets. One solution is to simply consider CSC spaces with a weak basis as is done in \cite[Definition 10.8.2]{Rmprp} in which the elements of the basis are the images of partial functions rather than sets. We present another approach to giving countable metric spaces a CSC structure.\\
\\
We note that if the ball is clopen, that is, $B(x,r)=\overline{B(x,r)}$, then
\begin{equation*}
B(x,r)=\{y\in x: d(x,y)<_\mathbb{R} r\}=\{y\in x: d(x,y)\leq_\mathbb{R} r\}
\end{equation*}
Since strict inequality between reals is $\Sigma^0_1$ and weak inequality is $\Pi^0_1$, we have that $B(x,r)$ is $\Delta^0_1$ definable. Working over $\textbf{RCA}_0$, any clopen ball exists by $\Delta^0_1$ comprehension. So, the problem of finding a CSC structure for a countable metric space is reduced to finding a basis of clopen balls.

\begin{prop}\label{metstruc} $\textbf{RCA}_0$ proves that for any countable metric space $(X,d)$ there exists a real $a\in \mathbb{R}$ such that
\begin{equation*}
\forall q\in\mathbb{Q}_{>0}\,\forall x,y\in X\, (d(x,y)\neq q\cdot a )
\end{equation*} 
and that there is a function $k$ such that $(X, (B(x,q\cdot a) )_{q\in \mathbb{Q}_{>0},x\in X},k)$ is a CSC space.\\
\\
\begin{proof} 
Informally, we will construct a real $a$ that is the diagonal of the reals $(q\cdot d(x,y))_{q\in\mathbb{Q}_{>0},x,y\in X}$. Let $(t_n)_{n\in \mathbb{N}}$ be an enumeration of all triples $(x,y,q)\in X\times X\times \mathbb{Q}_{>0}$. We define a real $a=(a_i)_{i\in\mathbb{N}}$ recursively as follows. Set $a_0=\frac{1}{2}$. At step $i+1$ let $t_i=(x,y,q)$ and let $q\cdot d(x,y)=(b_j)_{j\in\mathbb{N}}$. If $|b_{2i+2}-a_{2i}|<2^{-2i-1}$ then:
\begin{enumerate}[label={(\arabic*)}]
\item $b_{2i+2}\geq a_{2i}$, in which case we define $a_{2i+1}=a_{2i+2}=b_{2i+2}- 2^{-2i-1}$.
\item $b_{2i+2}< a_{2i}$, in which case we define $a_{2i+1}=a_{2i+2}=b_{2i+2}+ 2^{-2i-1}$.
\end{enumerate}
If instead $|b_{i+1}-a_i|\geq 2^{-2i-1}$ then we define $a_{i+1}=a_i$. By induction, the triangle inequality, and using the fact that
\begin{equation*} \forall i \in\mathbb{N}\;(|a_{2i}-a_{2i+1}|=|a_{2i}-a_{2i+1}|\leq 2^{-2i-1})
\end{equation*} we have that $(a_i)_{i\in\mathbb{N}}$ is a real number. Since $a=(a_{i})_{i\in\mathbb{N}}$ is recursive relative to $(X,d)$ it exists by $\Delta^0_1$ comprehension. By construction of $a$, we have that for all $x,y\in X$, $q\in \mathbb{Q}$ $a\neq_{\mathbb{R}} q\cdot d(x,y)$. In particular we have that for all $q\in\mathbb{Q}$ and all $x\in X$ the ball $B(x,q\cdot a)$ is clopen. Since the clopen balls are uniformly $\Delta^0_1$ definable relative to $(X,d)$, we have that the collection $(B(x,q\cdot a))_{q\in \mathbb{Q}_{>0},x\in X}$ exists by $\Delta^0_1$ comprehension. Let $(q_i)_{i\in\mathbb{N}}$ be an enumeration of $\mathbb{Q}_{>0}$ and let $k:X\times (\mathbb{Q}_{>0} \times X)^2\rightarrow \mathbb{Q}_{>0} \times X$ be the function given by
\begin{equation*}
k(x,(y,p),(z,r)) =(x,q_j)
\end{equation*}
where
\begin{equation*}
j=\min \{ i\in \mathbb{N} : p\neq q_i\neq r\wedge d(x,y) <_\mathbb{R} (p-q_i)\cdot a \wedge d(x,z) <_\mathbb{R} (r-q_i)\cdot a\}
\end{equation*}
or rather $j$ is the least such that $q_j\cdot a <_\mathbb{R} p\cdot a-d(x,y)$ and $q_j\cdot a <_\mathbb{R} r\cdot a-d(x,z)$. By construction of $a$, we also have that
\begin{equation*}
j=\min \{ i\in \mathbb{N} : p\neq q_i\neq r\wedge d(x,y) \leq_\mathbb{R} (p-q_i)\cdot a \wedge d(x,z) \leq_\mathbb{R} (r-q_i)\cdot a\}
\end{equation*}
So $j$ can be searched effectively. By the triangle inequality, we have
\begin{equation*}
B(x,q_j\cdot a)\subseteq B(y,p\cdot a)\cap B(z,r\cdot a)
\end{equation*}
By construction of $a$, we have that $k$ is recursive and so it exists by $\Delta^0_1$ comprehension.
\end{proof}
\end{prop}
\noindent 
We note that we did not need to use the fact that $d$ separates points in the previous proof.
The previous proposition tells us that any countable metric space has a basis of clopen sets. We will see later that any CSC space that distinguishes points (is $T_0$) and  has a basis of clopen sets is homeomorphic to a countable metric space.

\begin{defin} Let $(X,d)$ be a countable metric space, and $a\in \mathbb{R}_{>0}$ such that
\begin{equation} \label{metproperty}
\forall q\in\mathbb{Q}_{>0}\,\forall x,y\in X\, (d(x,y)\neq q\cdot a )
\end{equation} 
The CSC structure of $(X,d)$ is $(X, (B(x,q\cdot a) )_{q\in \mathbb{Q}_{>0},x\in X},k)$. It is straightforward to show that for any $a,b\in \mathbb{R}_{>0}$ satisfying \eqref{metproperty} the spaces $(X, (B(x,q\cdot a) )_{q\in \mathbb{Q}_{>0},x\in X},k)$ and $(X, (B(x,q\cdot b) )_{q\in \mathbb{Q}_{>0},x\in X},k)$ will be effectively homeomorphic. In particular, the CSC structure on $(X,d)$ is unique up to effective homeomorphism. We say that a CSC space is metrizable if it is effectively homeomorphic to the CSC structure of a countable metric space.
\end{defin}

\begin{defin}
Let $((X^j,(U_i^j)_{i\in\mathbb{N}},k^j))_{j\in\mathbb{N}}$ be a sequence of CSC spaces we define the topological disjoint union of the sequence of spaces as being the space
\begin{equation*}
\left(\coprod_{j\in \mathbb{N}} X^j\,,\,( \{j\}\times U_i^j )_{i,j\in\mathbb{N}}\;,\, k^* \right)
\end{equation*}
Where $\coprod_{j\in \mathbb{N}} X^j=\{(j,x):x\in X^j\}$ and $k^*((j,x),(j,a),(j,b))=(j,k^j(x,a,b))$. We have that for every $j\in\mathbb{N}$ $\{j\}\times X^j $ is a clopen set and the map from $X^j$ to $\coprod_{j\in\mathbb{N}} X^j$ given by $x\mapsto (j,x)$ is an effective embedding. Furthermore, we observe that the disjoint union of the spaces is computable relative to the sequence of CSC spaces.
\end{defin}

\begin{defin}
Given a linear order $(L,<_L)$ the order topology on $L$ is the space $(L,(\, ]a,b[\,)_{a,b\in L\cup\{+\infty,-\infty\}},k)$ where the function $k$ is given by
\begin{equation*}
k(z,(a,b),(x,y)\,)=\;(\max\{a,x\},\min\{b,y\})
\end{equation*}
By the standard indexing of the open intervals of $L$ we mean the sequence $(\, ]a,b[\,)_{a,b\in L\cup\{+\infty,-\infty\}}$. We will assume that every linear order has the order topology unless specified. We define the upper limit topology on $(L,<_L)$ to be the order topology where sets of the form $]a,b]$, with $a,b\in L$,  are added to the basis.\\
\\
Given a sequence of linear orders $(L_i,<_i)_{i\in\mathbb{N}}$ we define the sum order $\sum_{i\in\mathbb{N}}(L_i,<_i)$, or simply $\sum_{i\in\mathbb{N}}L_i$, to be the order with field $\bigcup_{i\in\mathbb{N}}\{i\}\times L_i$ with the lexicographic order, that is $(i,k)<_{lex}(j,l)$ if $i<j$ or,  $i=j$ and $k<_i l$. Given a linear order $(L,<_L)$ we denote its dual by $L^*$ or $(L,<^*_L)$, where for all $i,j\in L$ we have $i<^*_L j\leftrightarrow j<_L i$.\\
\\
A CSC space  $X$ is orderable if it is effectively homeomorphic to the order topology of a linear order and we say that $X$ is well-orderable if it is effectively homeomorphic to the order topology of a well-order.
\end{defin}

\begin{obs} If $L$ is well-ordered, its order topology coincides with its upper limit topology since $\mathopen]a,b\mathclose]=\mathopen]a, b+1\mathclose[$ where $b+1$ is either the successor of $b$ or $+\infty$ if $b$ is maximal. However, we will see later that over $\textbf{RCA}_0$, the two topologies are not, in general, effectively homeomorphic.
\end{obs}

\begin{obs} Let $(L,<_L)$ be a linear order $S\subseteq L$ be a set that is downwards closed with respect to $<_L$ then $S$ is open in the upper limit topology since $S=\bigcup_{s\in S}\;\mathopen]-\infty,s\mathclose]$.
\end{obs}

\section{Separation axioms}
We will introduce some of the basic topological separation axioms. Similar to how the function $k$ in the definition of CSC space witnesses the intersection property of the basis and the function $v$ witnesses the pointwise continuity of an effectively continuous function, we will also want to consider spaces which have a function that witnesses a given separation axiom. We will also show that assuming this additional structure exists for all CSC spaces satisfying the appropriate separation axioms requires arithmetic comprehension.
\begin{defin} A topological space $X$ is said to be Kolmogorov or $T_0$ if for every $x,y\in X$ there exists an open set $U$ such that $x\in U\leftrightarrow y\notin U$.
\end{defin}

\begin{defin} The Kolmogorov quotient of a CSC space $(X,(U_i)_{i\in\mathbb{N}},k)$ is the quotient space $X/\sim_K$ where the equivalence relation $\sim_K$ is given by:
\begin{equation*}
x\sim_K y\leftrightarrow (\forall i\in \mathbb{N}\; x\in U_i\leftrightarrow y\in U_i)
\end{equation*}
We may identify the Kolmogorov quotient with the subspace of $X$ given by the least elements of the $\sim_K$-equivalence classes. The Kolmogorov quotient is always a $T_0$ space. In general, the existence of the Kolmogorov quotient or the equivalence relation $\sim_K$ requires arithmetic comprehension. To see this, let $A\subseteq \mathbb{N}$ be a set and consider the space $(\mathbb{N},(U_i)_{i\in \mathbb{N}},k)$ where:
\begin{equation*}
U_{(n,t)}=\begin{cases} \{2n,2n+1\} \text{ if } \neg \Phi^A_e(e){\downarrow}_{\leq t}\\
\{2n\} \text{ if } \Phi^A_e(e){\downarrow}_{\leq t}
\end{cases}
\end{equation*} 
the Kolmogorov quotient of this CSC space will compute the jump of $A$.
\end{defin}

\begin{defin} A CSC space is said to be $T_1$ if every singleton is closed.
\end{defin}

\begin{defin} A CSC space $(X,(U_i )_{i\in\mathbb{N}},k)$ is said to be Hausdorff or $T_2$ if for every distinct $x,y\in X$
 there exists $i,j\in \mathbb{N}$ such that:
\begin{equation*}
x\in U_i \quad \quad y\in U_j \quad\quad U_i\cap U_j=\emptyset
\end{equation*}
\end{defin}

\begin{defin} \textbf{(Dorais \cite[Definition 6.1]{Dorais})}
We say that a space is effectively $T_2$ or $eT_2$ if there exists a pair of functions  $H_0,H_1:X^2 \rightarrow \mathbb{N}\times \mathbb{N}$ such that
\begin{equation*}
(x\in U_{H_0(x,y)})\wedge ( y\in U_{H_1(x,y)}) \wedge (U_{H_0(x,y)}\cap U_{H_1(x,y)}=\emptyset)
\end{equation*}
We observe that the functions $H_0$ and $H_1$ are arithmetically definable. So over $\textbf{ACA}_0$ every $T_2$ CSC space is $eT_2$. 
\end{defin}

\begin{prop} \label{embT2} If $(X,(U_i)_{i\in\mathbb{N}},k)$ has an effectively continuous injection to an $eT_2$ CSC space then $X$ is $eT_2$.\\
\\
\begin{proof} Let $f:(X,(U_i)_{i\in\mathbb{N}},k)\rightarrow (Y,(V_i)_{i\in\mathbb{N}},k')$ be an effectively continuous injection and let $H^Y_0$ and $H^Y_1$ witness that $(Y,(V_i)_{i\in\mathbb{N}},k')$ is $eT_2$. Let $v:X\times \mathbb{N}\rightarrow \mathbb{N}$ verify that $f$ is continuous, meaning that for all $x\in X$ and $i\in \mathbb{N}$ such that $f(x)\in V_i$ we have that
\begin{equation*}
x\in U_{v(x,i)}\subseteq f^{-1}(V_i)
\end{equation*}
For all $x,y\in X$ we define
\begin{equation*} H^X_0(x,y)=v(x,(H^Y_0(f(x),f(y)))
\end{equation*}
\begin{equation*} H^X_1(x,y)=v(y,(H^Y_1(f(x),f(y))).
\end{equation*}
The functions  $H^X_0$ and $H^X_1$ exist by $\Delta^0_1$ comprehension and witness that $X$ is $eT_2$.
\end{proof}
\end{prop}

\begin{thm}\label{eT2>aca} \textbf{(Dorais \cite[Example 7.4]{Dorais})} If there exists a function $f:\mathbb{N}\rightarrow \mathbb{N}$ whose range does not exist as a set then there is a $T_2$ CSC space which is not $eT_2$. Equivalently, over $\textbf{RCA}_0$, arithmetic comprehension is equivalent to every $T_2$ CSC space being $eT_2$.\\
\\
\begin{proof} Assume there is an injective function $f:\mathbb{N}\rightarrow \mathbb{N}$ whose range does not exist as a set. Let $<_f$ be the order on $\mathbb{N}\cup\{\infty\}$ given by:
\begin{itemize}
\item $\forall n\in\mathbb{N}\; n<_f\infty$.
\item $\forall n,m\in \mathbb{N}\; n<_f m\leftrightarrow f(n)<f(m)$.
\end{itemize}
Seeking a contradiction, assume that $(\mathbb{N}\cup\{\infty\},<_f)$ with the order topology is $eT_2$ and let $H_0$ and $H_1$ be the functions that witness it.  We have that $n$ is the $<_f$ successor of $m$ if and only if $m<_f n$ and $H_0(m,n)=(y,n)$ for some $y<_f n$ and $H_1(m,n)=(m,x)$ for some $x>_f m$. Since the image of $f$ is unbounded, every $n\in \mathbb{N}$ has an $<_f$ successor and the successor function exists since it is computable relative to $H_0$ and $H_1$. We can therefore define recursively $g$ where $g(0)$ the $<_f$-minimal element and $g(n+1)$ is the $<_f$-successor of $g(n)$. So $f\circ g$ is strictly increasing and has the same image as $f$, which implies that $\text{rng}(f)$ exists, contradicting our initial assumption.
\end{proof}
\end{thm}

\begin{obs} The space constructed in the proof above is a well-order with the order topology. While the upper limit topology and the order topology coincide on a well-order, they may not be effectively homeomorphic. This is because the upper limit topology is always $eT_2$, and arithmetic comprehension is equivalent to the order topology of any well-order is $eT_2$. However, it is straightforward to show that over $\textbf{RCA}_0$, the identity on a well-order is an effective homeomorphism between the upper limit topology and the order topology if and only if the successor function for the order exists.
\end{obs}

\begin{prop}\label{wellmess} Over $\textbf{RCA}_0$ the following are equivalent:
\begin{enumerate}[label={(\arabic*)}]
\item Arithmetic comprehension.
\item Every well-order with the upper limit topology is homeomorphic to the order topology of some linear order.
\end{enumerate}
\begin{proof} $(1\rightarrow 2)$ follows from the fact that the successor function on a well-order is arithmetically definable.\\
\\
$(2\rightarrow 1)$ Let $f:\mathbb{N}\rightarrow \mathbb{N}$ be an injection, $(X,<_f)$ be as in the proof of Theorem~\hyperref[eT2>aca]{\ref{eT2>aca}} and let $\infty_X$ denote the maximal element of $X$. We show that $\text{rng}(f)$ exists, which by \hyperref[ACAnot]{\ref{ACAnot}} implies arithmetic comprehension. By assumption, $X$ with the upper limit topology is effectively homeomorphic to some linear order $(L,<_L)$ with the order topology. Since $X$ is $eT_2$, we also have that $L$, with the order topology, is $eT_2$. As in the proof of Theorem~\hyperref[eT2>aca]{\ref{eT2>aca}}, we have that the successor and predecessor functions for $L$ will exist.
Since $X$ has one limit point and doesn't contain an infinite closed discrete subspace, we have by induction that $L$ will be order isomorphic to $\mathbb{N}+n$ or its dual for some $n\in \mathbb{N}$ or to $\mathbb{N}+1+\mathbb{N}^*$. The order topology of all these orders is effectively homeomorphic to that of $\mathbb{N}+1$ so without loss of generality we may assume $L=\mathbb{N}+1$. Let $g:X\rightarrow \mathbb{N}+1$ be an effective homeomorphism and let $v$ the function which witnesses that $g$ is effectively open. Let $\infty_\mathbb{N}$ denote the maximal element of $\mathbb{N}+1$, we have that
\begin{equation*}
\forall n\in\mathbb{N}\;(\mathopen]v(\infty_X,n),\infty_\mathbb{N}\mathclose]_{<_\mathbb{N}}\subseteq g(\, \mathopen] n,\infty_X \mathclose]_{<_f}))
\end{equation*}
Which implies that $\forall n\; \forall m\;(m\leq_f n\rightarrow g^{-1}(m)\leq_\mathbb{N} v(\infty_X,n))$. So in particular we have that $n$ is the $<_f$ successor of $m$ if and only if
\begin{equation*}
\forall y\leq_\mathbb{N} v(\infty_X,n)\;(g^{-1}(y)\leq_\mathbb{N} m\vee g^{-1}(y)\geq_\mathbb{N} n)
\end{equation*}
The $<_f$ successor function is $\Delta^0_1$ definable, so it exists by $\Delta^0_1$ comprehension. Using the successor function as in the proof of Theorem~\hyperref[eT2>aca]{\ref{eT2>aca}}, we get that the range of $f$ exists.
\end{proof}
\end{prop}

\begin{defin}
A CSC space $(X,(U_i)_{i\in\mathbb{N}},k)$  is said to be regular if for every point $x\in X$ and any $i\in\mathbb{N}$ such that  $x\in U_i$ there are open sets $V$ and $W$ such that
\begin{equation*}
x\in V\subseteq X\setminus W\subseteq U_i
\end{equation*}
\end{defin}

\begin{defin}
A CSC space that is $T_0$ and regular is said to be $T_3$. So over $\textbf{ACA}_0$, the study of regular spaces is reduced to the study of $T_3$ spaces since we may always restrict ourselves to the Kolmogorov quotient.
\end{defin}

\begin{obs} All $T_3$ spaces are $T_2$.
\end{obs}

\begin{defin}\label{effregdor}\textbf{(Dorais \cite[Definition 6.4]{Dorais})} A CSC space $(X,(U_i)_{i\in\mathbb{N}},k)$ is said to be effectively regular if for every effectively closed set $C$ of $X$ and any $x\notin C$ there exist open collections $U_0$ and $U_1$ for open sets such that $x\in U_0 \subseteq X\setminus U_1\subseteq X\setminus C$.
\end{defin}
The issue with Dorais' notion of effective regularity is that it lacks uniformity. There may not be a map which given a point $x$ and a code for an effective closed  set  $C$ gives codes for the open collections separating $x$ and $C$. We introduce a stronger notion of effective regularity which is equivalent to Schr{\"o}der's notion of effective regularity for represented spaces \cite{Schrod} and Louveau's notion of recursive regularity for basic spaces \cite{Louveau}. Such notion of regularity was shown by both authors to be the necessary and sufficient condition for metrizability. We show that over $\mathbf{RCA}_0$ the same is true for CSC spaces.

\begin{defin} We say that a CSC space $(X,(U_i)_{i\in\mathbb{N}},k)$ is uniformly regular or uniformly effectively regular if there is a pair of functions $R_0$ and $R_1$ such that for all $x,y\in X$ and $i\in\mathbb{N}$ such that $x\in U_i$ and $y\notin U_i$ the following holds
\begin{equation*}
x\in U_{R_0(x,i)}\subseteq U_i
\end{equation*}
\begin{equation*}
y\in U_{R_1(x,i,y)}\subseteq X\setminus U_{R_0(x,i)}
\end{equation*}
We will say that a CSC space is uniformly $T_3$ or $uT_3$ if it is $T_0$ and uniformly regular. Intuitively, $R_0(x,i)$ is the index of a neighborhood of $x$ whose closure is contained in $U_i$ and the sequence $(R_1(x,i,y))_{y\notin U_i}$ codes an open collection containing $X\setminus U_i$ that is disjoint from $U_{R_0(x,i)}$. That is, we have
\begin{equation*}
x\in U_{R_0(x,j)}\subseteq  X\setminus \left( \bigcup_{y\notin U_{j}}U_{R_1(x,j,y)}\right) \subseteq X\setminus U_i
\end{equation*}
\end{defin}

\begin{obs} Subspaces of regular CSC spaces are regular and subspaces of uniformly effectively regular CSC space are uniformly effectively regular.
\end{obs}

\begin{prop}\label{Regstab} $\textbf{RCA}_0$ proves that for any pair of CSC spaces $(X,(U_i)_{i\in \mathbb{N}},k)$ and $(Y,(V_i)_{i\in\mathbb{N}},k')$ such that $X$ effectively embeds into $Y$, if $Y$ is uniformly regular then $X$ is uniformly regular .\\
\\
\begin{proof}
Let $f:X\rightarrow Y$ be an effective embedding and $v_0$ verify that $f$ is continuous and $v_1$ verify that $f$ is open with respect to its image (See Definition~\hyperref[effhom]{\ref{effhom}}). Let $R^Y_0, R^Y_1$ be the functions given by uniform effective regularity of $Y$. We define
\begin{equation*}
R^X_0(x,i)=v_0(x,R^Y_0(f(x),v_1(f(x),i)))
\end{equation*} 
We have that
\begin{equation*}
x\in U_{R^X_0(x,i)}\subseteq f^{-1}(V_{R^Y_0(f(x),v_1(f(x),i))})\subseteq U_i
\end{equation*}
We define for $y\notin U_i$
\begin{equation*}
R^X_1(x,i,y)=v_0(y,R^Y_1(f(x),v_1(f(x),i),f(y)))
\end{equation*}
Since $y\notin U_i$ we have $f(y)\notin V_{v_1(f(x),i)}$ so $R^Y_1(f(x),v(f(x),i),f(y))$ is well defined. By definition of $R^Y_1$ we have
\begin{equation*}
f(y)\in V_{R^Y_1(f(x),v_1(f(x),i),f(y))}\subseteq Y\setminus V_{R^Y_0(f(x),v_1(f(x),i))}\subseteq Y\setminus f(U_{R^X_0(x,i)}).
\end{equation*}
In particular
\begin{equation*}
y\in U_{R^X_1(x,i,y)}\subseteq f^{-1}(V_{R^Y_1(f(x),v_1(f(x),i),f(y))})\subseteq f^{-1}(Y\setminus f(U_{R^X_0(x,i)}))=X\setminus U_{R^X_0(x,i)}
\end{equation*}
This means that the function $R^X_0$ and $R^X_1$ have the desired properties, and so $X$ is uniformly effectively regular.
\end{proof}
\end{prop}

\begin{prop} $\textbf{RCA}_0$ proves that every $uT_3$ CSC space is $eT_2$.\\
\\
\begin{proof}\label{label} Let $(X,(U_i)_{i\in\mathbb{N}},k)$ be a $uT_3$ CSC space and let $x,y\in X$ be distinct points. Since $X$ is $T_0$, let $i$ be the least number such that $x\in U_i\leftrightarrow y\notin U_i$.  If $x\in U_i$ then we let $H_0(x,y)=R_0(x,i)$ and $H_1(x,y)=R_1(x,i,y)$ otherwise if $y\in  U_i$ we set $H_0(x,y)=R_1(x,i,y)$ and $H_1(x,y)=R_0(x,i)$.
\end{proof}
\end{prop}

\begin{cor} Over $\textbf{RCA}_0$ arithmetic comprehension is equivalent to every regular CSC space is uniformly regular.\\
\\
\begin{proof} Follows from Theorem~\hyperref[eT2>aca]{\ref{eT2>aca}} and Proposition~\hyperref[label]{\ref{label}}.
\end{proof}
\end{cor}
\noindent
We would like to show that over $\textbf{RCA}_0$ $uT_3$ is in general a stronger notion than $eT_2$ for $T_3$ CSC spaces.  However, the next proposition will show us that we cannot use linear orders to separate these properties.
\begin{prop}\label{pain} $\textbf{RCA}_0$ proves that every $eT_2$ order topology is also $uT_3$.\\
\\
\begin{proof}
Let $(I_e)_{e\in\mathbb{N}}$ be the standard indexing of the open intervals of $(L,<_L)$ and $H_0$, $H_1$ be the functions that witness that the order topology on $L$ is $eT_2$. Given $y<_L x<_L z$ we define
\begin{equation*}
R_0(x,(y,z))=k(x,H_0(x,y),H_0(x,z))
\end{equation*}
and
\begin{equation*}
R_1(w,x,(y,z))=\begin{cases} \text{ index for: } I_{H_1(x,y)}\cup\mathopen]-\infty,y\mathclose] \text{ if } w\leq_L y \\
\text{ index for: } I_{H_1(x,z)}\cup\mathopen[z,+\infty \mathclose[ \text{ if } z\leq_L w
\end{cases}
\end{equation*}
The functions $R_0$ and $R_1$ witness that the order topology on $L$ is $uT_3$.
\end{proof}
\end{prop}

\begin{prop}\label{worst proof} Over $\textbf{RCA}_0$ arithmetic comprehension is equivalent to all regular $eT_2$ CSC spaces being $uT_3$.\\
\\
\begin{proof} We show arithmetic comprehension proves all regular CSC spaces are uniformly regular. Given a regular CSC space $(X,(U_i)_{i\in\mathbb{N}},k)$, for each $(x,j)\in X\times\mathbb{N}$ such that $x\in U_j$ we define
\begin{equation*}
R_0(x,j)=\min\{s\in\mathbb{N}: x\in U_s\subseteq \overline{U_s}\subseteq U_j\}
\end{equation*}
and for all $y\in X\setminus U_j$ we define
\begin{equation*}
R_1(x,j,y)=\min\{s\in \mathbb{N}: y\in U_s\wedge U_s\cap U_{R_0(x,j)}=\emptyset\}
\end{equation*}
The $R_0$ and $R_1$ have the desired properties and are arithmetically definable. So $\textbf{ACA}_0$ proves that every regular CSC space is uniformly regular.\\
\\
We show the converse. Let $A$ be a set. We show that the Turing jump of $A$ exists, which will imply arithmetic comprehension. We define
\begin{equation*}
X^e=\{-\infty\}\cup((\mathbb{Z}\cup\{+\infty\})\times\{t\in\mathbb{N}:\Phi^A_e(e){\downarrow}_{\leq t}\})
\end{equation*}
We define on $X^e$ the following topology
\begin{equation*}
U_0=\{-\infty\}\cup (\mathbb{Z}\times\{t\in\mathbb{N}:\Phi^A_e(e){\downarrow}_{\leq t}\})
\end{equation*}
For all $(m,n)\in (\mathbb{Z}\times\{t\in\mathbb{N}:\Phi^A_e(e){\downarrow}_{\leq t}\})$ we define
\begin{equation*}
U_{4(m,n)+1}=\{(m,n)\}
\end{equation*}
\begin{equation*}
U_{4(m,n)+2}=\{(+\infty,n)\}\cup\{(l,n):l\geq m\}
\end{equation*}
\begin{equation*}
U_{4(m,n)+3}=X^e\setminus\{(l,s):(s=+\infty)\vee(l\geq m\wedge s\leq n)\}
\end{equation*}
\begin{equation*}
U_{4(m,n)+4}=\begin{cases} X^e \quad \quad\quad  \text{ if } \quad \quad \quad \neg\Phi^A_e(e){\downarrow}_{\leq n}\\
\{-\infty\}\cup\{(l,s):l\leq m \wedge \Phi^A_e(e){\downarrow}_{\leq s}\}\quad \text{ if }\quad \Phi^A_e(e){\downarrow}_{\leq n}
\end{cases}
\end{equation*}
We define $H_0$ and $H_1$ on $X^e$ by:
\begin{itemize}
\item $H_0(-\infty,(m,n))=4(m,n)+3\;$ and $\;H_1(-\infty,(m,n))=4(m,n)+1$
\item $H_0(-\infty,(+\infty,n))=4(m,n)+3\;$ and $\;H_1(-\infty,(+\infty,n))=4(m,n)+2$
\item $H_0((+\infty,l),(+\infty,n))=4(0,l)+2\;$ and $\;H_1((+\infty,l),(+\infty,n))=4(0,n)+2$
\item $H_0((+\infty,l),(m,n))=4(m+1,l)+2\;$ and $\;H_1((+\infty,l),(m,n))=4(m,n)+1$
\item $H_0((p,l),(m,n))=4(p,l)+1\;$ and $\;H_1((p,l),(m,n))=4(m,n)+1$
\end{itemize}
Given $a=4(i,j)+r$ and $b=4(p,l)+s$ let $c = \max\{|i|,j, |p|,l\}+1$ and define:
\begin{itemize}
\item $k((m,n),a,b)=4(m,n)+1$
\item $k((\infty,n),a,b)= 4(c,n)+2$
\item $k(-\infty,a,b)=\begin{cases}4(-c,c)+4 \text{ if } \Phi^A_e(e){\downarrow}_{\leq c}\\
4(-c,c)+3 \text{ otherwise}
\end{cases}$
\end{itemize}
It is straightforward to show that the functions $H_0$, $H_1$, and $k$ have the desired properties and that $X^e$ is an $eT_2$ CSC space. To show that $X^e$ is regular, we first observe that every basic open set, except perhaps $U_0$ and $U_{4(m,n)+3}$, is clopen and $-\infty$ is the only non isolated point of $U_0$ and $U_{4(m,n)+3}$. If $\Phi^A_e(e){\uparrow} $ then $X^e=\{-\infty\}$ otherwise let $t$ be the least number such that $\Phi^A_e(e){\downarrow}_{\leq t}$ then the sequence $(U_{4(-n,t)+4})_{n\in\mathbb{N}}$ is a basis of clopen neighborhoods of $-\infty$.
So every point of $X^e$ has a clopen neighborhood, so $X^e$ is regular.
\\
\\
The CSC space $X=\coprod_{e\in \mathbb{N}} X^e$ will be regular and $eT_2$ so by assumption $X$ is also $uT_3$. Let $R_0$ and $R_1$ be functions that witness the uniform regularity of $X$. For $e\in\mathbb{N}$ we have $\Phi^A_e(e){\downarrow}_{\leq n}$ if and only if for all $m\in \mathbb{N}$ $U_{4(m,n)+4}\subseteq U_0$ in the space $X^e$. For all $e$ there is some $m\in \mathbb{N}$ such that $R_0((e,-\infty),(e,0))=(e,4(m,n)+4)$, where by $(e,-\infty)$ we mean $-\infty$ of the space $X^e$, and $(e,0)$ is the index for $U_0$ in $X^e$. So by construction of $X^e$ we have $\Phi^A_e(e){\downarrow}_{\leq n}$ if for some $m\leq R_0((e,-\infty),(e,0))$ we have $R_0((e,-\infty),(e,0))=(e,4(m,n)+4)$. In particular, $A'\leq_T R_0$, which means the Turing jump of $A$ exists by $\Delta^0_1$ comprehension.
\end{proof}
\end{prop}

\begin{obs}\label{scatter metriz>aca} We observe that the space in the previous proof is scattered. So, arithmetic comprehension is equivalent to all $T_3$ scattered spaces being $uT_3$. A CSC space is scattered if every subspace has an isolated point, we will look at $T_3$ scattered CSC space in the next part. We also point out that the space in the previous proof is effectively regular in the sense of Dorais \hyperref[effregdor]{\ref{effregdor}}. This shows that effective regularity is a strictly weaker property than uniform regularity  over $\mathbf{RCA}_0$.
\end{obs}

\noindent
The following technical result will be helpful in the next part.
\begin{thm}\label{theworst} Over $\textbf{RCA}_0$, every scattered CSC space that is $eT_2$ and $T_3$ effectively embeds into a linear order implies arithmetic comprehension.
\\
\\
\begin{proof}
Fix a set $A\subseteq \mathbb{N}$, we show that $A'$ exists. Let $(X,(U_i)_{i\in\mathbb{N}},k)$ be as in the proof of Proposition~\hyperref[worst proof]{\ref{worst proof}}.  Let $(L,<_L)$ be a linear order, $f:X\rightarrow  L$ an effective embedding, and $v_0$ and $v_1$ be functions that witness respectively that $f$ is effectively continuous and effectively open. For each $e$, we define $m_e$ to be the unique number such that
\begin{equation*}
 v_0((-\infty,e),v_1((-\infty,e),0))=4(m_e,n)+r=s_e
\end{equation*}
where $r\in\{3,4\}$ and $n\in \mathbb{N}$, or $s_e=0$. We have that the map $e\mapsto m_e$ is recursive relative to $v_0$ and $ v_1$, so it exists by $\Delta^0_1$ comprehension. Assume that $\Phi^A_e(e){\downarrow}$, we have that $v_1((-\infty,e),0)$ will be the index of some interval $\mathopen]a,b\mathclose[_{<_L}$. By construction we have
 $U_{s_e}\subseteq f^{-1}(\mathopen]a,b\mathclose[_{<_L})\subseteq U_0$. 
We also have that $f^{-1}(\mathopen[a,b\mathclose]_{<_L})$ is 
a closed set contained in $U_0$ plus at most two other points. If $\Phi^A_e(e){\downarrow}$ and either $r=3$ or $s_e=0$, then the closure of $U_{s_e}$ contains infinitely many points which are not in $U_0$. But the closure of $U_{s_e}$ is contained in $U_0$ plus at most two other points. So if $\Phi^A_e(e){\downarrow}$ then $r=4$. By construction of $X$, we have $\Phi^A_e(e){\downarrow}$ and $U_{4(m_e,n)+4}\subseteq U_0$ if and only if $\Phi^A_e(e){\downarrow}_{\leq m_e}$. So we have that $\Phi^A_e(e){\downarrow} \leftrightarrow \Phi^A_e(e){\downarrow}_{\leq m_e}$, which means that $A'$ exists since it is computable from the map $e\mapsto m_e$. Since every set has a Turing jump by Theorem~\hyperref[ACAnot]{\ref{ACAnot}}, we have arithmetic comprehension.
\end{proof}
\end{thm}
\noindent
We will use a construction similar to that in Theorem~\hyperref[worst proof]{\ref{worst proof}} to prove another technical lemma.
\begin{thm}\label{horrible} Over $\textbf{RCA}_0$, for every $T_3$ scattered CSC space there is a continuous bijection from $X$ to the order topology of a well-order implies arithmetic comprehension.\\
\\
\begin{proof} 
Let $A$ be a set, we show that $A'$ exists. Define
\begin{equation*}
X^e=\{-\infty\}\cup((\mathbb{Z}\times\{t\in\mathbb{N}:\Phi^A_e(e){\downarrow}_{\leq t}\})\cup\{(n,+\infty):n\in\mathbb{N}\}
\end{equation*}
with basic open sets:
\begin{equation*}
U_0=\{-\infty\}\cup (\mathbb{Z}\times\{t\in\mathbb{N}:\Phi^A_e(e){\downarrow}_{\leq t}\})
\end{equation*}
and for all $(m,n)\in (\mathbb{Z}\times\{t\in\mathbb{N}:\Phi^A_e(e){\downarrow}_{\leq t}\})$ 
\begin{equation*}
U_{3(m,n)+1}=\{(m,n)\}
\end{equation*}
\begin{equation*}
U_{3(m,n)+2}=\{(+\infty,n)\}\cup\{(l,n):l\geq m\}
\end{equation*}
\begin{equation*}
U_{3(m,n)+3}=\begin{cases} U_0 \quad \quad\quad  \text{ if } \quad \quad \quad \neg\Phi^A_e(e){\downarrow}_{\leq n}\\
\{-\infty\}\cup\{(l,s):l\leq m \wedge \Phi^A_e(e){\downarrow}_{\leq s}\}\quad \text{ if }\quad \Phi^A_e(e){\downarrow}_{\leq n}
\end{cases}
\end{equation*}
Given $a=4(i,j)+r$ and $b=4(p,l)+s$ let $c =\max\{|i|,j, |p|,l\}+1$ and define:
\begin{itemize}
\item $k((m,n),a,b)=3(m,n)+1$
\item $k((\infty,n),a,b)= 3(c,n)+2$
\item $k(-\infty,a,b)=3(-c,c)+3$
\end{itemize}
Take $X=\coprod_{e\in\mathbb{N}} X^e$, we have that it is $T_3$ and scattered. We note that for each $e$ the subspace $X^e\times\{e\}$ is clopen in $X$. By assumption, there is an effectively continuous bijection $f:X\rightarrow W$ where $W$ has the order topology of a well-order $<_W$. Let $(I_i)_{i\in\mathbb{N}}$ be the standard indexing of the basic open intervals of $(W,<_W)$ and $v$ be a function that witnesses that $f$ is effectively continuous. We have that all intervals of a well-order are clopen. Define for each $e\in\mathbb{N}$
\begin{equation*}
i_e=\begin{cases} \text{ index for: } \mathopen]f((e,+\infty,0)),+\infty\mathclose[_{<_W}\quad \text{ if }\quad f((e,+\infty,0))<_W f((e,-\infty)) \\
\text{ index for: } \mathopen]-\infty,f((e,+\infty,0))\mathclose[_{<_W} \quad\text{ if }\quad  f((e,-\infty)) <_W f((e,+\infty,0))
\end{cases}
\end{equation*}
We have that $f^{-1}(I_{i_e})$ will be a clopen neighborhood of $(e,-\infty)$. So, in particular, we have that $v((e,-\infty),i_e)$ will be the index for a basic open set of $(e,-\infty)$ that is contained in a clopen set that does not contain $(e,+\infty,0)$. But by construction of $X$ we have that either $\Phi^A_e(e){\downarrow}_{\leq v((e,-\infty),i_e)}$ or $\Phi^A_e(e){\uparrow}$. So we have that $A'$ is computable relative to $v,(i_e)_{e\in\mathbb{N}},$ and $X$, and so it exists by $\Delta^0_1$ comprehension.
\end{proof}
\end{thm}

\section{Characterizations of dimension zero} 
We already have seen that countable metric spaces have a basis of clopen sets, that is, they are zero dimensional. In this section we will study the relation between uniform regularity and dimension zero. We will also prove that over $\textbf{RCA}_0$, the $uT_3$ CSC spaces are metrizable. Instead of formalizing some form of the Urysohn metrization theorem for uniformly regular CSC spaces in $\mathbf{RCA}_0$ we will instead show that every uniformly regular CSC space is homeomorphic to a closed subset of $\mathbb{Q}$.
\begin{defin} In general topology, a space is said to be zero dimensional or of dimension zero if it has a basis of clopen sets. 
\end{defin}

\begin{defin} Given a CSC space $(X,(U_i)_{i\in \mathbb{N}},k)$ we say that:
\begin{itemize}
\item $X$ has a basis of clopen sets if for all $i\in \mathbb{N}$ the set $U_i$ is also closed. 
\item $X$ has an algebra of clopen sets if there exists $Int:\mathbb{N}^2\rightarrow \mathbb{N}$ such that $U_{Int(i,j)}=U_i\cap U_j$ and $Comp:\mathbb{N}\rightarrow \mathbb{N}$ such that $U_{Comp(i)}=X\setminus U_i$. 
\end{itemize}
\end{defin}

\begin{obs}\label{theobs} For any sequence of sets $(A_i)_{i\in\mathbb{N}}$ the algebra generated by them is $\Delta^0_1$ definable and so are functions $Int$ and $Comp$ that compute respectively intersections and complements. A similar construction is done in \cite[Definition 2.11]{Dorais}. So for any zero dimensional CSC space $(X,(U_i)_{i\in\mathbb{N}},k)$ there is a CSC space with an algebra of clopen sets $(X,(V_i)_{i\in\mathbb{N}},Int,Comp)$ such that $(U_i)_{i\in\mathbb{N}}$ is a subsequence of $(V_i)_{i\in\mathbb{N}}$ and every $V_i$ is open in $(X,(U_i)_{i\in\mathbb{N}},k)$. This means that the identity $X\rightarrow X$ is a homeomorphism. It will not, in general, be an effective homeomorphism. What is needed to ensure the identity is an effective homeomorphism is a function that encodes complements in some way.
\end{obs}

\begin{defin} Let $(X,(U_i)_{i\in\mathbb{N}},k)$ be a CSC space, we say that $X$ is effectively zero dimensional if there is a function $G:X\times\mathbb{N}\rightarrow \mathbb{N}$ such that if $x\in X\setminus U_i$ then $x\in U_{G(x,i)}\subseteq X\setminus U_i$.
\end{defin}

\begin{obs} $\textbf{RCA}_0$ proves that any $T_0$ effectively zero dimensional CSC space is $eT_2$.
\end{obs}

\begin{obs}\label{wozd} Let $(L,<_L)$ be a linear order, the partial function
\begin{equation*}
G(z,(x,y))=\begin{cases} (-\infty,x)\quad \quad \text{ if }\quad \quad z\leq x \\
(y,+\infty)\quad \quad \text{ if }\quad \quad z> y
\end{cases}
\end{equation*}
witnesses that the upper limit topology is effectively zero dimensional.
\end{obs}

\begin{prop}\label{G>Reg} $\textbf{RCA}_0$ proves that every effectively zero dimensional CSC space is uniformly regular.\\
\\
\begin{proof} Let $(X,(U_i)_{i\in\mathbb{N}},k)$ be an effectively zero dimensional CSC space and $G$ a function that witnesses that $X$ is effectively zero dimensional. For all $x\in X$ and $i\in\mathbb{N}$ such that $x\in U_i$ we define $R_0(x,i)=i$ and for all $y\notin U_i$ we define $R_1(x,i,y)=G(y,i)$. It is easy to verify that $R_0$ and $R_1$ witness that $X$ is uniformly regular.
\end{proof}
\end{prop}

\begin{defin} Let $(X,(U_i)_{i\in\mathbb{N}},k)$ be a CSC space, $x\in X$, and a set $I=(i_n)_{n\in\mathbb{N}}$ such that for all $n\in \mathbb{N}$ we have $x\in U_{i_n}$ then we define $K_I$ as:
\begin{enumerate}[label={(\arabic*)}]
\item $K_I(x,0)=i_0$
\item $K_I(x,n+1)=k(x,K(x,n,I),i_{n+1})$
\end{enumerate} 
$\textbf{RCA}_0$ proves the existence of $K_I$ for all CSC spaces and all $I$ since it proves functions are closed under primitive recursion. If $F$ is a finite set such that $\forall i\in F\; x\in U_i$ then we define $K(x,F)=K_F(x,|F|)$. We observe that
\begin{equation*}
x\in U_{K(x,F)}\subseteq \bigcap_{i\in F} U_i
\end{equation*}
This will allow us to find effectively sufficiently small neighborhoods of $x$.
\end{defin}
\noindent
We show that uniformly regular CSC spaces are effectively zero dimensional. To do so, we first show that a strong form of normality holds for uniformly regular CSC spaces. 
\begin{thm}\label{SNORM} Over $\textbf{RCA}_0$, every uniformly regular  CSC space $(X,(U_i)_{i\in\mathbb{N}},k)$ and any two sequences of coded closed collections $(C_{0,e})_{e\in\mathbb{N}}$ and $(C_{1,e})_{e\in\mathbb{N}}$ such that for all $e\in\mathbb{N}$ we have that $C_{0,e}\cap C_{1,e}=\emptyset$, then there exists a sequence of clopen sets $(D_e)_{e\in\mathbb{N}}$ such that for all $e\in\mathbb{N}$ $C_{0,e}\subseteq D_e\subseteq X\setminus C_{1,e}$.\\
\\
\begin{proof} We prove the case for one pair $(C_0,C_1)$ of disjoint closed collections coded by $f_0$ and $f_1$ respectively. The general case will follow from the fact that the construction can be carried out uniformly. To construct $D$ we will define recursively two increasing sequences of closed sets $(C^n_0)_{n\in\mathbb{N}}$ and $ (C^n_1)_{n\in\mathbb{N}}$ such that $C^0_0=C_0$ and $C^0_1=C_1$ and for all $n$ the sets $C^n_0$ and $C^n_1$ are disjoint. We also define two increasing sequences of open sets $(A^n_0)_{n\in\mathbb{N}}$ and $(A^n_1)_{n\in\mathbb{N}}$ such that for all $n\in\mathbb{N}$ and $j\leq 1$ we have $A^n_j\subseteq C^n_j$. We would like that $\bigcup_{n\in\mathbb{N}}A^n_0\cup A^n_1 =X$ and in such case we define $D=\bigcup_{n\in\mathbb{N}}A^n_0$ which by construction will be clopen and $C_0\subseteq D\subseteq X\setminus C_1$. At every step, we will pick an $x\notin A^n_0\cup A^n_1$ and find a small enough neighborhood $U_{l}$ such that for some $j\leq 1$ we have $x\in U_l\subseteq X\setminus C^n_j$. We then set $A^{n+1}_{1-j}=A^n_{1-j}\cup U_{R_0(x,l)}$ and $C^{n+1}_{1-j}=C^n_{1-j}\cup \bigcap_{y\notin U_l}(X\setminus U_{R_1(x,l,y)})$.\\
\\
To be able to construct the sequences $(C^n_0)_{n\in\mathbb{N}},(C^n_1)_{n\in\mathbb{N}},(A^n_0)_{n\in\mathbb{N}},$ and $(A^n_1)_{n\in\mathbb{N}}$ we will want them to be coded as finite objects. We will code $A^n_j$ with a finite set $F^n_j$ of indices such that $A^n_j=\bigcup_{i\in F^n_j} U_i$. We will code $C^n_j$ as a finite set of pairs $E^n_j$ such that
\begin{equation*}C^n_j=\bigcap_{i\in\mathbb{N}}(X
\setminus U_{f_j(i)})\cup \bigcup_{(x,l)\in E^n_j}\left(\bigcap_{y\notin U_l}(X\setminus U_{R_1(x,l,y}))\right)
\end{equation*}
We note that, in general, $C^n_j$ may not exist as a set. \\
\\
Let $X=(x_n)_{n\in\mathbb{N}}$ and let $F^0_0=F^0_1=E^0_0=E^0_1=\emptyset$. 
Given $C^n_0,C^n_1$ and $A^n_0,A^n_1$, if $A^n_0\sqcup A^n_1=X$ then we are done. Otherwise, let $m$ be the least number such that $x_m\in X\setminus (A^n_0\cup A^n_1)$. Let $(i,j)$ be the least pair such that $j\leq 1$ and
\begin{equation*}
(\exists s\leq i\;x_m\in U_{f_j(s)})\wedge \forall (z,l)\in E^n_j\; \exists y\leq i\; x_m \in U_{R_1(z,l,y)}
\end{equation*}
The formula above states that in $i$ steps, we verify $x_m$ is not in $C^n_j$. Since $C^n_0$ and $C^n_1$ are disjoint we have that the search will eventually terminate.
Let $s$ be the least number such that $x_m\in U_{f_j(s)}$ and define
\begin{equation*} I=\{R_1(z,l,y):(z,l)\in E^n_j \wedge y=\min\{ w\leq i :x_m \in U_{R_1(z,l,w)}\}\}\cup\{f_j(s)\}
\end{equation*}
which will be a finite set by $\textbf{B}\Sigma^0_1$ and will be coded. We then define
\begin{equation*}
F^{n+1}_{1-j}=F^n_{1-j}\cup\{R_0(x_m,K(x_m,I))\}\quad \text{ and }\quad F^{n+1}_j=F^n_j
\end{equation*}
and
\begin{equation*}
E^{n+1}_{1-j}=E^n_{1-j}\cup\{(x_m,K(x_m,I))\}\quad \text{ and }\quad E^{n+1}_j=E^n_j
\end{equation*}
\\
We now verify that the sets $A^{n+1}_j$, $C^{n+1}_j$,  $A^{n+1}_{1-j}$, and $C^{n+1}_{1-j}$ have the wanted properties.
By construction $A^{n+1}_j=A^n_j\subseteq C^n_j=C^{n+1}_j$. We have that
\begin{equation*}A^{n+1}_{1-j}=A^n_{1-j}\cup U_{R_0(x_m,K(x_m,I))}\subseteq C^{n+1}_{1-j}\cup\bigcap_{y\notin U_{K(x_m,I)}} X\setminus U_{R_1(x_m,K(x_m,I),y)} =C^{n+1}_{1-j}
\end{equation*} since $A^n_{1-j}\subseteq C^n_{1-j}$ and by definition of $R_0$ and $R_1$ we have
\begin{equation*}
U_{R_0(x_m,K(x_m,I))}\subseteq \bigcap_{y\notin U_{K(x_m,I)}} X\setminus U_{R_1(x_m,K(x_m,I),y)}
\end{equation*}
We prove that $C^{n+1}_j\cap C^{n+1}_{1-j}=\emptyset$. We have by construction that
\begin{equation*} C^{n+1}_j=C^n_j\quad \text{ and }\quad C^{n+1}_{1-j}=C^n_{1-j}\cup \bigcap_{y\notin U_{K(x_m,I)}} X\setminus U_{R_1(x_m,K(x_m,I),y)}
\end{equation*}
By definition of the function $R_1$ and $K$ we have
\begin{equation*} \bigcap_{y\notin U_{K(x_m,I)}} X\setminus U_{R_1(x_m,K(x_m,I),y)}\subseteq U_{K(x_m,I)}\subseteq \bigcap_{i\in I}U_i
\end{equation*}
It therefore suffices to show that $C^n_j\cap \bigcap_{i\in I}U_i=\emptyset$. Let $r\in C^n_j$ using the definition of $C^n_j$ we have the following cases:\\
\\
\textbf{Case 1:} $r\in \bigcap_{i\in\mathbb{N}}X\setminus U_{f_j(i)}$ so $r\notin U_{f_j(s)}$ and $r\notin \bigcap_{i\in I} U_i$ since $f_j(s)\in I$.\\
\\
\textbf{Case 2:} There is $(z,l)\in E^n_j$ such that $r\in \bigcap_{y\notin U_l} X\setminus U_{R_1(z,l,y)}$. By definition of $I$ there is some $y$ such that $R_1(z,l,y)\in I$ and so $r\notin U_{R_1(z,l,y)}$.\\
\\
Therefore, the sequences $(C^n_0)_{n\in\mathbb{N}},(C^n_1)_{n\in\mathbb{N}},(A^n_0)_{n\in\mathbb{N}},$ and $(A^n_1)_{n\in\mathbb{N}}$ have the desired properties. Since the sequences $(F^n_0)_{n\in\mathbb{N}}$, $(F^n_1)_{n\in\mathbb{N}}$, $(E^n_0)_{n\in\mathbb{N}}$, and $(E^n_1)_{n\in\mathbb{N}}$ are defined recursively we have that they exist by $\Delta^0_1$ comprehension. As observed before we have that the set \begin{equation*}D=\bigcup_{n\in\mathbb{N}}A^n_0=X\setminus \bigcup_{n\in\mathbb{N}}A^n_1
\end{equation*}
will be the wanted clopen set. Both $D$ and its complement are $\Sigma^0_1$ so $D$ exists by $\Delta^0_1$ comprehension.
\end{proof}
\end{thm}

\begin{obs} Applying the previous result to the space $(\mathbb{N},(\{i\})_{i\in\mathbb{N}},k)$, where $k(i,i,i)=i$, gives us the classic recursion theory result that any pair of disjoint $\Pi^0_1$ sets can be separated by a recursive set.
\end{obs}

\begin{cor}\label{REG>G} $\textbf{RCA}_0$ proves that any uniformly regular CSC space is effectively homeomorphic to an effectively zero dimensional CSC space.\\
\\
\begin{proof} Let $R_0$ and $R_1$ be the functions that witness that $X$ is $uT_3$. For each pair $(x,i)$, where $x\in U_i$, let $C_0^{(i,x)}$ be the closure of the point $x$; which by regularity is contained in $U_{i}$, and $C^{(x,i)}_1=X\setminus U_i$. We point out that the closure of $x$ may not be a set, but it will be a coded closed collection. By Theorem~\hyperref[SNORM]{\ref{SNORM}} we have that there exists a sequence of clopen sets $(D_{(x,i)})_{x\in U_i}$ such that for all $x\in U_i$ we have $x\in D_{(x,i)}\subseteq X\setminus(X\setminus U_i)=U_i$. We define $k'$. Given $x\in D_{(y,j)}\cap D_{(z,l)}$, by the construction carried out in Theorem~\hyperref[SNORM]{\ref{SNORM}} we have that
\begin{equation*}D_{(y,j)}=\bigcup_{n\in\mathbb{N}}A^n_{0,(y,j)}=\bigcup_{n\in\mathbb{N}}\bigcup_{s\in F^n_{0,(y,j)}} U_s\quad \quad \text{ and } D_{(z,j)}=\bigcup_{n\in\mathbb{N}}A^n_{0,(z,l)}=\bigcup_{n\in\mathbb{N}}\bigcup_{t\in F^n_{0,(z,l)}} U_t
\end{equation*}
Let $m=\min\{n\in\mathbb{N}: x\in A^n_{0,(y,j)}\cap A^n_{0,(z,l)}\}$ and define
\begin{equation*}
s=\min \{s'\in F^m_{0,(y,j)}: x\in U_{s'}\} \quad \quad \text{ and }\quad \quad t=\min \{t'\in F^m_{0,(z,l)}: x\in U_{t'}\}
\end{equation*}
and let $k'(x,(y,i),(z,j))=(x,k(x,s,t))$. It is straightforward, using $R_0$, to show that $(X,(D_{(x,i)})_{x\in U_i},k')$ is a CSC space with a basis of clopen sets and it is effectively homeomorphic to $(X,(U_i)_{i\in\mathbb{N}},k)$. We define a function $G$ witnessing that $X$ is effectively zero dimensional. Let $x\in X$ and $i\in \mathbb{N}$ be such that $x\in U_i$. Observing the construction carried out in Theorem~\hyperref[SNORM]{\ref{SNORM}} we have that $X\setminus D_{(x,i)}=\bigcup_{n\in\mathbb{N}} A^n_1=\bigcup_{n\in\mathbb{N}}\bigcup_{j\in F^n_1} U_j$. Given $y\notin D_{(x,i)}$ let $m=\min\{n\in\mathbb{N}: y\in A^n_1\}$ and define $G(y,(x,i))=(y,\min \{j\in F^n_1: y\in U_j\})$. We have that $G$ witnesses that $(X,(D_{(x,i)})_{x\in U_i},k')$ is effectively zero dimensional and  exists by $\Delta^0_1$ comprehension.
\end{proof}
\end{cor}

\begin{prop}\label{G>Alg} An effectively zero dimensional CSC space $(X,(U_i)_{i\in\mathbb{N}},k)$ is effectively homeomorphic to a CSC space with an algebra of clopen sets.\\
\\
\begin{proof} Let $(V_i)_{i\in\mathbb{N}}$ be the sequence of sets given by $V_{2j}=U_i$ and $V_{2i+1}=X\setminus U_i$.
We then define $A_m=\bigcap_{i\in m}\bigcup_{j\in i} V_i$. We define  $Int(n,m)$ to be the number that codes the union of $m$ and $n$. For $m\in\mathbb{N}$ coding the set
\begin{equation*}
\{\{a_{0,0},\dots,a_{0,n_0}\},\dots,\{a_{r,0},\dots,a_{r,n_r}\}\}
\end{equation*}
then $Comp(m)$ will code the set of all sets of the form
\begin{equation*}
\{a_{0,f(0)}^c,a_{1,f(1)}^c,\dots,a^c_{r,f(r)}\}
\end{equation*}
where $f:r\rightarrow \max\{n_j:j\leq r\}$ is a coded function such that for all $j\leq r $ we have $f(j)\leq n_j$ and where $a^c=a-1$ if $a$ is odd and $a^c=a+1$ if $a$ is even. One observes that  $A_{Comp(m)}=X\setminus A_m$.\\
\\
Let $Id:(X,(U_i)_{i\in\mathbb{N}},k)\rightarrow (X,(A_m)_{m\in\mathbb{N}},Int,Comp)$ be the identity function from $X$ to $X$, we show it's an effective homeomorphism. The function  $(x,i)\mapsto \{\{2i\}\}$ witnesses that $Id$ is effectively continuous. We need to show that there is a function $v$ such that for any $(x,m)$  such that $x\in A_m$ then $x\in U_{v(x,m)}\subseteq A_m$.  We define $v$ on the complexity of the set $m$ codes. If $m=\{\{2i\}\}$ then we set $v(x,m)=i$. If $m=\{\{2i+1\}\}$ then $v(x,m)=G(x,i)$. For $m=\{\{a_0,\dots,a_n\}\}$ where $a_0<a_1<\dots<a_n$ then we define $v(x,m)=v(x,\{\{ a_0\}\})$. Finally for any $m=\{m_0,\dots,m_r\}$ then let $I=\{v(x,\{m_i\}): i\leq r\}$ we define
\begin{equation*}
v(x,m)=K(x,I)
\end{equation*}
The function $v$ is defined primitively recursively, and so it exists by $\Delta^0_1$ comprehension, and it verifies that the identity is effectively open between $(X,(U_i)_{i\in\mathbb{N}},k)$ and $(X,(A_m)_{m\in\mathbb{N}},Int,Comp)$.
\end{proof}
\end{prop}

\begin{thm}\label{Big4} $\textbf{RCA}_0$ proves that for every CSC space  $(X,(U_i)_{i\in \mathbb{N}},k)$ the following are equivalent:
\begin{enumerate}[label={(\arabic*)}]
\item $X$ is $uT_3$.
\item $X$ is $T_0$ and effectively homeomorphic to a CSC space with an algebra of clopen sets.
\item $X$ is effectively homeomorphic to a subspace of $\mathbb{Q}$.
\item $X$ is metrizable.
\end{enumerate}
\begin{proof}
The case in which $X$ is finite is trivial, so we assume $X$ is infinite and let $(x_i)_{i\in\mathbb{N}}$ be the elements of $X$ enumerated increasingly. $(3\rightarrow 4)$ is trivial. \\
\\
$(1\rightarrow 2)$ follows from Corollary~\hyperref[REG>G]{\ref{REG>G}} and Proposition~\hyperref[G>Alg]{\ref{G>Alg}}. $(2\rightarrow 1)$ follows from the fact that any CSC space with an algebra of clopen sets is effectively zero dimensional, which by Proposition~\hyperref[G>Reg]{\ref{G>Reg}} is $uT_3$. \\
\\
$(4\rightarrow 2)$ We have already seen in Proposition~\hyperref[metstruc]{\ref{metstruc}} that all countable metric spaces have a CSC structure that consists of clopen balls. In particular given a metric space $(X,d)$ there exists an $a\in \mathbb{R}_{>0}$ such that $(X, (B(x,q\cdot a) )_{q\in \mathbb{Q}_{>0},x\in X},k)$ is a CSC space with a basis of clopen sets. We define $G:X\times (X\times\mathbb{Q}_{>0})\rightarrow (X\times \mathbb{Q}_{>0})$ which will send $(x,(y,q))$, where $x\notin B(y,q\cdot a)$, to $(x,\frac{1}{n})$ where $n$ is the least $n$ such that $\frac{1}{2^n}\cdot a < d(x,y)- q\cdot a$. Since the CSC structure of $X$ is effectively zero dimensional, by Proposition~\hyperref[G>Alg]{\ref{G>Alg}}, it is effectively homeomorphic to a space with an algebra of clopen sets.\\
\\
$(2\rightarrow 3)$ Assume that $X$ has an algebra of clopen sets. We will show there is an embedding $f$ from $X$ to the rational points of the Cantor space, which we will denote as $2^{\mathbb{N}}\cap \mathbb{Q}$. We observe that $2^{\mathbb{N}}\cap \mathbb{Q}$ is a CSC space with an algebra of clopen sets of the form $[[\tau]]=\{g\in 2^{\mathbb{N}}\cap\mathbb{Q}:\tau\sqsubseteq g\}$, and it effectively embeds into $[0,1]\cap \mathbb{Q}$.\\
\\
We construct recursively the following things
\begin{itemize} 
\item A total surjection $\phi:\mathbb{N}\rightarrow \mathbb{N}$ such that for all $n\in\mathbb{N}$ and for all $i<j<2n$ we have
\begin{equation*}\{m<2n:x_i\in U_{\phi(m)}\}\neq \{m< 2n:x_j\in U_{\phi(m)}\}
\end{equation*}
The function $\phi$ will define a sequence of basic open sets $(U_{\phi(n)})_{n\in\mathbb{N}}$.
For $\tau\in 2^{<\mathbb{N}}$ we write
\begin{equation*}
N_\tau=\{x\in X:\forall m<|\tau|\;(x\in U_{\phi(m)}\leftrightarrow \tau(m)=1)\}
\end{equation*}
Since we are assuming that $X$ has an algebra of sets, we have that for every sequence $\tau$ that
\begin{equation*}
N_\tau=\left(\bigcap_{i< |\tau|\wedge \tau(i)=1}U_{\phi(i)}\right)\cap\left(\bigcap_{i< |\tau|\wedge \tau(i)=0}U_{Comp(\phi(i))}\right)
\end{equation*}
So there is some $j$ such that $N_\tau=U_j$. Furthermore, there is a uniformly effective procedure to find an index $j$ from $\tau$. In general, $N_\tau$ might be empty, but this will not affect the construction.
\item $f(x_0),\dots,f(x_n)\in 2^{<\mathbb{N}}\cap \mathbb{Q}$.
\item A function $r_n:2^{\leq 2n}\rightarrow 2^{\leq 2n}$ which is monotone, level preserving and for all $\sigma \in \text{rng}(r_n)$ there is an $i\leq n$ such that $f(x_i)\in N_{\sigma}$. We also require that for all $i\leq n$ we have $f(x_i)\in [[r_n(seq_{2n}(x_i))]]$ where $seq_{2n}(x_i)$ is a sequence of length $2n$ such that
\begin{equation*}seq_{2n}(x_i)(j)=\begin{cases} 1 \quad \text{ if }\quad x\in U_{\phi(j)} \\
 0 \quad  \text{ if }\quad x\notin U_{\phi(j)}
\end{cases}
\end{equation*}
This ensures us that for all $\tau$ of length $n$ $f(N_\tau)=f(X)\cap [[r_n(\tau)]]$ and that $f$ will be injective.
Finally, we will require for all $m<n$ that $r_m\subseteq r_n$.
\end{itemize}
Assume we have $f(x_0),\dots,f(x_{n}),r_n,\phi|_{< 2n}$ and that for all $\sigma\in 2^{=2n}$ $[[r_n(\sigma)]]$ contains a unique $f(x_i)$.\\
\\
\textbf{Step $1$:} We set $\phi(2n)=\min\mathbb{N}\setminus \text{rng}(\phi|_{<2n})$, this ensures us that $\phi$ will be surjective. 
Given $\tau\in 2^{<\mathbb{N}}$ of length $2n-1$ there exists exactly one $j\leq n$ such that $f(x_j)\in [[r_n(\tau)]]$. Let $i\leq 1$ be such that $f(x_j)\in [[r_n(\tau)^\frown(i)]]$,  we define $r_n'(\tau^\frown(0))=r_n'(\tau^\frown(1))=r_n(\tau)^\frown(i)$.\\
\\
\textbf{Step 2:} There is, by  construction of $r'_n$, exactly one $j\leq n$ such that $r'_n(seq_{2n+1}(x_{n+1}))=r'_n(seq_{2n+1}(x_j))$. We define
\begin{equation*}
\phi(2n+1)=\min \{ l\in\mathbb{N}\setminus\text{rng}(\phi|_{<2n+1}):x_{n+1}\in U_l\leftrightarrow x_j\notin U_l\}
\end{equation*}
That is, the first index for a basic open set not in the range of $\phi|_{<2n+1}$ that separates $x_{n+1}$ and $x_j$. The search for such $l$ will eventually terminate since we are assuming $X$ is $T_0$. Assume that $x_j\in U_{\phi(2n+1)}$ and $x_{n+1}\notin U_{\phi(2n+1)}$, the proof for the other case is analogous. Let $\rho$ be the unique sequence of length $2n+2$ such that $f(x_j)\in [[\rho]]$ and let $i=\rho(2n+1)$. For all $\sigma\in 2^{=2n+1}$ such that $r'_n(\sigma)= r'_n(seq_{2n+1}(x_{n+1}))=r'_n(seq_{2n+1}(x_j))$ we define
\begin{equation*}
r_{n+1}(\sigma^\frown(0))=r'_n(\sigma)^\frown (1-i) 
\end{equation*}
\begin{equation*}
r_{n+1}(\sigma^\frown(1))=r'_n(\sigma)^\frown (i)=\rho
\end{equation*}
For all other sequences of length $2n+2$ we define $r_{n+1}$ as in step $1$
.\\
\\
In either case, we are ensured that for all $j\leq n$, we have
\begin{equation*}
f(x_j)\in [[r_{n+1}(seq_{2n+2}(x_j))]]
\end{equation*}
\textbf{Step 3:} We define $f(x_{n+1})$ to be the least Cantor rational in the open set $[[r_{n+1}(seq_{2n+2}(x_{n+1}))]]$. This ends the construction.\\
\\
Given the functions $f,\phi$ and $r=\bigcup_{n\in\mathbb{N}}r_n$ we show that $f$ is an effective homeomorphism with its range.\\
\\
We first show that the range of $f$ exists. By definition $y\in \text{rng}(f)\leftrightarrow \exists x\; f(y)=x$, and so the range is $\Sigma^0_1$ definable relative to $f$. We also have that $y\in \text{rng}(f)$ if and only if for all $n\in \mathbb{N}$ there is a $\sigma\in 2^{=n}$ such that $y\in [[r_n(\sigma)]]$, which is a $\Pi^0_1$ definition. So we have that $\text{rng}(f)$ is $\Delta^0_1$ definable and therefore exists by $\Delta^0_1$ comprehension. We note that $\Pi^0_1$ definition of $\text{rng}(f)$ tells us it is a closed subspace of $2^{\mathbb{N}}\cap \mathbb{Q}$.\\
 \\
Let $x_j\in X$ be a point and $U_{\phi(i)}$ an open set containing it. We wish to effectively find an open neighborhood of $f(x)$ that is contained in $f(U_{\phi(i)})$. Let $n>\max\{i,j\}$, we have that $seq_{2n}(x_j)(i)=1$ and so
\begin{equation*} x_j\in N_{seq_{2n}(x_j)}\subseteq U_{\phi(i)}
\end{equation*}
Which by construction implies
\begin{equation*}
f( x_j)\in f(N_{seq_{2n}(x_j)})=f(X)\cap [[r_n(seq_{2n}(x_j))]]\subseteq f(U_{\phi(i)})
\end{equation*} 
This implies that $f:X\rightarrow f(X)$ is effectively open.\\
\\
Let $[[\tau]]$ be a basic open set of the Cantor space and $f(x_j)\in [[\tau]]$. By construction we have that $\tau\in\text{rng}(r)$, since $r$ is level preserving we can effectively find $\rho=r^{-1}(\tau)$. We   have that $x\in N_{\rho}$ and $f(N_\rho)\subseteq[[\tau]]$ by construction. This shows that $f$ is effectively continuous. Therefore $f$ is an effective embedding of $X$ into $2^{<\mathbb{N}}\cap \mathbb{Q}$.
\end{proof}
\end{thm}
\noindent
Sierpi{\' n}ski \cite{Sierp} proved using a back and forth construction that any pair of countable sets of $\mathbb{R}^n$ without isolated points are homeomorphic. That the back and forth method, used by Sierpinski, can be carried out in $\textbf{RCA}_0$ for effectively zero dimensional $T_0$ spaces was first observed by Soldà (unpublished).  

\begin{thm}\label{sol} \textbf{(Sold{\`a} \cite{Solda})}\footnote{Sold{\`a} proved this for spaces which are supereffectively zero dimensional. A CSC space $(X,(U_i)_{i\in \mathbb{N}},k)$ is said to be supereffectively zero dimensional if there exists a function $\chi:\mathbb{N}\rightarrow \mathbb{N}$ such that $X\setminus U_i=U_{\chi(i)}$. Sold{\`a} conjectures that having the complement be a basic set may not be necessary but rather it suffices that $\chi(i)$ is an open code for the collection $U_{\chi(i)}$. This weaker definition of supereffectively zero dimensional CSC space coincides with our notion of being effectively zero dimensional.} $\textbf{RCA}_0$ proves that any two $T_0$ CSC spaces with an algebra of clopen sets and without isolated points are effectively homeomorphic.
\end{thm}

\begin{cor}\label{Qchar} Any non empty $uT_3$ CSC space without isolated points is effectively homeomorphic to $\mathbb{Q}$ with the order topology.\\
\\
\begin{proof}
By Corollary~\hyperref[Reg>G]{\ref{REG>G}} and Proposition~\hyperref[G>Alg]{\ref{G>Alg}} any $uT_3$ CSC space has an algebra of clopen sets, so by Theorem~\hyperref[sol]{\ref{sol}}, any $uT_3$ space without isolated points will be homeomorphic to $\mathbb{Q}$.
\end{proof} 
\end{cor}

\begin{obs}\label{ACAT3Q} Since over $\textbf{ACA}_0$ all $T_3$ CSC space are $uT_3$ we have that $\textbf{ACA}_0$ proves that every $T_3$ space is homeomorphic to a subset of $\mathbb{Q}$ and is homeomorphic to a space with an algebra of clopen sets. 
\end{obs}

\begin{thm}\label{stuff>aca} The following are equivalent over $\textbf{RCA}_0$:
\begin{enumerate}[label={(\arabic*)}]
\item Arithmetic comprehension.
\item Every $T_3$ $eT_2$ CSC space is effectively homeomorphic to a metric space.
\item Every $T_3$ $eT_2$ CSC space is effectively homeomorphic to a CSC space with an algebra of clopen sets.
\item Every $T_3$ $eT_2$ CSC space is effectively homeomorphic to a subset of the rationals.
\end{enumerate}
\begin{proof} By Theorem~\hyperref[Big4]{\ref{Big4}}, $(2)$, $(3)$, and $(4)$ are equivalent over $\textbf{RCA}_0$. By Proposition~\hyperref[G>Reg]{\ref{G>Reg}}, all effectively zero dimensional CSC spaces, and therefore all CSC spaces with an algebra of clopen sets, are uniformly regular. By Proposition~\hyperref[Regstab]{\ref{Regstab}}, $(4)$ is equivalent to every regular $eT_2$ CSC space is uniformly regular, which is equivalent to arithmetic comprehension over $\textbf{RCA}_0$ by Theorem~\hyperref[worst proof]{\ref{worst proof}}.
\end{proof}
\end{thm}
\noindent
It is clear that a CSC space that is not $T_0$ cannot be homeomorphic to a subspace of the rationals. However, for non $T_0$ uniformly regular CSC spaces we can consider a pseudometric.

\begin{defin} A pseudometric on a set $X$ is a function $d:X\times X\rightarrow \mathbb{R}_{\geq 0}$ that is symmetric and satisfies the triangle inequality. We observe that similar to the case of a metric space \hyperref[metruc]{\ref{metstruc}}, we can associate a CSC structure to a set with a pseudometric. As in the case for metric spaces, we get that the CSC structure of a pseudometric will be unique up to effective homeomorphism and effectively zero dimensional. We say that a CSC space is pseudometrizable if it is effectively homeomorphic to the CSC structure of a pseudometric space. 
\end{defin}

\begin{prop}\label{merenda} $\textbf{RCA}_0$ proves that every uniformly regular CSC space is effectively homeomorphic to a CSC space with a pseudometric.\\
\\
\begin{proof} Let $(X,(U_i)_{i\in\mathbb{N}})$ be a uniformly regular CSC space by Corollary~\hyperref[REG>G]{\ref{REG>G}} and Proposition~\hyperref[G>Alg]{\ref{G>Alg}} we can assume without loss of generality that $X$ has an algebra of clopen sets. Let $U[x]=\{i\in\mathbb{N}:x\in U_i\}$ we define  $d$ by
\begin{equation*}
d(x,y)=\sum_{k\in U[x]\triangle U[y]} \frac{1}{2^k}
\end{equation*}
The function $d$ is symmetric, and the triangle inequality follows from the following property of the symmetric difference
\begin{equation*}
\forall A,B,C\;(A\triangle C)\subseteq ((A\triangle B)\cup(B\triangle C))
\end{equation*}
By definition of $d$ we have
\begin{equation*}B(x,\frac{1}{2^k})=N_{U^k[x]}=\{y\in X: \forall i<k\;(x\in U_i\leftrightarrow y\in U_i)\}
\end{equation*}
So we have that $B(x,\frac{1}{2^k})$ is the boolean combination of clopen sets and is, therefore, clopen. We show that the identity $Id:(X,(U_i)_{i\in\mathbb{N}}, Int, Comp)\rightarrow (X,(B(x,\frac{1}{2^k}))_{k\in\mathbb{N}},k)$ is an effective homeomorphism. 
For $x\in X$ and $i\in \mathbb{N}$ such that $x\in U_i$ then we have that $x\in B(x,\frac{1}{2^i})\subseteq U_i$ and so the map $(x,i)\mapsto i$ verifies that the identity is effectively open. The map that sends $(x,k)$ to the code for the open set
\begin{equation*}
\left(\bigcap_{i< k\wedge x\in U_i} U_i\right)\cap \left(\bigcap_{i< k\wedge x\notin U_i} X\setminus U_{i}\right)
\end{equation*}
is computable relative to the functions $Comp$ and $Int$ and it verifies that the identity is effectively continuous.
\end{proof}
\end{prop}
\noindent
Just as in the $T_0$ case, we get the following result:
\begin{thm}\label{noT0} $\textbf{RCA}_0$ proves that for any CSC space $(X,(U_i)_{i\in\mathbb{N}},k)$ the following are equivalent:
\begin{enumerate}[label={(\arabic*)}]
\item $X$ is uniformly regular.
\item $X$ is effectively zero dimensional.
\item $X$ is effectively homeomorphic to a CSC space with an algebra of clopen sets.
\item $X$ is pseudometrizable.
\end{enumerate}
\begin{proof}
$(1\rightarrow 2)$ is Corollary~\hyperref[REG>G]{\ref{REG>G}}, $(2\rightarrow 3)$ is Proposition~\hyperref[G>Alg]{\ref{G>Alg}}, and $(3\rightarrow 4)$ is Proposition~\hyperref[merenda]{\ref{merenda}}. $(4\rightarrow 1)$ follows from the fact that the fact that the CSC structure on a pseudometric space is effectively zero dimensional.
\end{proof}
\end{thm}

\section{Compactness}
In this section, we will look at compact CSC spaces, especially $T_2$ compact CSC spaces. It is a classic result of general topology that compact $T_2$ spaces are $T_3$ and normal. We also have that compact $T_2$ CSC spaces are well-orderable. We will show that these characterizations of $T_2$ compact spaces can be carried out in $\textbf{ACA}_0$. We will also lay out a few technical lemmas that will be needed for the subsequent sections.
\begin{defin} A CSC space $(X,(U_i)_{i\in\mathbb{N}},k)$ is said to be compact if for every $I\subseteq \mathbb{N}$ such that $X=\bigcup_{i\in I} U_i$ there exists a finite $a\subseteq I$ such that $X=\bigcup_{i\in a} U_i$. The covering relation for the space $X$ is the set $C=\{a\in\mathbb{N}:X=\bigcup_{i\in a} U_i\}$. We say that $X$ is effectively compact if it is compact and has a covering relation. We say that $X$ is sequentially compact if every sequence contains a converging subsequence. A subset of a CSC space is said to be compact or effectively compact if its subspace topology is respectively compact or effectively compact. 
\end{defin}

\begin{lemma}\label{coverlem} \textbf{(Dorais \cite[Proposition 7.5]{Dorais})} $\textbf{RCA}_0$ proves that every linear order with the order topology has a covering relation.
\end{lemma}

\begin{obs}\textbf{(Dorais  \cite[Example 5.4]{Dorais})} Over $\textbf{ACA}_0$ a CSC space is compact if and only if it is sequentially compact.
\end{obs}

\begin{obs} $\textbf{RCA}_0$ proves that the union of two effectively compact sets is effectively compact. In particular, Given $C_0,C_1$ the covering relations of $K_0,K_1$ respectively, then $C_0\cap C_1$ is the covering relation for $K_0\cup K_1$.
\end{obs}

\begin{thm}\label{effcomp<>ACA} Over $\textbf{RCA}_0$ the following are equivalent:
\begin{enumerate}[label={(\arabic*)}]
\item Arithmetic comprehension.
\item Every compact CSC space is effectively compact (Dorais \cite[Example 3.5]{Dorais}).
\item Every compact $uT_3$ CSC space is effectively compact.
\end{enumerate}
\begin{proof} $(1\rightarrow 2)$ since the covering relation for a CSC space is arithmetically definable. $(2\rightarrow 3)$ is trivial. For $(3\rightarrow 1)$ consider the CSC space $(X,(U_i)_{i\in\mathbb{N}},k)$ where $X=\mathbb{N}\cup\{\infty\}$:
\begin{equation*}
U_{2i}=\{i\}
\end{equation*}
\begin{equation*}
U_{2(n,s)+1}= \{\infty\}\cup\{j\in\mathbb{N}: j\geq n\wedge \forall e\in s\; \Phi^A_e(e){\downarrow}_{\leq j}\rightarrow \exists r<j \Phi^A_e(e){\downarrow}_{\leq r}\}
\end{equation*}
that is all the $j\geq n$ such that for all $e\in s$ the Turing machine of index $e$ does not halt for the first time at step $j$. We define $k$ that sends $(n,a,b)\mapsto 2n$ and $(\infty,2(i,s)+1,2(j,t)+1)$ to $2(\max\{i,j\},r)+1$ where $r$ is the code for the union of $s$ and $t$. $X$ is $uT_3$ follows from the fact that the basis consists of clopen sets and all points except $\infty$ are isolated. By assumption, $X$ is effectively compact and, therefore, has a covering relation $C$. We have that $\{2(0,\{e\})+1\}\in C\leftrightarrow \Phi^A_e(e){\uparrow}$ and so $A'\leq_T C$.
\end{proof}
\end{thm}

\begin{prop}\label{effhomcom} \textbf{(Dorais  \cite[Proposition 3.2]{Dorais})}  $\textbf{RCA}_0$ proves that being effectively compact is preserved under effective homeomorphism.
\end{prop}

\begin{prop}\label{closedincomp} \textbf{(Dorais  \cite[Proposition 3.6]{Dorais})} $\textbf{RCA}_0$ proves that for any compact CSC space $(X,(U_i)_{i\in\mathbb{N}},k)$ any effectively closed $C\subseteq X$ subset is a compact subspace.
\end{prop}

\begin{prop}\label{Kcode} \textbf{(Dorais \cite[Proposition 6.3]{Dorais})}  $\textbf{RCA}_0$ proves that for any $eT_2$ CSC space $(X,(U_i)_{i\in\mathbb{N}},k)$ and $K\subseteq X$ an effectively compact set, then $K$ is effectively closed.
\end{prop}
\noindent
Modifying Dorais' proof of the previous proposition, we get the following results:

\begin{prop}\label{Kcodes} $\textbf{RCA}_0$ proves that for any $eT_2$ CSC space $(X,(U_i)_{i\in\mathbb{N}},k)$ and $(K_n,C_n)_{n\in\mathbb{N}}$ sequence where $K_n$ is an effectively compact subset of $X$ and $C_n$ is the covering relation for $K_n$, then $(K_n)_{n\in\mathbb{N}}$ is uniformly effectively closed.\\
\\
\begin{proof} Let  $H_0,H_1$ the functions that witness that $X$ is $eT_2$. Define $p:\mathbb{N}\times X \rightarrow \mathbb{N}$ such that for any $n\in \mathbb{N}$ and $x\notin K_n$ we have
\begin{equation*}
p(n,x)=\min\{F \text{ code for a finite subset of }K_n: \{H_0(x,y):y\in F\}\in C_n\}
\end{equation*}
We observe that $p$ is computable relative to $(K_n,C_n)_{n\in\mathbb{N}}$ and so it exists by $\Delta^0_1$ comprehension. For each $x\notin K_n$ we define $s(n,x)=K(x,\{H_1(x,y):y\in p(n,x)\})$. We have that by construction
\begin{equation*}
x\in U_{s(n,x)}\subseteq \bigcap_{y\in p(n,x)} U_{H_1(x,y)}\subseteq X\setminus \bigcup_{y\in p(n,x)}U_{H_0(x,y)}\subseteq X\setminus K_n
\end{equation*}
We have that $s$ exists by $\Delta^0_1$ comprehension and, for each $n$, we have that $s(n,\cdot)$ is a code for the closed set $K_n$, so $(K_n)_{n\in\mathbb{N}}$ is uniformly effectively closed.
\end{proof}
\end{prop}

\begin{thm}\label{compact T3} $\textbf{RCA}_0$ proves that every effectively compact $eT_2$ CSC space  is $uT_3$.\\
\\
\begin{proof}
Let $(X,(U_i)_{i\in\mathbb{N}},k)$ be an effectively compact $eT_2$ CSC space. Let $C$ be the covering relation of $X$ and $H_0$ and $H_1$ be functions that witness that $X$ is $eT_2$. We define $R_0$ and $R_1$ on $X$. For any $x\in U_i$ we have that $\{U_i\}\cup \{U_{H_1(x,y)}:y\notin U_i\}$ is a covering of $X$. Let 
 $F\subseteq X\setminus U_i$ be the least finite set such that $\{i\}\cup\{H_1(x,y):y\in F\}\in C$. Define
\begin{equation*}
R_0(x,i)=K(x,\{H_0(x,y):y\in F\})
\end{equation*}
and for each $z\notin U_i$ we define
\begin{equation*} R_1(x,i,z)=\min\{y\in F: z \in U_{H_1(x,y)}\}
\end{equation*}
By definition, we have
\begin{equation*}
x\in U_{R_0(x,i)}\subseteq\bigcap_{y\in F} U_{H_0(x,y)}\subseteq X\setminus \bigcup_{y\in F} U_{H_1(x,y)}=X\setminus\bigcup_{z\notin U_i} U_{R_1(x,i,z)}\subseteq  U_i
\end{equation*}
So the functions $R_0$ and $R_1$ witness that $X$ is $uT_3$ and are recursive relative to the sets and functions $C,(X,(U_i)_{i\in\mathbb{N}},k)$, $H_0$, and $H_1$ so they exist by $\Delta^0_1$ comprehension.
\end{proof}
\end{thm}

\begin{cor} \label{acacomp} $\textbf{ACA}_0$ proves that every $T_2$ compact CSC space is $T_3$.
\end{cor}

\begin{prop}\label{wocomp} (See \cite[Proposition 6.18]{Greenberg}) $\textbf{RCA}_0$ proves that every well-order with a maximal element is compact with respect to its order topology.
\end{prop}

\begin{obs}  $\textbf{RCA}_0$ proves that every linear order with its order topology has a covering relation, so $\textbf{RCA}_0$ proves that every well-order with maximal element is effectively compact. Using the same proof, we can show that the upper limit topology on a well-order with a maximal element is effectively compact. We will see later that over $\textbf{ACA}_0$ all compact $T_2$ CSC spaces will be homeomorphic to a well-order with maximum element.
\end{obs}

\begin{defin}
We say that a CSC space $(X,(U_i)_{i\in\mathbb{N}},k)$ is strongly compact if for every sequence of open set $(A_n)_{n\in\mathbb{N}}$ there exists a finite set $F\subseteq\mathbb{N}$ such that $(A_n)_{n\in F}$ covers $X$. A CSC space is strongly effectively compact if it is strongly compact and has a covering relation. We note that every sequentially compact CSC space is strongly compact by the usual topological argument (see \cite[Proposition 4.1]{Dorais}).
\end{defin}
\noindent
Clearly, over $\textbf{ACA}_0$, every compact CSC space is strongly compact. We would like to show that over $\textbf{RCA}_0$, strong compactness is a stronger property than compactness. Shafer proved in the absence of arithmetic comprehension, there can be compact spaces which are not strongly compact \cite[Proposition 3.10]{Shafer}. The CSC space considered by Shafer is discrete and therefore $T_2$, however, it is not effectively $T_2$ nor effectively compact. We would like to extend this result to show that in absence of arithmetic comprehension being strongly compact is a stronger property even for effectively compact $eT_2$ spaces. We need the following definition and lemma.

\begin{defin} Over $\textbf{RCA}_0$ a linear order $(L,<_L)$ is complete if there does not exist a pair of non empty subsets $A$ and $B$ such that $L=A\sqcup B$, $\forall a\in A\;\forall b\in B\;(a<_L b)$, $A$ does not have a maximum and $B$ does not have a minimum.
\end{defin}

\begin{lemma}\label{Shaf} \textbf{(Shafer \cite[Lemma 4.4]{Shafer})} For any tree $T\subseteq \mathbb{N}^{<\mathbb{N}}$, if $T$ is well-founded, then $(T,<_\text{KB})$ is complete.
\end{lemma}
\noindent
Making use of Shafer's \cite{Shafer} result that an infinite finitely branching tree that is well-founded will fail to be strongly compact we get the following result.
\begin{prop}\label{strongcomp} The following are equivalent over $\textbf{RCA}_0$:
\begin{enumerate}[label={(\arabic*)}]
\item Arithmetic comprehension.
\item Every $eT_2$ effectively compact CSC space is strongly compact.
\end{enumerate}
\begin{proof}
$(1\rightarrow 2)$ follows from the fact that over $\textbf{ACA}_0$ compact is equivalent to strongly compact.\\
\\
$(2\rightarrow 1)$ Let $f:\mathbb{N}\rightarrow \mathbb{N}$ be an injection, we wish to prove that $\text{rng}(f)$ exists as a set. Define
\begin{equation*}
T=\{\sigma: \forall m,n<|\sigma |\;((f(n)=m\leftrightarrow\sigma(m)=n+1)\wedge(\sigma(m)>0\leftrightarrow f(\sigma(m)-1)=m))\}
\end{equation*}
That is, the tree $T$ contains all sequences that approximate the inverse of the function $g$.
The same construction is used in \cite[Theorem III.7.2]{Simp} to prove that K{\H o}nig's lemma implies arithmetic comprehension. We have that $T$ has the following properties:
\begin{enumerate}[label={(\arabic*)}]
\item Every $\sigma\in T$ has at most two immediate $\sqsubseteq$-successors.
\item $T$ has at most one branch, and the unique branch of $T$ can compute $\text{rng}(f)$.
\end{enumerate}
We prove that the Kleene-Brouwer order on $T$ is of order type $\omega+\omega^*$, that is, $T$ is infinite and for every $\sigma\in T$ either $\mathopen]-\infty,\sigma\mathclose[_{<_\text{KB}}$ is finite and we say $\sigma$ is in the $\omega$ part or $\mathopen]\sigma,+\infty[_{<_\text{KB}}$ is finite and we say $\sigma$ is in the $\omega^*$ part. Given $\sigma\in T$, we consider the following cases:\\
\\
\textbf{Case 1:} If $\forall j<|\sigma|\;(\exists n\; f(n)=j\rightarrow \sigma(j)=n+1)$ then we have that there is no element to the right of $\sigma$ in $T$ and so the only elements that are $<_\text{KB}$ greater than $\sigma$ are its initial segments. This means that in $\mathopen]\sigma,+\infty\mathclose[_{<_\text{KB}}$ has $|\sigma|$ many elements and so $\sigma$ is in the $\omega^*$ part.\\
\\
\textbf{Case 2:} If there is a  $j<|\sigma|$ and $n\in\mathbb{N}$ such that $f(n)=j$ and $ \sigma(j)=0$ then we have that $\sigma$ and any sequence to the left of $\sigma$ in $T$ cannot have an extension of length more than $n$ in $T$. This in particular means that $\mathopen]-\infty,\sigma\mathclose]_{<_\text{KB}}$ has at most $2^n$ elements and so $\sigma$ is in the $\omega$ part.
\\
\\
We note that by the fact that $(T,<_{\text{KB}})$ has order type $\omega+\omega^*$ we have that every element of $T$ which isn't the maximum has a successor and every element other than the minimum will have a predecessor.
In particular, we have that the upper limit topology on $(T,<_\text{KB})$ is $eT_2$ and discrete so it cannot be strongly compact. $T$ does however have a covering relation, we prove that it is compact. Assume that $T$ is not compact, and let $(\mathopen]a_n,b_n\mathclose[)_{n\in\mathbb{N}}$ be an infinite covering of $T$ which does not admit a finite subcovering. We have that every interval in $(T,<_\text{KB})$ will be either finite or cofinite, so for every $n\in\mathbb{N}$ we have that $\mathopen]a_n,b_n\mathclose[$ will be finite.
We show that both the $\omega$ part and the $\omega^*$ part of $(T,<_{\text{KB}})$ will be $\Sigma^0_1$ definable. We have that $x$ is in the $\omega$ part if there exists a finite sequence $\sigma$ such that $a_{\sigma(0)}=-\infty$ and for all $j<|\sigma|-1$ we have $b_j\leq_{\text{KB}} a_{j+1}$ and $x\in \mathopen]a_{\sigma(|\sigma|-1)},b_{\sigma(|\sigma|-1)}[$. That is, there exists a finite sequence of intervals of the covering whose union is an initial segment of $(T,<_{\text{KB}})$ containing $x$. In a similar fashion we may give a $\Sigma^0_1$ definition of the $\omega^*$ part of $(T,<_{\text{KB}})$. By $\Delta^0_1$ comprehension, we have that both the $\omega$ and the $\omega^*$ part of $T$ exist. So the Kleene-Brouwer order of $T$ is not complete. By Lemma~\hyperref[Shaf]{\ref{Shaf}} it will have a branch which will compute $\text{rng}(f)$.
\end{proof}
\end{prop}

\begin{prop}\label{univcomp} The following are equivalent over $\textbf{RCA}_0$:
\begin{enumerate}[label={(\arabic*)}]
\item $\Pi^1_1$ comprehension.
\item For every sequence of CSC spaces $(X^j)_{j\in\mathbb{N}}$ the set $\{j\in\mathbb{N}: X^j \text{ is compact}\}$ exists.
\end{enumerate}
\begin{proof} $(1\rightarrow 2)$ follows from the fact that being compact is a $\Pi^1_1$ definable property. We first show that $(2)$ implies arithmetic comprehension. Let $A$ be a set, we define for each $e\in \mathbb{N}$ the space $X^e=\{t\in\mathbb{N}:\Phi^A_e(e){\downarrow}_{\leq t}\}$ with the discrete topology. We observe that $X^e$ is compact if $X^e$ is finite which holds if and only if $\Phi^A_{e}(e){\uparrow}$. So $\mathbb{N}\setminus A'=\{e\in\mathbb{N}:X^e \text{ is compact}\}$. So $(2)$ implies that the Turing jump of every set exists which implies arithmetic comprehension. It therefore suffices to show  $(2\rightarrow 1)$ over $\textbf{ACA}_0$.\\
\\
We show that $(2)$ implies that for any sequence of trees $(T_i)_{i\in\mathbb{N}}$ the collection $
\{i\in\mathbb{N}: T_i \text{ is well-founded}\,\}$
exists, which by Theorem~\hyperref[Pioneone]{\ref{Pioneone}} is equivalent to $\Pi^1_1$ comprehension. It suffices to prove that a tree $T$ is well-founded if and only if $T$ with the upper limit topology with respect to the Kleene-Brouwer ordering is compact. If $T$ is well-founded, then  $(T,<_\text{KB})$ is well-ordered, which implies that its order topology is equal to its upper limit topology. By Proposition \hyperref[wocomp]{\ref{wocomp}} and the fact that $\emptyset\in T$ is $<_{\text{KB}}$-maximum, we have that $(T,<_\text{KB})$ with the upper limit topology is compact. If $T$ is not well-founded then $(T,<_\text{KB})$ has an infinite descending chain $(x_j)_{j\in\mathbb{N}}$. The space $T$ admits a partition of basic open sets
\begin{equation*}
\{T\setminus{\uparrow}\{x_j:j\in\mathbb{N}\},\mathopen]x_0,+\infty\mathclose[\,\}\cup\{\,\mathopen]x_{j+1},x_j]:j\in\mathbb{N}\}  
\end{equation*}
so $T$ is not compact.
\end{proof}
\end{prop}

\begin{obs} In the previous proof, we saw that for every tree, there is a canonically defined CSC space which is compact if and only if the tree is well-founded. This implies that being a compact CSC space is a universal $\Pi^1_1$ formula over $\textbf{ACA}_0$, so in particular, compactness for CSC space cannot be expressed by a $\Sigma^1_1$ formula. Furthermore, Shafer proved that over $\textbf{WKL}_0$, every linear order is complete if and only if it is compact \cite[Corollary 4.2]{Shafer}. This means that being a complete linear order is also a $\Pi^1_1$ universal formula over $\textbf{ACA}_0$.
\end{obs}

\section{Linear orders}
It is tempting to say that since we showed that over $\textbf{ACA}_0$ every $T_3$ space is homeomorphic to a subset of the rationals, we have proved that over $\textbf{ACA}_0$ every $T_3$ space is linearly ordered. The issue is that, in general, given $X$ a subspace of a linear order $(Y,<)$,  the subspace topology on $X$ is not the same as the topology given by the order $<|_{X\times X}$. For example, the set $S=\{0\}\cup\{1+\frac{1}{n}:n\in\mathbb{N}\}\subseteq \mathbb{Q}$ with the subspace topology will be a countable discrete space while with the order topology it will be homeomorphic to the $1$ point compactification of a countable discrete space. The issue here is that $\{0\}=S\cap\mathopen]-\frac{1}{2},\frac{1}{2}\mathclose[$ but there isn't a $b\in S$ such that $\mathopen]-\infty,b\mathclose[\cap S=\{0\}$ and so $\{0\}$ is not open in the order topology of $S$. In general topology, the subspaces of linearly ordered spaces with the order topology are called Generalized Ordered spaces or GO spaces. There are GO spaces that are not orderable, for example, the space $\left] 0,1\right[\cup\{2\}\subseteq \mathbb{R}$. However, in the countable case, we have that all GO spaces are orderable. 
\\
\begin{obs}\label{Path} Let $(L,<_L)$ be a linear order and $S\subseteq L$. Then the order topology on $S$ does not coincide with the subspace topology if and only if there is an $x\in S$ such that either:
\begin{enumerate}[label={(\arabic*)}]
\item There is a $y\in L\setminus S$ such that $y>x$, $\mathopen]x,y\mathclose[\subseteq L\setminus S$, and $\{z\in S: x<_L z\}$ does not have a least element.
\item There is a $y\in L\setminus S$ such that $y<x$, $\mathopen]y,x\mathclose[\subseteq L\setminus S$, and $\{z\in S: z<_L x\}$ does not have a greatest element.
\end{enumerate}
In general, the subspace topology on $S$ is finer than the topology induced by $<_L$. In the case where $L$ is a well-order, the order topology on $S$ coincides with the subspace topology  if and only if  every $x$ which is isolated in the subspace topology is also isolated in the order topology. Equivalently, we have that the order topology of $S$ coincides with the subspace topology if and only if every $s\in S$ which has a predecessor in $(L,<_L)$ has a predecessor in $(S,<_L)$.
\end{obs}

\begin{thm}\textbf{(Lynn \cite{Lynn})} Any separable zero dimensional metric space is orderable.
\end{thm}
\noindent
We modify Lynn's proof so that it can be carried out in $\textbf{ACA}_0$.

\begin{thm}\label{Lynns} $\textbf{ACA}_0$ proves that any $T_3$ space is homeomorphic to a linear order with its order topology.\\
\\
\begin{proof}
Let $X=(x_i)_{i\in\mathbb{N}}$ be a $T_3$ CSC space. By Theorem~\hyperref[stuff>aca]{\ref{stuff>aca}} every $T_3$ space is homeomorphic to a subset of the rationals, so without loss of generality, we may assume $X$ is a subspace of $\mathopen]0,\sqrt{2}\mathclose[$. For each $\sigma\in 2^{<\mathbb{N}}$ we define $I_\sigma$ to be the clopen subinterval of $\mathbb{Q}$
\begin{equation*}
\left]\sqrt{2} \left(\sum_{i<|\sigma|}\frac{\sigma(i)}{2^{i+1}}\right),\sqrt{2}\left(\sum_{i<|\sigma|}\frac{\sigma(i)}{2^{i+1}}+\frac{1}{2^{|\sigma|}}\right)\right[
\end{equation*}
By arithmetic comprehension the sequence $(I_\sigma)_{\sigma\in 2^{<\mathbb{N}}}$ exists. 
We observe that $\sigma\sqsubseteq \tau\leftrightarrow I_\tau\subseteq I_\sigma$ and that for all $n$ the family $\{I_\sigma: |\sigma|=n\}$ is a partition of $\mathopen]0,\sqrt{2}\mathclose[$ into $2^n$ many clopen sets. Given an $x\in X$ and $n\in \mathbb{N}$ we write $
seq_n(x)$ to be the unique sequence of length $n$ such that $x\in I_{seq_n(x)}$. \\
\\
The set $T=\{\sigma\in 2^{<\mathbb{N}}:\exists x\in X\;x\in I_\sigma\}$ exists by arithmetic comprehension. Inductively we define a partial function $a_{(\cdot)}:2^{<\mathbb{N}}\rightarrow X$. We set $a_{(0)}=x_0$, let $i\leq 1$ be such that $x_0\in I_{(i)}$ and if there is an element in $I_{(1-i)}\cap X$ then let $a_{(1)}$ be the least element of $I_{(1-i)}\cap X$.\\
\\ 
Assume $a_\sigma$ is defined. Let $i=\sigma(|\sigma|-1)$, that is $i$ is the last digit in the sequence $\sigma$, we define $a_{\sigma^\frown(i)}=a_\sigma$. For ease of notation, let
\begin{equation*} \tau=seq_n(a_\sigma)=seq_n(a_{\sigma^\frown(i)})\quad \text{ and }\quad j=seq_{n+1}(a_\sigma)(n)
\end{equation*}
If $\tau^\frown(1-j)\in T$ let $a_{\sigma^\frown(1-i)}$ be the least element of $X\cap I_{\tau^\frown(1-j)}$.\\
\\
We note that for a sequence $\sigma$ in the domain of $a_{(\cdot)}$ if $i=\sigma(|\sigma|-1)$ then $a_\sigma= a_{\sigma^\frown(i)}=\dots= a_{\sigma^\frown (i,\dots,i)}$.
\\
\\
By construction we have that for all $\sigma \in T$ exists a unique $\tau$ such that $|\tau|=|\sigma|$ and $a_\tau\in I_\sigma$, let $f:T\rightarrow 2^\mathbb{N}$ denote the map that sends $\sigma$ to $\tau$. We have $f(seq_n(a_\sigma))=\sigma$ and $f$ is total, order preserving, level preserving, and injective. The function $f$ is arithmetically definable relative to $X$, so it exists by arithmetic comprehension.
\\
\\
We define an order $<_f$ where
\begin{equation*} a<_f b\leftrightarrow \exists n\in\mathbb{N}\;(f(seq_n(a))<_{lex}f(seq_n(b)))
\end{equation*}
It is easy to verify that $<_f$ defines a total order on $X$.  Under this order, we have
\begin{equation*} I_{f^{-1}(\sigma)}=I_{seq_{|\sigma|}(a_{\sigma})}=[a_{\sigma^\frown(0)},a_{\sigma^\frown(1)}]_{<_f}
\end{equation*}
We show that the order topology induced by $<_f$ on $X$ is the same topology as the subspace topology. Given $x\in\mathopen]a,b\mathclose[_{<_f}$ we have that by definition there exists an $n$ such that
\begin{equation*}
f(seq_n(a))<_{lex}f(seq_n(x))<_{lex}f(seq_n(b))
\end{equation*}
So $I_{seq_n(x)}\subseteq \mathopen]a,b\mathclose[_{<_f}$ which implies that the topology on $X$ is finer than the topology induced by $<_f$.\\
\\
We now show the converse. We have that the set $\{I_{\sigma}\cap X:\sigma\in T\}$ generates the subspace topology on $X$, so it suffices to prove that they are open in the topology generated by $<_f$. Given $x\in I_\sigma$, let $\tau_0$ and $\tau_1$ be such that $|\tau_0|=|\sigma|=|\tau_1|$ and $f(\tau_0)$ is the lexicographic predecessor of $f(\sigma)$ and $f(\tau_1)$ the lexicographic successor of $f(\sigma)$ in the set $\{f(\rho):\rho\in T\wedge |\rho|=|\sigma|\}$, the case where $\sigma$ is the lexicographic maximum or minimum are proved similarly. By definition of $<_f$ we have that
\begin{equation*} I_{\tau_0}=[a_{f(\tau_0)^\frown (0)},a_{f(\tau_0)^\frown (1)}]_{<f}\quad \text{ and }\quad I_{\tau_1}=[a_{f(\tau_1)^\frown (0)},a_{f(\tau_1)^\frown (1)}]_{<f}
\end{equation*}
So in particular $\min I_{\tau_1} =a_{f(\tau_1)^\frown(0)}$ and $\max I_{\tau_0}=a_{f(\tau_0)^\frown(1)}$. By definition of $<_f$ we have \begin{equation*}I_\sigma\subseteq\mathopen]a_{f(\tau_0)^\frown(1)},a_{f(\tau_1)^\frown(0)}\mathclose[_{<_f}
\end{equation*} Let $x\in \mathopen]a_{f(\tau_0)^\frown(1)},a_{f(\tau_1)^\frown(0)}\mathclose[_{<_f}$, setting $n=|\sigma|$ we have
\begin{equation*}f(\tau_0)<_{lex} f(seq_n(x))<_{lex}f(\tau_1)
\end{equation*}
Since $f(\sigma)$ is the only sequence of length $n$ that is lexicographically between $f(\tau_0)$ and $f(\tau_1)$ we have that $f(seq_n(x))=f(\sigma)$ which by injectivity of $f$ means that $x\in I_\sigma$. So $I_\sigma=\mathopen]a_{f(\tau_0)^\frown(1)},a_{f(\tau_1)^\frown(0)}\mathclose[_{<_f}$ which implies the topology induced by $<_f$ is the same as the subspace topology of $X$.
\end{proof}
\end{thm}
\noindent
For compact $T_2$ CSC space, we can do better. Sierpi{\'n}ski and Masurkiewicz \cite{Sierp3} showed that any closed countable subset of $\mathbb{R}^m$ is well orderable. Friedman and Hirst essentially showed that over $\textbf{ATR}_0$ any compact $T_2$ CSC space is homeomorphic to a well-order with the order topology \cite[Lemma 4.7]{First}. We will show that this characterization can be carried out in $\textbf{ACA}_0$ and for $eT_2$ strongly effectively compact CSC spaces it can be carried out over $\textbf{RCA}_0+\text{WF}(\omega^\omega)$. The usual proof of this result is by a back and forth argument which requires considerations over the Cantor-Bendixson analysis of the space, which, as we will see later, requires $\mathbf{ATR}_0$. Furthermore, we know that any compact $T_2$ CSC space $X$ will be homeomorphic to $(\omega^\alpha+1)\cdot n$, where $\alpha$ is some well order and $n\in \omega$, but over $\mathbf{RCA}_0$ the exponential of a well order may not be a well order (see \cite[Theorem II.5.4.1]{Girard} and \cite[Theorem 2.6]{Hirst}). So over $\mathbf{RCA}_0$ we cannot simply formalize the standard argument as in $\mathbf{ATR}_0$.
\begin{defin}
    By $\text{WF}(\omega^\omega)$ we mean the statement that the lexicographic order on the non increasing sequences of $\mathbb{N}^{<\mathbb{N}}$ is well ordered. $\mathbf{RCA}_0+\mathbf{I}\Sigma_2^0\vdash \text{WF}(\omega^\omega)$ and $\mathbf{RCA}_0+\text{WF}(\omega^\omega)\vdash Con(\mathbf{RCA}_0)$, so in particular, $\text{WF}(\omega^\omega)$ is not provable in $\mathbf{RCA}_0$. Over $\mathbf{RCA}_0$ the statement that for every $n\in \mathbb{N}$ the Kleene-Brouwer order on the tree of all sequences of length at most $n$ is well founded is equivalent to $\text{WF}(\omega^\omega)$.
\end{defin}
\begin{obs}\label{incdec} Over $\textbf{RCA}_0$ let $(X,(U_i)_{i\in\mathbb{N}},k)$ be a compact effectively zero dimensional CSC space with covering relation $C$ and $G$ the function witnessing the effective zero dimensionality of $X$. We have that the inclusion relation on indices of basic open sets $\{(i,j): U_i\subseteq U_j\}$ is defined by
\begin{equation*} \forall x\in X\;(x\in U_i\rightarrow x\in U_j)
\end{equation*}
which is a $\Pi^0_1$ formula, and by
\begin{equation*} \exists F\subseteq X\setminus U_i\;(F\text{ is finite}\wedge\{G(x,i): x\in F\}\cup \{j\}\in C)
\end{equation*}
which is a $\Sigma^0_1$ formula. That is the inclusion relation $\{(i,j): U_i\subseteq U_j\}$ is $\Delta^0_1$ definable relative to $G$ and $C$ and therefore exists by $\Delta^0_1$ comprehension. If $X$ is an effectively compact space with an algebra of clopen sets, then we have
\begin{equation*}\{i\in\mathbb{N}: U_i\neq \emptyset\}=\{i: Comp(i)\notin C\}
\end{equation*}
that is, being an empty basic open subset is decidable relative to $Comp$ and $C$. Similarly, we have that for effectively compact CSC space with an algebra of clopen sets, strict inclusion will be decidable relative to the additional structure.
\end{obs}

\begin{thm}\label{compwell} 
\begin{enumerate}[label={(\arabic*)}]
\item $\mathbf{RCA}_0+\text{WF}(\omega^\omega)$ proves that if $X$ is $eT_2$ and strongly effectively compact then $X$ is effectively homeomorphic to the upper limit topology of a well-order with maximal element.
\item $\mathbf{RCA}_0$ proves that if $X$ is effectively homeomorphic to a well-order with a maximum element with the upper limit topology then $X$ is effectively compact and $eT_2$.
\end{enumerate}
\end{thm}
\noindent
\begin{proof} 
For the second implication, we observe that the upper limit topology on a well-order with maximum element is $eT_2$ and effectively compact and being $eT_2$, and effectively compact is preserved under effective homeomorphism.\\
\\
Let $(X,(U_i)_{i\in\mathbb{N}},k)$ be a strongly effectively compact $eT_2$ CSC space and let $C$ be the cover relation of $X$. The case in which $X=\emptyset$ is trivial, so we may assume $X\neq \emptyset$. By Theorem~\hyperref[Big4]{\ref{Big4}} and \hyperref[compact T3]{\ref{compact T3}} we may assume without loss of generality that $X$ has an algebra of clopen sets and for every $x\in X$ there are infinitely many $j\in \mathbb{N}$ such that $x\in U_j$. Let $Int$ and $Comp$ code the intersection and complement, respectively. As observed, we can effectively determine inclusion and being empty relative to the cover relation and the functions $Comp$ and $Int$.\\
\\
Let $F:X\times \mathbb{N}\rightarrow \mathbb{N}$ be a partial function such that for all $x$ in $X$ we have that
\begin{equation*} F(x,n)=\min\{s\in \mathbb{N}:
\forall m< n \;U_s\subseteq U_{F(x,m)}\}
\end{equation*}
We have that $F$ exists by $\Delta^0_1$ comprehension. Informally, $F$ lists out a weakly descending sequence of neighborhoods for every point. For ease of notation, we will write $U(x,n)=U_{F(x,n)}$. We observe that the sets of the form $U(x,n)$ form a basis of clopen sets for $X$.\\
\\
We construct inductively a tree $\mathcal{A}$ and we assign to each $\sigma \in \mathcal{A}$ a label $(x,i)\in X\times\mathbb{N}$ where $x$ is the $<_\mathbb{N}$ least element of $U_i$. For the first step, we add the empty sequence to $\mathcal{A}$ with label $(x,i)$, where $i$ is such that $X=U_i$ and $x$ is the $<_\mathbb{N}$ least element of $X$. Let $\sigma$ be a sequence in $\mathcal{A}$ with label $(y,j)$ and let $m$ be the least number such that such that $U(y,m)\subseteq U_j$. For each $k\geq m$, such that $U(x,k)\setminus U(x,k+1)$ is non empty, we add $\sigma^\frown(k+1)$ to $\mathcal{A}$ with label $(z,l)$ where $l$ is an index for $U(y,k)\setminus U(y,k+1)$ and $z$ is the $<_\mathbb{N}$ least element of $U_l$. If $U(y,m)\neq U_j$ then we add $\sigma^\frown (0)$ to $\mathcal{A}$ with label $(z,l)$ where $l$ an index for $U_j\setminus U(y,m)$ and $z$ is the $<_\mathbb{N}$ least element of $U_l$. We have that $\mathcal{A}$ is $\Delta^0_1$ definable relative to $(X,(U_i)_{i\in\mathbb{N}},k), C$, and $Comp$ so it exists by $\Delta^0_1$ comprehension.

\begin{obs} Given $\sigma,\tau \in \mathcal{A}$ with labels $(y,j)$ and $(z,l)$ respectively, we have:
\begin{enumerate}[label={(\arabic*)}]
\item $\tau\sqsubseteq \sigma$ if and only if $U_{l}\subseteq U_{j}$.
\item If $\tau$ properly extends $\sigma$ then $y\notin U_{l}$.
\item If $\tau$ properly extends $\sigma$ then $y<_\mathbb{N}z$.
\item If $\sigma\in \mathcal{A}$ is terminal then $y$ is isolated.
\item $\tau$ and $\sigma$ are incomparable if and only if $U_{j}$ and $ U_{l}$ are disjoint.
\item for all $w\in U_j\setminus \{y\}$ there exists exactly one $k\in\mathbb{N}$ such $\sigma^\frown (k)\in \mathcal{A}$, $\sigma^\frown (k)$ has label $(u,r)$ for some $u\in X$, and $w\in U_r$.
\end{enumerate} 
\end{obs}
\noindent
We prove that $\mathcal{A}$ is well-founded. Seeking a contradiction, assume that $\mathcal{A}$ is not well-founded. Let $f$ be a branch of $\mathcal{A}$ and let $(x_n,i_n)$ be the label on $f|_{<n}$. We have that $(U_{i_n})_{n\in\mathbb{N}}$ is a descending sequence of non empty clopen sets, so its intersection is also non empty otherwise $(X\setminus U_{i_n})_{n\in\mathbb{N}}$ will be a covering of $X$ without a finite subcovering. Let $x$ be an element of $\bigcap_{n\in\mathbb{N}}U_{i_n}$, then for all $n\in \mathbb{N}$ we have that $x_n<_\mathbb{N}x$. This is absurd since the sequence $(x_n)_{n\in\mathbb{N}}$ is strictly increasing and therefore is unbounded in $\mathbb{N}$. We note that since we are working over $\textbf{RCA}_0$, we cannot conclude that the Kleene-Brouwer order on $\mathcal{A}$ is a well-order. However, we will later use the strong compactness of $X$ to show that $\mathcal{A}$ is in fact well-ordered by the Kleene-Brouwer order.\\
\\
We show that every $x\in X$ appears in the label of exactly one sequence in $\mathcal{A}$. By the observations above, we have that $x$ cannot appear in the label of more than one sequence in $\mathcal{A}$. Given an $x\in X$, define
\begin{equation*}
S=\{\sigma\in \mathcal{A}: (y,j) \text{ is the label on } \sigma \text{ and }x\in U_{j}\}
\end{equation*}
$S\neq \emptyset$ since it contains the empty sequence. By the observations made above, $S$ is a chain, and since $\mathcal{A}$ is well-founded, we have $S$ must have a maximal element $\sigma$ with respect to inclusion. Let $(y,j)$ be the label on $\sigma$. If $x\neq y$ then by construction of $\mathcal{A}$ there will be an extension $\tau$ of $\sigma$ with label $(z,l)$ such that $x\in U_l$ contradicting the maximality of $\sigma$.\\
\\
We show that $X$ is homeomorphic to $\mathcal{A}$ with the upper limit topology induced by the Kleene-Brouwer order. Let $\alpha_x$ denote the unique sequence of $\mathcal{A}$ that has $x$ in its label. We show the map $G:X\rightarrow \mathcal{A}$ given by $x\mapsto \alpha_x$ is a homeomorphism between $X$ and $\mathcal{A}$ with the upper limit topology induced by the Kleene-Brouwer order.\\
\\
We show that $G$ is effectively continuous. Given $x\in X$ then a basic neighborhood for $\alpha_x$ is of the form $\mathopen] \alpha_y,\alpha_x\mathclose]$ with $\alpha_y<_{\text{KB}} \alpha_x$. Let $(x,j)$ be the label of $\alpha_x$ and $m$ the least number such that $U(x,m)\subseteq U_j$. 
Define
\begin{equation*}
l=\begin{cases} \min\{ n\geq m: y\notin U(x,n)\} \quad\text{ if } \quad \alpha_x\sqsubseteq \alpha_y \\
m \quad \quad\quad\quad\quad\;\, \text{ otherwise }
\end{cases}
\end{equation*} 
If $y\in U(x,m)$ then $\alpha_x\sqsubseteq \alpha_y$ and $y\in U(x,l-1)\setminus U(x,l)$ so $\alpha^\frown(l)\in \mathcal{A}$. We have that $\alpha_y\leq_\text{KB} \alpha_x^\frown (l)<_\text{KB} \alpha_x$ 
by construction and so $G(U(x,l))\subseteq  \mathopen]\alpha_y,\alpha_x \mathclose]$. If $y\notin U(x,m)$ then for all $z\in U(x,m)$  we have
\begin{equation*} \alpha_y <_\text{KB} \alpha_x^\frown(m)\leq_{\text{KB}}\alpha_z<_\text{KB} \alpha_x
\end{equation*}
and so $G(U(x,m))=G(U(x,l))\subseteq \mathopen]\alpha_y,\alpha_x\mathclose]$ (note that $\alpha_x^\frown (m)$ may not be a member of $\mathcal{A}$). The function  $(x,y)\mapsto l$ exists by $\Delta^0_1$ comprehension and witnesses that $G$ is effectively continuous.
\\
\\
We show that $G$ is effectively open. Let $x\in X$ be a point and $U_i$ a basic neighborhood of $x$. Let $(x,j)$ be the label on $\alpha_x$, $m$ the least number such that $U(x,m)\subseteq U_j\cap U_i$ and $m_0$ the least number such that $U(x,m_0)\subseteq U_j$. We consider the following cases:\\
\\
\textbf{Case 1:} If $U(x,m)\subsetneq U_j$ and $U(x,m)=U(x,m_0)$, let $y$ be the least element in $U_j\setminus U(x,m)$. We have that $\alpha_y=\alpha_x^\frown (0)$ and for any $z\in X$ if $\alpha_y<_{KB}\alpha_z\leq_{\text{KB}} \alpha_x$ then $z\in U(x,m)\subseteq U_i$. We set $I=\mathopen]\alpha_y,\alpha_x\mathclose]$.\\
\\
\textbf{Case 2:} If $U(x,m)\subsetneq U_j$ and $U(x,m)\neq U(x,m_0)$, set $t=\max\{s< m: U(x,s)\setminus U(x,s+1)\neq \emptyset\} $ and $y$ the $<_\mathbb{N}$ least element of $U(x,t)\setminus U(x,t+1)$. We have that $\alpha_y=\alpha_x^\frown (t+1)$ and for any $z\in X$ if $\alpha_y<_{KB}\alpha_z\leq_{\text{KB}} \alpha_x$ then $z\in U(x,t+1)\subseteq U(x,m)\subseteq U_i$. We set $I=\mathopen]\alpha_y,\alpha_x\mathclose]$.\\
\\
\textbf{Case 3:} If $U(x,m)=U_j$ and there are no sequences to the left of $\alpha_x$ we set $I=\mathopen]-\infty,\alpha_x\mathclose]= G(U_j)\subseteq G(U_i)$.\\
\\
\textbf{Case 4:} If $U(x,m)=U_j$ and there are sequences to the left of $\alpha_x$ in $\mathcal{A}$ then define
\begin{equation*}
s=\max\{ t<|\alpha_x|: \exists v<\alpha_x(t)\;(\alpha_x|_{<t}^\frown (v)\in \mathcal{A})\}
\end{equation*}
\begin{equation*}
u=\max\{ v<\alpha_x(s): \alpha_x|_{<s}^\frown (v)\in \mathcal{A}\}
\end{equation*}
We have that $\alpha_x|_{<s}^\frown (u)\in\mathcal{A}$ will be equal to some $\alpha_y$ and for all $z\in U(x,m)$ we have $\alpha_y<_{KB} \alpha_z\leq_{KB} \alpha_x$ and so let $I=\mathopen]\alpha_y,\alpha_x\mathclose]= G(U(x,m))$.
In all four cases we have that $I\subseteq G(U(x,m))$, we can also effectively find an index for $I$ and so $G$ is effectively open.\\
\\
This proves $G$ is an effective homeomorphism from $X$ to the upper limit topology of $(\mathcal{A},<_\text{KB})$. We now prove that $\mathcal{A}$ is a well-order. Assume that there exists an infinite sequence $(x(n))_{n\in\mathbb{N}}$ such that $(\alpha_{x(n)})_{n\in\mathbb{N}}$ is strictly $<_{\text{KB}}$-decreasing in $\mathcal{A}$. Without loss of generality, we may assume that:
\begin{enumerate}[label={(\arabic*)}]
\item The sequences $(\alpha_{x(n)})_{n\in\mathbb{N}}$ are pairwise incomparable. We can do this since we showed that $\mathcal{A}$ is well-founded.
\item For all $n\in \mathbb{N}$ we have $|\alpha_{x(n)}|<|\alpha_{x(n+1)}|$. Here we are using the assumption $\text{WF}(\omega^\omega)$ in the form for all $m$ the Kleene-Brouwer order of the tree of all sequences of length at most $m$ is a well-order.
\end{enumerate} 
Let $(x(n),j(n))$ be the label on $\alpha_{x(n)}$ we would like to show that $O=\bigcup_{n\in\mathbb{N}} U_{j(n)}$ exists as a set and is closed. Since the $(\alpha_{x(n)})_{n\in\mathbb{N}}$ are pairwise incomparable we have that the sets $U_{j(n)}$ are pairwise disjoint. We have that for all $y\in X$ that $y\in O$ if and only if there exists an $n\in\mathbb{N}$ such that $\alpha_{x(n)}\sqsubseteq \alpha_y$. Since the lengths of the $\alpha_{x(n)}$ are strictly increasing we have that
\begin{equation*} 
y\in O\leftrightarrow \exists n<|\alpha_y|\:(\alpha_{x(n)}\sqsubseteq \alpha_y)
\end{equation*}
In particular, we have that $O$ exists by $\Delta^0_1$ comprehension. We show that $O$ is closed. Given $y\notin O$, let $(y,l)$ be the label on $\alpha_y$. If $U_l\cap O=\emptyset$ we are done, otherwise there exists $x(n)\in U_l$ and $\alpha_{x(n)}<_\text{KB}\alpha_y$. Let $n=\min\{j\in\mathbb{N}: x(j)\in U_l\}$ and let $m=\min\{j\in\mathbb{N}: x(n)\notin U(y,j)\}$. We have that $U(y,m)\cap O=\emptyset$, so $O$ is closed. The sequence $\{U_{j(n)}:n\in\mathbb{N}\}\cup\{ X\setminus O\}$ is an infinite partition of $X$ in clopen sets contradicting our assumption that $X$ was strongly compact.
\end{proof}

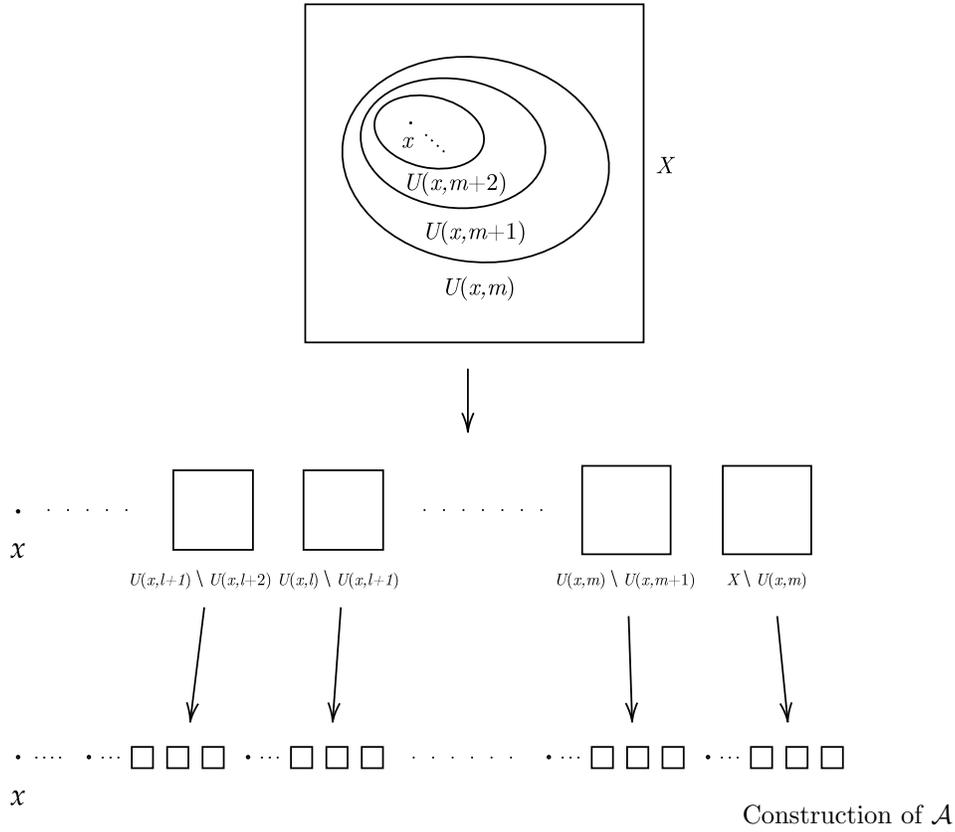
\[
\centering
\begin{tikzpicture}
\tikzstyle{every node}=[font=\normalsize]
\draw [ line width=1pt ] (8.75,7) rectangle (8.75,7);
\draw [ line width=1pt ] (4.5,7) rectangle (4.5,7);
\draw [ line width=1pt ] (3.75,10.25) rectangle (8.75,6.25);
\draw [ line width=1pt ] (6.25,6.5) rectangle (6.25,6.5);
\node [font=\normalsize] at (9.25,8) {$X$};
\draw [ line width=1pt ] (6,8.5) ellipse (2cm and 1.5cm);
\draw [ line width=1pt ] (5.75,8.75) ellipse (1.25cm and 1cm);
\draw [ line width=1pt ] (5.5,8.75) circle (0cm);
\node [font=\normalsize] at (6.75,6.75) {$U(x,m)$};
\node [font=\normalsize] at (6,7.5) {\textit{U(x,m+1)}};
\draw [line width=1pt, ->, >=Stealth] (6.25,5.75) -- (6.25,4.75);
\draw [line width=1pt, ->, >=Stealth] (7.75,5.75) -- (7.75,4.75);
\node [font=\normalsize] at (5.5,9.2) {\textit{x}};
\draw [ line width=1pt ] (8.25,4.25) rectangle (9.5,3);
\draw [ line width=1pt ] (5.5,4.25) rectangle (6.75,3);
\draw [line width=1pt, dotted] (4.1,3.5) -- (3.4,3.5);
\node [font=\normalsize] at (9,2.5) {$X\setminus  U(x,m)$};
\node [font=\normalsize] at (6,2.5) {$U(x,m)\setminus U(x,m+1)$};
\node [font=\normalsize] at (2,2.8) {$x$};
\draw [ line width=1pt ] (4.75,9.25) ellipse (0cm and 0cm);
\draw [line width=1pt, dotted] (6.2,8.42) -- (5.67,8.75);
\filldraw[black] (5.5,8.9) circle (1pt);
\filldraw[black] (2,3.5) circle (1pt);
\draw [line width=1pt, ->, >=Stealth] (4.75,5.75) -- (4.75,4.75);
\end{tikzpicture}
\]
Construction of $\mathcal{A}$.
\\
\\
\begin{cor}\label{corwell} $\textbf{ACA}_0$ proves that every sequence of compact $T_2$ CSC spaces  $(K_i)_{i\in\mathbb{N}}$ the disjoint sum $\coprod_{i\in\mathbb{N}} K_i$ is homeomorphic to a well-order with the order or the upper limit topology.\\
\\
\begin{proof} Over $\textbf{ACA}_0$, the upper limit topology on a well-order is homeomorphic to the order topology. Following the proof of Theorem~\hyperref[compwell]{\ref{compwell}} we can uniformly define a sequence of well-orders $(<_i)_{i\in\mathbb{N}}$ such that $<_i$ has $K_i$ as its field and the order topology induced by $<_i$ is the same as the topology of $K_i$. Let $<_K$ be the order on $\coprod_{i\in\mathbb{N}}K_i$ given by
\begin{equation*}
(x,j)<_K(y,i)\leftrightarrow (j<_{\mathbb{N}} i \vee (j=i\wedge x<_i y))
\end{equation*}
It is straightforward to show that $<_K$ is a well-order and its order topology is the same as the topology of $\coprod_{i\in\mathbb{N}}K_i$.

\end{proof}
\end{cor}

\section{Summary Part I}
\begin{thm}
Over $\textbf{RCA}_0$ the following are equivalent:
\begin{enumerate}[label={(\arabic*)}]
\item Arithmetic comprehension.
\item Every compact $T_2$ CSC space is $eT_2$ (Dorais \cite[Example 7.4]{Dorais}).
\item Every well-order with the upper limit topology is effectively homeomorphic to the order topology of some linear order (Proposition ~\hyperref[wellmess]{\ref{wellmess}}).
\item Every $eT_2$ scattered $T_3$ space is $uT_3$ (Proposition~\hyperref[worst proof]{\ref{worst proof}} and Observation~\hyperref[scatter metriz>aca]{\ref{scatter metriz>aca}}).
\item Every $T_3$ CSC space is effectively homeomorphic to a linear order with the order topology (Theorem~\hyperref[Lynns]{\ref{Lynns}}).
\item Every compact CSC space is effectively compact (Dorais \cite[Example 3.5]{Dorais}).
\item Every $uT_3$ compact CSC space is effectively compact (Theorem~\hyperref[effcomp<>ACA]{\ref{effcomp<>ACA}}).
\item Every $T_2$ compact CSC space is effectively homeomorphic to a well-order with the upper limit topology or the order topology (Corollary~\hyperref[corwell]{\ref{corwell}}).
\item Every $eT_2$ $T_3$ scattered CSC space effectively embeds into a linear order (Theorem~\hyperref[theworst]{\ref{theworst}}).
\item  For every $T_3$ scattered CSC space there is a continuous bijection from $X$ to the order topology of a well-order (Theorem~\hyperref[horrible]{\ref{horrible}}).
\item Every $eT_2$ effectively compact CSC space is strongly compact (Proposition~\hyperref[strongcomp]{\ref{strongcomp}}).
\end{enumerate}
\end{thm}
\noindent
Along with the following results, one can produce many more statements equivalent to arithmetic comprehension over $\textbf{RCA}_0$ (Proposition~\hyperref[G>Alg]{\ref{G>Alg}} and Theorem~\hyperref[Big4]{\ref{Big4}}). $\textbf{RCA}_0$ proves that for any CSC space $X$ the following are equivalent:
\begin{enumerate}[label={(\arabic*)}]
\item $X$ is $uT_3$.
\item $X$ is metrizable.
\item $X$ is $T_0$ and effectively zero dimensional.
\item $X$ is $T_0$ effectively homeomorphic to a space with an algebra of clopen sets.
\item $X$ is homeomorphic to a subspace of $\mathbb{Q}$.
\item $X$ is homeomorphic to a closed subspace of $\mathbb{Q}$.
\end{enumerate}
and by Theorem~\hyperref[noT0]{\ref{noT0}} we have that the following are equivalent:
\begin{enumerate}[label={(\arabic*)}]
\item $X$ is uniformly regular.
\item $X$ is effectively zero dimensional.
\item $X$ is effectively homeomorphic to a CSC space with an algebra of clopen sets.
\item $X$ is pseudometrizable.
\end{enumerate}

\newpage
\part{Topological characterizations of $\mathbf{ATR}_0$ }
 
\noindent
We saw that $T_3$ CSC spaces admit nice characterizations. Namely, they are all orderable, metrizable, and homeomorphic to a subspace of the rationals, and that these characterizations can be carried out in $\textbf{ACA}_0$.  We will restrict our attention to locally compact CSC spaces and $T_3$ scattered CSC spaces. Countable scattered spaces have been studied in general topology since every second countable scattered space is also countable (see \cite[Proposition 8.5.5]{Semadeni}). Kat{\v e}tov \cite{Katetov} showed that the $T_3$ scattered CSC spaces are precisely those which admit a continuous bijection to a compact $T_2$ CSC space. Bel'nov \cite{Belnov} and Knaster and Urbanik \cite[Theorem IV]{Knaster} showed that every scattered subspace of the Cantor space is homeomorphic to a subspace of a well order.  We will show that both these theorems, formulated in terms of CSC spaces, are equivalent to arithmetic transfinite recursion over $\mathbf{RCA}_0$. The $T_3$ scattered CSC spaces are precisely those which are completely
metrizable\footnote{This result is also attributed to Bel'nov and Knaster, Urbanik by Kannan and Rajagopalan \cite{KannanRajagop}. However, Knaster and Urbanik do not explicitly mention complete metrics and instead their results are stated in terms of being homeomorphic to a closed subspace of the Cantor space.  They in fact cite Young \cite[Page 65]{Young} for the result that every CSC space is $G_\delta$ if and only if it is scattered. 
The theorem that every $G_\delta$ set of a Polish space is Polish was in part shown by Hausdorff \cite{Hausdorff}
and the equivalence is stated in a footnote of Alexandroff and Urysohn \cite[Page 96]{Alexandrov}. 
These two theorems combined give the result that a CSC space is completely metrizable if and only if it is $T_3$ and scattered. This result, with time, seems to have become folklore (see \cite{Michael} the remark after Theorem 1.3).},
we will show that this characterization theorem is also equivalent to arithmetic transfinite recursion over $\mathbf{RCA}_0$.  We will also consider the $T_2$ locally compact CSC spaces which are a strict subclass of the $T_3$ scattered CSC spaces. Mazerkiewicz and Sierpi{\'n}ski \cite{Sierp3} essentially showed that $T_2$ locally compact CSC spaces are precisely those which are homeomorphic to a countable well order. Unlike for $T_2$ compact CSC spaces, which can be well ordered in the system $\mathbf{ACA}_0$, the well orderability of locally compact $T_2$ CSC spaces is equivalent to arithmetic transfinite recursion. We will also find a series of related characterization theorems, all of which will be equivalent to arithmetic transfinite recursion. Finally we will show that over $\mathbf{ACA}_0$ several properties of CSC spaces are $\Pi^1_1$ universal.

\section{Locally compact and scattered spaces}

The study of $T_2$ locally compact CSC spaces was done indirectly by Hirst, who considered locally totally bounded closed sets of complete separable metric spaces \cite{Hirst}. He proved that over $\textbf{ACA}_0$ any such space is the disjoint union of compact balls and used this decomposition to define the one point compactification. In this section, we will reformulate these theorems for CSC spaces and we will give a definition for the one point compactification which is suitable for $\textbf{RCA}_0$. 

\begin{defin} A topological space is said to be dense in itself if it does not contain isolated points. By Corollary~\hyperref[Qchar]{\ref{Qchar}} $\mathbf{ACA}_0$ proves that $\mathbb{Q}$ is the only non empty dense in itself $T_3$ CSC space up to homeomorphism.
\end{defin}

\begin{defin}
A CSC space is said to be scattered if every subspace has an isolated point or, rather, it does not contain a non empty dense in itself subspace. We observe that over $\textbf{ACA}_0$, being a  $T_3$ scattered CSC space is equivalent to not having a subspace homeomorphic to $\mathbb{Q}$. Furthermore, over $\textbf{RCA}_0$, we have that a uniformly $T_3$ CSC space is scattered if and only if it does not have a subspace which is effectively homeomorphic to $\mathbb{Q}$. 
\end{defin}

\begin{defin} Let $(X,(U_i)_{i\in\mathbb{N}},k)$ be a CSC space. We say that $X$ is locally compact (or l.c.\ for short) if, for every $x\in X$, there exists a compact neighborhood of $x$. We say that $X$ is locally effectively compact if, for every $x\in X$, there exists an effectively compact neighborhood for $x$. We say that a CSC space $X$ has a choice of compact neighborhoods or $CCN$ for short, if there exists a sequence $(K_x,i(x))_{x\in X}$ such that for all $x\in X$ we have $x\in U_{i(x)}\subseteq K_x$ and $K_x$ is compact. We say that a CSC space $X$ has a choice of effective compact neighborhoods or an effective $CCN$ if there is a sequence $(K_x, i(x), C_x)_{x\in X}$ such that $x\in U_{i(x)}\subseteq K_x$, $K_x$ is compact, and $C_x$ is a covering relation for $K_x$.
\end{defin}

\begin{obs} In the case in which $(X,(U_i)_{i\in \mathbb{N}},k)$ is l.c.\ with a basis of clopen sets, then a choice of compact neighborhoods is equivalent to having a sequence of indices $(i(x))_{x\in X} $ such that $U_{i(x)}$ is compact. Since we will mostly be working with zero dimensional spaces we will usually use this definition of $CCN$.
\end{obs}

\begin{lemma}\label{lc are scattered} $\textbf{ACA}_0$ proves that every $T_2$ l.c.\ CSC space is scattered.\\
\\
\begin{proof} Let $(X,(U_i)_{i\in\mathbb{N}})$ be a $T_2$ l.c.\ CSC space, and let $S$ be a non empty dense in itself subset of $X$. Enumerate $X=(x_n)_{n\in\mathbb{N}}$, let $y_0\in S$, and let $V_0$ be a compact neighborhood of $y_0$. Since we are working in $\textbf{ACA}_0$, we have that $V_0$ is also sequentially compact. Given $y_n$ and $V_n$ such that $y_n\in V_n$, let $y_{n+1}$ be the least element in $S\cap V_n\setminus\{y_n\}$ which will be infinite since $S$ does not have isolated points. Let $V_{n+1}$ be the first basic open set such  that $y_{n+1}\in V_{n+1}\subseteq V_n$ and $x_{n+1}\notin V_{n+1}$. The sequence $(y_n,V_n)_{n\in\mathbb{N}}$ exists by arithmetic comprehension. Assume that the subsequence $(y_{n_k})_{k\in\mathbb{N}}$ converges to $x_m\in X$,  by construction $(y_{n_k})_{k\in\mathbb{N}}$ is eventually in $V_m$ which does not contain $x_m$ by construction. This contradicts the fact that $V_0$ is sequentially compact. 
\end{proof}
\end{lemma}

\begin{obs} The converse does not hold. For a counterexample, take the set
\begin{equation*}S=\omega^2+1\setminus\{\omega\cdot n:n\in\omega\}\subseteq \omega^2+1
\end{equation*} with the subspace topology. We observe that the point $\omega^2\in S$ does not have a compact neighborhood in $S$. Since $S$ is order isomorphic to $\omega^2+1$, its subspace topology cannot coincide with its order topology. In general, the subspaces of a l.c.\ CSC space may not be l.c.\ This counterexample can be easily formalized in $\textbf{RCA}_0$. We will see later on that all scattered $T_3$ CSC spaces, over a sufficiently strong theory, are precisely the subspaces of well orders. It is also worth noting that, in general, l.c.\ $T_2$ second countable spaces are not scattered. For example, $\mathbb{R}$ is locally compact, but it is not scattered since it is connected.
\end{obs}

\begin{prop}\label{lc>T3} $\textbf{RCA}_0$ proves that every $eT_2$ locally effectively compact CSC space is $T_3$.\\
\\
\begin{proof} Let $(X,(U_i)_{i\in\mathbb{N}},k)$ be a $T_2$ locally effectively compact CSC space and let $x\in X$ and $i\in\mathbb{N}$ be such that $x\in U_i$. Let $K$ be an effectively compact neighborhood of $x$ and let $j\in\mathbb{N}$ be such that $x\in U_j\subseteq K$. By Proposition~\hyperref[Kcode]{\ref{Kcode}} there is an $f$ which codes $K$. We have that $K$ is effectively $T_2$ since it is the subspace of an effectively $T_2$ space. By Theorem~\hyperref[compact T3]{\ref{compact T3}} $K$ is uniformly $T_3$. Let $R^K_0$ and $R^K_1$ witness that $K$ is uniformly regular and let $s=k(x,i,j)$, we have
\begin{equation*}
x\in U_{R^K_0(x,s)}\subseteq K\setminus \bigcup_{y\notin U_{s}}U_{R^K_1(x,s,y)} \subseteq U_{s}\subseteq K\cap U_i
\end{equation*}
So, in particular
\begin{equation*}
x\in U_{R^K_0(x,s)}\subseteq X\setminus \left(\bigcup_{y\notin U_{s}}U_{R^K_1(x,s,y)}\cup \bigcup_{n\in\mathbb{N}}U_{f(n)}\right)\subseteq U_{s}\subseteq U_i
\end{equation*}
which shows that $X$ is regular.
\end{proof}
\end{prop}
\noindent
Using the same proof, we get
\begin{cor} $\textbf{RCA}_0$ proves that every $eT_2$ l.c.\ CSC space with an effective $CCN$ is uniformly $T_3$.
\end{cor}

\begin{cor} $\textbf{ACA}_0$ proves that every l.c.\ $T_2$ CSC space is $T_3$.
\end{cor}

\begin{obs} Let $(W,<_W)$ be a well order, then we have that the collection $(\mathopen]-\infty,w+1\mathclose[_{<_W})_{w\in W}$ is a $CCN$ for $W$ with the order topology.
\end{obs}

\begin{prop}\label{CCNstab} $\textbf{RCA}_0$ proves that having a $CCN$ and an effective $CCN$ are preserved under effective homeomorphism. 
\\
\\
\begin{proof} Let $(X,(U_i)_{i\in\mathbb{N}},k)$ and $(Y,(V_i)_{i\in\mathbb{N}},k')$ be a CSC spaces and $f:X\rightarrow Y$ a homeomorphism and let $v$ be the function witnessing that $f$ is effectively open. Assume that $(K_x,i(x))_{x\in X}$ is a $CCN$ for $X$. Then $(f(K_{f^{-1}(y)}),v(f^{-1}(y),i(f^{-1}(y))))_{y\in Y}$ is a $CCN$ for $Y$. If there exists a sequence $(C_x)_{x\in X}$ such that for all $x$ the set $C_x$ is the covering relation for $K_x$ then by Proposition~\hyperref[effhomcom]{\ref{effhomcom}} there exists a sequence $(D_y)_{y\in Y}$ such that for all $y\in Y$ $D_y$ is the covering relation for $f(K_{f^{-1}(y)})$.
\end{proof}
\end{prop}

\begin{obs} This definition of local compactness is quite weak. The standard definition of local compactness that is usually used in general topology is that for every point $x$ and open neighborhood $U$ of $x$, there is a compact neighborhood of $x$ contained in $U$. It is a classic result in general topology that both definitions of local compactness coincide for $T_2$ spaces.
\end{obs} 

\begin{prop}\label{compbasis} $\textbf{ACA}_0$ proves that any $T_2$ l.c.\ CSC space with $CCN$ is homeomorphic to a space with a basis of compact sets.\\
\\
\begin{proof}  Let $(X,(U_i)_{i\in \mathbb{N}},k)$ to be a $T_2$ CSC space with a $CCN$. By Propositions~\hyperref[lc>T3]{\ref{lc>T3}},  \hyperref[CCNstab]{\ref{CCNstab}}, and Theorem~\hyperref[stuff>aca]{\ref{stuff>aca}} we may assume that $X$ has a basis of clopen sets. Let $(i(x))_{x\in X}$ be a $CCN$ for $X$ and let $J=\{j\in\mathbb{N}:\exists x\in X\, U_j\subseteq U_{i(x)}\}$. We have that for all $j\in J$, the basic clopen set $U_j$ is compact by \hyperref[closedincomp]{\ref{closedincomp}}. We have that for all $n\in \mathbb{N}$ the set $U_n$ is open with respect to the basis $(U_j)_{j\in J}$ since for all $x\in U_n$ we have
\begin{equation*}
x\in U_{ k(x,i(x),n)}\subseteq U_n\cap U_{i(x)}\subseteq U_n
\end{equation*} and  $k(x,i(x),n)\in J$ since $U_{ k(x,i(x),n)}\subseteq U_{i(x)}$. So  $(U_j)_{j\in J}$ is a basis for $X$ consisting of compact sets.
\end{proof}
\end{prop}

\begin{prop}\label{2 for 1} $\textbf{ACA}_0$ proves that every l.c.\ $T_2$ CSC space with a $CCN$ is the disjoint union of open compact sets and is well-orderable.\\
\\
\begin{proof} Let $(X,(U_i)_{i\in\mathbb{N}},k)$ be a l.c.\ $T_2$ CSC space with a $CCN$. By \hyperref[compbasis]{\ref{compbasis}}, we may assume that $X$ has a basis of compact sets. The sequence of sets
\begin{equation*}
V_i=U_j\setminus \left(\bigcup_{i<j} U_i\right)
\end{equation*}
defines a partition of $X$ into clopen compact sets. In particular, we have that $X$ is homeomorphic to the disjoint sum of $(V_i)_{i\in\mathbb{N}}$. By Corollary~\hyperref[corwell]{\ref{corwell}} we have that $\coprod_{i\in \mathbb{N}}V_i$ is well-orderable and therefore $X$ is well-orderable.
\end{proof}
\end{prop}

\begin{obs} Since every well order has a $CCN$, we have that over $\textbf{ACA}_0$ a $T_2$ CSC space is well-orderable if and only if it has a $CCN$. This means that the strength of the well orderability of $T_2$ l.c.\ CSC spaces lies in being able to choose a compact neighborhood for every point of a $T_2$ l.c.\ CSC space. In the next section we will show that this is equivalent to arithmetic transfinite recursion over $\mathbf{RCA}_0$.
\end{obs}

\begin{prop}\label{in comp>CCN} $\textbf{ACA}_0$ proves that for any l.c.\ $T_2$ CSC space $(X,(U_i)_{i\in\mathbb{N}})$ any l.c. subspace $Y\subseteq X$ also has a $CCN$.\\
\\
\begin{proof} By Proposition~\hyperref[embT2]{\ref{embT2}} $Y$ is $T_2$ since it is the subspace of a $T_2$ space. Since $CCN$ are preserved under homeomorphism, we may assume without loss of generality that $X$ has a basis of clopen sets. Let $(i(x))_{x\in X}$ be a $CCN$ for $X$. For each $y\in Y$ and $j\in \mathbb{N}$ such that $U_j\cap Y$ is closed in $X$, if $U_j\subseteq U_{i(y)}$ then $U_j\cap Y$ is compact by Proposition~\hyperref[closedincomp]{\ref{closedincomp}}. Conversely, if $U_j\cap Y$ is compact, then by Proposition~\hyperref[Kcode]{\ref{Kcode}} $Y\cap U_j$ is closed in $X$ since it is a compact set in a $T_2$ space.  For each $y\in Y$ let $\widehat{i}(y)$ be the first $i$ such that  $y\in U_i\subseteq U_{i(x)}$ and $U_i\cap Y$ is closed in $X$. The sequence $(\widehat{i}(y))_{y\in Y}$ is arithmetically definable, so it exists by arithmetic comprehension and is a $CCN$ for $Y$.
\end{proof}
\end{prop}
\noindent
We will end this section by giving a definition of one point compactification which is suitable over $\textbf{RCA}_0$. 
\begin{defin} Let $X$ be a topological space, the Alexandroff compactification or one point compactification of $X$ is the space $X\cup\{\infty\}$ where $U\subseteq X\cup\{\infty\}$ is open if $\infty \notin U$ and $U$ is open in $X$ or $\infty \in U$ and $X\setminus U$ is compact in $X$. It is a classical result in general topology that the Alexandroff compactification of $X$ is $T_2$ if and only if $X$ is $T_2$ and l.c.\ 
\end{defin} 

\begin{obs} In general the one point compactification of a CSC space may not be a CSC space. For example, the one point compactification of $\mathbb{Q}$ is not  second countable. We would like to consider CSC spaces that have one point compactification which is also a CSC space.
\end{obs}

\begin{defin} A space $X$ is hemicompact if there exists an increasing sequence of compact sets $(K_n)_{n\in \omega}$ such that every compact set of $X$ is contained in some $K_n$. Equivalently, a space is hemicompact if the poset of compact subsets has a countable cofinal sequence. The Alexandroff or one point compactification of a space $X$ is first countable if and only if $X$ is first countable and hemicompact. \\
\\
A space $X$ has an exhaustion by compact sets if there is a sequence of compact sets $(K_n)_{n\in \mathbb{N}}$ such that for all $n$ $K_n\subseteq int(K_{n+1})$ and $X=\bigcup_{n\in\mathbb{N}}K_n$. We observe that a space $X$ with an exhaustion by compact sets is hemicompact and l.c.\
\end{defin}

\begin{defin} Over $\textbf{RCA}_0$ we say that a CSC space $X$ has an exhaustion by compact sets if there exists a sequence $(K_i,V_i,f_i,g_i)_{i\in\mathbb{N}}$ such that:
\begin{itemize}
\item $\bigcup_{i\in\mathbb{N}} K_i= X$
\item For each $i\in \mathbb{N}$ $K_i$ is compact, effectively closed, and $g_i$ is a code for $K_i$.
\item For each $i\in \mathbb{N}$ $V_i$ is an effectively open set and $f_i$ is a code for $V_i$. 
\item For each $i\in\mathbb{N}$ $K_i\subseteq V_{i+1}\subseteq K_{i+1}$.
\end{itemize}
\end{defin}

\begin{prop}\label{Chonk} $\textbf{RCA}_0$ proves that every $T_2$ locally effectively compact CSC space with an effective $CCN$ has an exhaustion of compact sets $(\widehat{K}_i,V_i,f_i,g_i)_{i\in\mathbb{N}}$. Furthermore there exists a sequence $(C_n)_{n\in\mathbb{N}}$ such that for all $n\in\mathbb{N}$ $C_n$ is the covering relation for $\widehat{K}_n$.\\
\\
\begin{proof}
Let $(X,(U_i)_{i\in \mathbb{N}})$ be a l.c.\ CSC  space and let $(K_x,i(x),C_x)_{x\in X}$ be an effective $CCN$. Enumerate the points of $X=(x_n)_{n\in\mathbb{N}}$. For each finite subset $a\subseteq X$ we have that $\widehat{C}_a=\bigcap_{i\in a} C_i$ is the covering relation for $\bigcup_{x\in a} K_x$. The sequence $(\widehat{C}_a)_{a\subseteq X}$ exists by $\Delta^0_1$ comprehension. Instead of defining directly the sequence $(\widehat{K}_i)_{i\in\mathbb{N}}$ we instead define recursively a sequence of finite sets $(F_n)_{n\in\mathbb{N}}$ such that $\widehat{K}_n=\bigcup_{x\in F_n} K_x$. Define $F_0=\{\min X\}$. Given $F_n$, define
\begin{equation*}
F_{n+1}=\min\{F\text{ finite set }:\{i(x):x\in F\}\in \widehat{C}_{F_n}\wedge\exists x\in F\;(x_n\in U_{i(x)})\}
\end{equation*} 
Equivalently $F_{n+1}$ is the least finite subset of $X$ such that $(U_{i(x)})_{x\in F_{n+1}}$ is a covering of $\bigcup_{y\in F_n} K_y\cup\{x_n\}$. We have that the sequence $(F_n)_{n\in\mathbb{N}}$ is well defined since for each $n$ $\bigcup_{y\in F_n} K_y$ is compact and it exists by $\Delta^0_1$ comprehension.\\
\\
Define $\widehat{K}_n=\bigcup_{x\in F_n} K_x$, we have that $\widehat{C}_{F_n}$ is the covering relation for $\widehat{K}_n$. By Proposition~\hyperref[Kcodes]{\ref{Kcodes}}, we have that the sequence $(\widehat{K}_n)_{n\in\mathbb{N}}$ will be uniformly effectively closed and so there exists a sequence $(g_n)_{n\in\mathbb{N}}$ such that for all $n\in\mathbb{N}$ $g_n$ is a closed code for $\widehat{K}_n$.\\
\\
Let $f_n$ be the partial function that enumerates $F_n$ increasingly. The sequences $(\widehat{K}_n,\bigcup_{x\in F_n}U_{i(x)},f_n,g_n)_{n\in\mathbb{N}}$ and $(C_{F_n})_{n\in\mathbb{N}}$ are $\Delta^0_1$ definable relative to $(F_n)_{n\in\mathbb{N}}$ and so they exist by $\Delta^0_1$ comprehension. By construction we have that $$(\widehat{K}_n,\bigcup_{x\in F_n}U_{i(x)},f_n,g_n)_{n\in\mathbb{N}}$$ is an exhaustion by compact sets of $X$ and for all $n\in\mathbb{N}$ $C_{F_n}$ is a covering relation for $\widehat{K}_n$.
\end{proof}
\end{prop}

\begin{defin}\label{defprop} Let $(X,(U_i)_{i\in\mathbb{N}},k) $ be a $T_2$ CSC space with an exhaustion of compact sets $(K_n,A_n,f_n,g_n)_{n\in \mathbb{N}}$ we define its one point compactification as the space
\begin{equation*}
(X\cup\{\infty\},(V_{i})_{i\in\mathbb{N}}\,\widehat{k})
\end{equation*} 
where $V_{2i}=U_i$ and $V_{2i+1}=(X\cup\{\infty\})\setminus K_i$.
\begin{enumerate}[label={(\arabic*)}]
\item $\widehat{k}(x,2i,2j)=k(x,i,j)$.
\item $\widehat{k}(x,2i,2j+1)=k(x,i,g(r))$ where $r=\min\{n\in \mathbb{N}: x\in U_{g(r)}\}$.
\item $\widehat{k}(x,2i+1,2j+1)= \max\{2i+1,2j+1\}$.
\end{enumerate}
We show that $X\cup\{\infty\}$ is compact. Let $I\subseteq \mathbb{N}$ be such that $(V_i)_{i\in I}$ is a covering of $X\cup\{\infty\}$. There exists some $2j+1\in I$ such that 
$\infty \in V_{2j+1}=(X\cup\{\infty\})\setminus K_j$. 
Since $K_j$ is compact and $(U_i)_{i\in I}$ covers
 $K_j$ we have that there exists a finite set $a\subseteq I$ such that $K_j\subseteq \bigcup_{i\in a} U_i$ and 
 so $(U_i)_{i\in a\cup\{2j+1\}}$ is a finite subcovering of $X\cup \{\infty\}$. So the one point compactification of $X$ is compact.\\
\\ 
Assume $X$ is effectively $T_2$, and $H_0,H_1$ the functions that witness that $X$ is effectively $T_2$ then the functions $\widehat{H}_0,\widehat{H}_1$ such that for all $x,y\in X$ we have
\begin{equation*}
\widehat{H}_0(x,y)=2\cdot H_0(x,y) \quad \quad \text{ and } \quad \quad \widehat{H}_1(x,y)=2\cdot H_1(x,y)
\end{equation*}
For all $x\in X$, let $j=\min\{i: x\notin K_i\}$. Define
\begin{equation*}
H_0(\infty,x)=2\cdot k+1
\end{equation*}
and
\begin{equation*}
\widehat{H}_1(\infty,x)=2\cdot f_j(r) \;\quad\text{ where }\;\quad r=\min\{t: x\in f_j(t)\}
\end{equation*}
The functions $\widehat{H}_0,\widehat{H}_1$ exist by $\Delta^0_1$ comprehension and witness that $X\cup\{\infty\}$ is effectively $T_2$.\\
\\
Assume that $X$ is effectively $T_2$ and there exists a sequence $(C_n)_{n\in\mathbb{N}}$ such that for all $n$ $C_n$ is the covering relation for $K_n$.  Given a finite set $F\subseteq I$, let $j=\max\{i: 2i+1\in F\}$, we have $\bigcup_{i\in F} V_i=X\cup\{\infty\}$ if and only if $\{i: 2i\in F\}\in C_j$, so the covering relation for $X\cup\{\infty\}$ is $\Delta^0_1$ definable relative to $(C_n)_{n\in\mathbb{N}}$ and so $X\cup\{\infty \}$ is effectively compact.
\end{defin}

\begin{prop}\label{alexunique} Over $\mathbf{RCA}_0$ let $(X,(U_i)_{i\in\mathbb{N}},k) $ be a $T_2$ CSC space and $(K^0_i,U^0_i,f^0_i,g^0_i)_{i\in\mathbb{N}}$ and $(K^1_i,U^1_i,f^1_i,g^1_i)_{i\in\mathbb{N}}$ be two exhaustions by compact sets of $X$. The CSC spaces
\begin{equation*}
(X\cup\{\infty\},(U_{i})_{i\in\mathbb{N}}\cup(X\cup\{\infty\}\setminus K^0_n)_{n\in\mathbb{N}},\widehat{k}^0)
\end{equation*} 
and 
\begin{equation*}
(X\cup\{\infty\},(U_{i})_{i\in\mathbb{N}}\cup(X\cup\{\infty\}\setminus K^1_n)_{n\in\mathbb{N}},\widehat{k}^1)
\end{equation*} 
are effectively homeomorphic. That is, the topology of the one point compactification does not depend on the choice of the exhaustion by compact sets.\\
\\
\begin{proof}
We show that the identity $Id:X\cup\{\infty\}\rightarrow X\cup\{\infty\}$ is an effective homeomorphism. Let $v:X\times \mathbb{N}\rightarrow \mathbb{N}$ be the function given by:
\begin{enumerate}[label={(\arabic*)}]
\item $v(x,2i)=2i$.
\item $v(x,2i+1)=2\cdot g(r)$  where  $ r=\min\{s: x\in U_{g(s)} \}$
\item $v(\infty, 2j+1)= \min\{r: \{ f_i(m): i,m<r\}\in C_i\}$
\end{enumerate}
We have that $v$ exists by $\Delta^0_1$ comprehension and witnesses that the identity $Id:X\cup\{\infty\}\rightarrow X\cup\{\infty\}$ is  effectively continuous. By symmetry, we have that the identity is also effectively open, and so the identity is an effective homeomorphism.
\end{proof}
\end{prop}

\begin{prop}\label{alexemb} $\textbf{RCA}_0$ proves every l.c.\ effectively $T_2$ CSC space with a effective $CCN$ $(X,(U_i)_{i\in \mathbb{N}},k)$  has a one point compactification $X\cup\{\infty\}$ which  effectively compact effectively $T_2$ and the inclusion  $X\rightarrow X\cup\{\infty\}$ is an effective embedding.\\
\\
\begin{proof} By Proposition~\hyperref[Chonk]{\ref{Chonk}} we have that $X$ has an exhaustion by compact sets $(K_i,U_i,f_i,g_i)_{i\in\mathbb{N}}$ and there is a sequence $(C_i)_{i\in\mathbb{N}}$ such that for all $i\in\mathbb{N}$ we have that $C_i$ is the covering relation for $K_i$. By the observations made in Definition~\hyperref[defprop]{\ref{defprop}}, we have that the one point compactification  of $X$ exists and is an effectively $T_2$ effectively compact CSC space. We show that the inclusion $X\rightarrow X\cup\{\infty\}$ is an effective homeomorphism. We define $v:X\times \mathbb{N}\rightarrow \mathbb{N}$ where $v(x,2i)=i$ and $v(x,2i+1)=g(r)$ where
 $r=\min \{s\in \mathbb{N}: x\in U_{g_i(s)}\}$. The function $v$ exists by $\Delta^0_1$ comprehension and witnesses that the inclusion $X\rightarrow X\cup\{\infty\}$ is effectively continuous. We have that the map $(x,i)\mapsto 2i$ witnesses that the inclusion is effectively open, so the inclusion of $X$ in $X\cup\{\infty\}$ is an effective homeomorphism.
\end{proof}
\end{prop}

\begin{prop}\label{awful} Over $\textbf{RCA}_0$, for all CSC spaces $(X,(U_i)_{i\in\mathbb{N}},k)$, if there is an effective embedding $f:X\rightarrow \widehat{X}$, where $(\widehat{X},(V_i)_{i\in\mathbb{N}},k')$ is compact, $eT_2$, and $\widehat{X}\setminus \text{rng}(f)$ contains exactly one point, then $X$ has a $CCN$. If $\widehat{X}$ is also effectively compact, then $X$ has an exhaustion by compact sets and the one point compactification of $X$ is effectively homeomorphic to $\widehat{X}$.\\
\\
\begin{proof} Let $p$ be the unique element of $\widehat{X}\setminus\text{rng}(f)$. Let $H_0,H_1$ be the functions that witness that $\widehat{X}$ is effectively $T_2$ and $v$ witness that $f$ is effectively continuous. The sequence
\begin{equation*}
(f^{-1}(\widehat{X}\setminus V_{H_0(p, f(x))}), v(x,H_1(p, f(x))))_{x\in X}
\end{equation*}
will be a $CCN$ for $X$. For the second part, assume that $\widehat{X}$ is also effectively compact and let $C$ be its covering relation. We have that $\widehat{X}$ will be $uT_3$ by Theorem~\hyperref[compact T3]{\ref{compact T3}}. By Theorem~\hyperref[Big4]{\ref{Big4}}, we can assume without loss of generality that $X$ has an algebra of clopen sets. Let $Comp$ and $Int$ be functions that witness $\widehat{X}$ has an algebra of clopen sets. As in the proof of Theorem~\hyperref[compwell]{\ref{compwell}}, let $(V_{i(n)})_{n\in\mathbb{N}}$  define a decreasing cofinal sequence of neighborhoods of $p$. For ease of notation let $j(n)=Comp(i(n))$,  we have
\begin{equation*} (f^{-1}( V_{j(n)}),f^{-1}( V_{j(n)}),(v(x,j(n))_{x\in f^{-1}( V_{j(n)})}, (v(x,i(n))_{x\in f^{-1}( V_{i(n)})})_{n\in\mathbb{N}},
\end{equation*}
will be an exhaustion by compact sets. Following the proof of Proposition~\hyperref[alexunique]{\ref{alexunique}} one shows that $\widehat{X}$ is effectively homeomorphic to the one point compactification of $X$.
\end{proof}
\end{prop}

\section{Choice of compact neighborhoods and arithmetic transfinite recursion}
In this section, we will show that every l.c.\ CSC space has a $CCN$ is equivalent to $\textbf{ATR}_0$ over $\textbf{RCA}_0$. $\textbf{ATR}_0$ proves the existence of $CCN$ for l.c.\ CSC space by the following lemma
\begin{lemma}\label{prcd} \textbf{(Simpson \cite[Lemma VIII.4.7]{Simp})} For any $\Pi^1_1$ formula $\varphi(x,i,X)$, where $X$ is the only set variable in $\varphi$, $\textbf{ATR}_0$ proves
\begin{equation*}
\forall X\,(\forall x\,\exists i\,\varphi(x,i,X)\rightarrow \exists f\;\forall x\,\varphi(x,f(x),X))
\end{equation*}
\end{lemma}
\noindent
We introduce a special case of \hyperref[prcd]{\ref{prcd}} and prove that it is equivalent to arithmetic transfinite recursion over $\textbf{RCA}_0$.
\begin{defin} An eventually well-founded tree array is a collection $(T_i^j)_{i,j\in \mathbb{N}}$ of trees such that for all $j\in \mathbb{N}$ the sequence $(T^j_i)_{i\in\mathbb{N}}$ is eventually well-founded. A modulus for an eventually well-founded array of trees $(T_i^j)_{i,j\in \mathbb{N}}$ is  a sequence $(n_j)_{j\in\mathbb{N}}$ such that for all $j$ and all $i\geq n_j$ the tree $T^j_i$ is well-founded. By $TAM$, we mean the statement that every eventually well-founded array of trees admits a modulus. By $1TAM$, we mean $TAM$ restricted to arrays of trees that have at most $1$ branch each.
\end{defin}
\begin{prop}\textbf{(Simpson \cite[Exercise VIII.4.25]{Simp})}\label{TSP>BTC} $\textbf{RCA}_0+1TAM$ implies arithmetic transfinite recursion.
\end{prop}

\begin{lemma} The following are equivalent over $\textbf{RCA}_0$:
\begin{enumerate}[label={(\arabic*)}]
\item $TAM$
\item For any sequence of sequences of linear orders $(L^j_i)_{i,j\in\mathbb{N}}$ such that for each $j$ the sequence $(L^j_i)_{i\in\mathbb{N}}$ is decreasing with respect to inclusion and eventually well-ordered then there is a sequence $(n_j)_{j\in\mathbb{N}}$ such that for all $j\in\mathbb{N}$ and all $i\geq n_j$ $L^j_i$ is a well order.
\end{enumerate}
\begin{proof} We have that $TAM$ and $(2)$ both imply arithmetic comprehension over $\textbf{RCA}_0$. So it suffices to show that $TAM$ is equivalent to $(2)$ over $\textbf{ACA}_0$.\\
\\
$(1\rightarrow 2)$ follows from the fact that a linear order $L$ is a well order if and only if
\begin{equation*}
T(L)=\{\sigma\in L^{<\mathbb{N}}:\forall j<|\sigma|\;\sigma(j+1)<_L\sigma(j)\}
\end{equation*}
is well-founded.
\\
\\
Assume $(2)$. Let $(T_i^j)_{i,j\in\mathbb{N}}$ be a sequence of trees  such that for all $j\in\mathbb{N}$ the sequence $(T^j_i)_{i\in\mathbb{N}}$ is eventually well-founded. Define
\begin{equation*}
S^j_n=\{\emptyset\}\cup \{(m)^\frown \sigma:m\geq n\wedge\sigma\in T^j_m\}=\coprod_{m\geq n} T^j_m
\end{equation*}
We observe that $T^j_n$ is well-founded if and only $(S^j_n,<_{\text{KB}})$ is a well order and that for all $m\geq n$ $S^j_m\subseteq S^j_n$. The family of linear orders $((S^j_n,<_{\text{KB}}))_{n,j\in\mathbb{N}}$
 satisfies the conditions of $(2)$ so there exists by assumption a sequence $(m_j)_{j\in\mathbb{N}}$ such that $\forall i\geq m_j$ $(S^j_i,<_{\text{KB}})$ is a well order. So $\forall i\geq m_j\; T^j_i$ is well-founded. 
\end{proof}
\end{lemma}

\begin{prop}\label{ATR>TSP} $\textbf{ATR}_0$ proves $TAM$.\\
\\
\begin{proof} Let $(L^j_i)_{i,j\in\mathbb{N}}$ be as in the previous lemma. Fix a $j\in\mathbb{N}$, we wish to show that there exists $m_j$ such that $\forall k\in\mathbb{N}$ $L^j_{m_j}\cong L^j_{m_j+k}$. By assumption, there is an $n_j$ such that $L^j_i$ is well-ordered for all $i\geq n_j$. Let $L=L^j_{n_j}+1$, by using the comparability of well orders and that $\forall i\;L^j_{n_j+i}\subseteq L^j_{n_j}$ we have that every $L^j_{n_j+i}$ is isomorphic to a proper initial segment of $L$. So
\begin{equation*}
\forall i\; \exists a\in L\; \exists f:L^j_{n_j+i} \xrightarrow{\sim} \{b\in L:b<a\}
 \end{equation*}
By $\Sigma^1_1$ choice, we have that
\begin{equation*}
\exists (a_i)_{i\geq n_j}\; \exists (f_i)_{i\geq{n_j}}\;\forall i\geq n_j\;  a_i\in L \wedge f_i:L^j_{i} \xrightarrow{\sim} \{b\in L:b<a_i\}
\end{equation*}
We have that for all $i_0<i_1$ that $a_{i_0}\geq a_{i_1}$ and since $L$ is a well order the sequence $(a_i)_{i\geq n_j}$ is eventually constant. Let $m_j$ be such that for all $i\in\mathbb{N}$ $a_{m_j}=a_{m_j+i}$ then we have for all $i\geq m_j$ that $L^j_{m_j}\cong L^j_{m_j+i}$. In particular we have that for all $i\geq m_j$ the function $f^{-1}_i\circ f_{m_j}$ is an isomorphism between $L^j_{m_j}$ and $L^j_{m_j+i}$.\\
\\
We have, therefore, that
\begin{equation*}
\forall j\;\exists m\; \forall k\; \exists f:L^j_m\xrightarrow{\sim} L^j_{m+k}
\end{equation*}
Using $\Sigma^1_1$ choice, we have
\begin{equation*}
\forall j\; \exists (f_k)_{k\in\mathbb{N}}\; \exists m\; \forall k\; f_k:L^j_m\xrightarrow{\sim} L^j_{m+k}
\end{equation*}
Using $\Sigma^1_1$ choice again, we get
\begin{equation*}
 \exists (f^j_k)_{j,k\in\mathbb{N}}\;\; \exists (m_j)_{j\in \mathbb{N}} \;\;\forall j\; \forall k\;f^j_k:L^j_m\xrightarrow{\sim} L^j_{m+k}
\end{equation*}
We observe that for all $j\in\mathbb{N}$ $L^j_{m_j}$ must be a well order so $(m_j)_{j\in\mathbb{N}}$ is a sequence with the desired property.
\end{proof}
\end{prop}

\begin{prop}\label{TSP>CCN} $\textbf{ATR}_0$ proves that every l.c.\ CSC space has a choice of compact neighborhoods.
\\
\\
\begin{proof} Let $(X,(U_i)_{i\in\mathbb{N}},k)$ be a l.c.\ CSC space and let $X=(x_n)_{n\in\mathbb{N}}$ be an enumeration of the points of $X$. Let $f:X\times\mathbb{N}\rightarrow \mathbb{N}$ be the function given by
\begin{equation*}
f(x,i)=\begin{cases} \min\{ j> f(x,i-1): x\in U_j\} \text{ if such } j \text{ exists}\\
\max\{j:x\in U_j\}\quad \quad \quad \text{otherwise}
\end{cases}
\end{equation*}
Intuitively, $f(x,\cdot)$ lists out the indices of basic neighborhoods of $x$.
For each $x\in X$ let $(T^x_i)_{i\in\mathbb{N}}$ be the sequence of trees where for all $i\in \mathbb{N}$
 $T^x_i$ is the set of strictly increasing sequences $\sigma$ such that
\begin{equation*}
\forall \tau  \sqsubsetneq \sigma\, (\overline{U_{f(x,i)}}\not\subseteq\bigcup_{n<|\tau|} U_{\tau(n)}) \wedge \forall n<|\sigma |\;(x_n\in \overline{U_{f(x,i)}}\rightarrow \exists m\leq n\;x_n\in  U_{\sigma(m)} )
\end{equation*}
That is, all sequences $\sigma$ such that none of its proper initial segments defines a covering for $\overline{U_{f(x, i)}}$ and for each $n<|\sigma|$ if $x_n\in \overline{U_{f(x,i)}}$ then there is an $m\leq n$ such that $x_n\in U_{\sigma(m)}$.
An infinite branch in $T^x_i$ defines a covering of $\overline{U_{f(x,i)}}$, which does not have finite subcovering. Similarly, an infinite cover of $\overline{U_{f(x,i)}}$ that does not admit a finite subcover defines a branch in $T^x_i$. Therefore, $\overline{U_{f(x,i)}}$ is compact if and only if $T^x_i$ is well-founded. Since $X$ is l.c.\, we have that for each $x\in X$, the sequence $(T^x_i)_{i\in\mathbb{N}}$ is eventually well-founded. Thus the modulus $(n_x)_{x\in X}$ given by $TAM$ is such that for all $x\in X$ the neighborhood $\overline{U_{f(x,n_x)}}$ is compact. So $(n_x,\overline{U_{f(x,n_x)}})_{x\in X}$ exists by arithmetic comprehension and is a $CCN$ for $X$.
\end{proof}
\end{prop}

\begin{prop}\label{CCN>ACA} Over $\textbf{RCA}_0$ every $T_2$ l.c.\ CSC space having a basis of compact neighborhoods implies arithmetic comprehension.\\
\\
\begin{proof} Let $A\subseteq \mathbb{N}$ be a set, we show that its Turing jump exists. Set
\begin{equation*}
X=\mathbb{N}\times(\mathbb{N} \cup\{\infty\})
\end{equation*}
and define the basic open sets of $X$ to be
\begin{equation*}
U_{2(e,t)}=\begin{cases}\{(e,\infty)\}\quad\quad\;\text{ if }\quad \Phi^A_e(e){\downarrow}_{\leq t}\\
\{(e,n):n\geq t\}\cup \{(e,\infty)\}\quad \text{ otherwise }
\end{cases}
\end{equation*}
and
\begin{equation*}
U_{2(e,s)+1}=\{(e,s)\} 
\end{equation*}
Define $k((e,n),a,b)=2(e,n)+1$ and $k((e,\infty),a,b)=\max\{a,b\}$.
$X$ can be viewed as the disjoint union of the spaces $X_e=\{e\}\times(\mathbb{N}\cup\{\infty\})$. It is straightforward to show that $X$ is l.c.\ We observe that all of the basic sets are clopen and that $X$ is $T_3$, and it is easy to show that it is, in fact, $uT_3$.\\
\\
Let $(i(x))_{x\in X}$ be a $CCN$ for $X$. Let $m_e$ be the unique number such that $2(e,m_e)=i((e,\infty))$. For each $e\in \mathbb{N}$, $\Phi^A_e(e){\downarrow}$ if and only if $U_{i(e,\infty)}$ is isolated. In particular by construction, we have that $\Phi^A_e(e){\downarrow}$ if and only if $\Phi^A_e(e){\downarrow}_{m_e}$. So $A'\leq_T (m_e)_{e\in\mathbb{N}}$ exists by $\Delta^0_1$ comprehension.
\end{proof}
\end{prop}

\begin{prop}\label{T2CCN>ATR} Over $\textbf{RCA}_0$ the statement that every $T_2$ l.c.\ CSC space has a $CCN$ implies arithmetic transfinite recursion.\\
\\
\begin{proof} By Proposition~\hyperref[CCN>ACA]{\ref{CCN>ACA}} we may work over $\textbf{ACA}_0$. We show that every $T_2 $ l.c.\ CSC space has a $CCN$ implies $1TAM$ which is equivalent to $\textbf{ATR}_0$ over $\textbf{RCA}_0$. Let $(T_j^i)_{i,j\in\mathbb{N}}$ be an eventually well-founded array of trees as in the condition of $1TAM$ and let $X$ be equal to the disjoint union of the spaces $X^j=(\coprod_{i\in\mathbb{N}}T^j_i,<_{\text{KB}})\cong \sum_{i\in\mathbb{N}} (T^j_i,<_{\text{KB}})+1$ with the upper limit topology where
\begin{equation*}
\coprod_{i\in\mathbb{N}}T^j_i=\{\emptyset\}\cup \bigcup\{(m)^\frown T^j_m:m\in\mathbb{N}\}
\end{equation*}
We show that for every $j\in\mathbb{N}$, the space $X^j$ is l.c.\ To do this, it suffices to show that a tree with at most one branch is l.c.\ with the upper limit topology induced by the Kleene-Brouwer order. If $T$ is well-founded, then by Proposition~\hyperref[wocomp]{\ref{wocomp}} $T$ with the topology induced by the Kleene-Brouwer order is compact and therefore l.c.\ Let $T$ be a tree with exactly one branch $f$ and $x\in T$ be a point. If $\forall n\in\mathbb{N}\;x<_{\text{KB}} f_{\leq n}$ then $\mathopen]-\infty,x\mathclose]$ is well-ordered with a maximal element and therefore is a compact neighborhood of $x$. Otherwise, if there is an $n$ such that $f_{\leq n} <_{\text{KB}} x$ then the interval $\mathopen]f_{\leq n},x\mathclose]$ is well-ordered with respect to the Kleene-Brouwer order and so it is a compact neighborhood of $x$.\\
\\
By assumption, there is a choice of compact neighborhoods for the space $X$. This implies that we have a choice of a compact neighborhood for $\emptyset\in X_j$, which we denote by $K_j$. We define
\begin{equation*}
m_j=\min\{m\in\mathbb{N}: \mathopen](m),\emptyset\mathclose]\subseteq K_j\}+1
\end{equation*}
Recall that the upper limit topology on a linear order with maximal element is compact if and only if it is well-ordered. So for all $i\geq m_j$ $(T^i_j,<_{\text{KB}})$ is a well order, and so $T^i_j$ is well-founded.
\end{proof}
\end{prop}

\section{Embedding scattered CSC spaces into well orders}
We have that $\textbf{ATR}_0$ proves that every l.c.\ $T_2$ CSC space has a $CCN$ and therefore embeds into a compact CSC space. We would like to show that we can embed $T_3$ scattered CSC spaces into a well order. For scattered linear orders, we have the following results.
\begin{thm}\label{Cl}\textbf{(Clote \cite{Clote})} Arithmetic transfinite recursion is equivalent to every countable scattered linear order having a countable set of initial segments. That is, for any scattered linear order $(L,<_L)$ there is a sequence of sets $(I_n)_{n\in\mathbb{N}}$ such that
\begin{equation*}
\forall X\subseteq L\; [(\forall x\in X\; \forall y<_L x( y\in X))\rightarrow \exists n\; X=I_n]
\end{equation*}
In particular, we consider the empty set and all of $L$ to be initial segments.
\end{thm}

\begin{obs} The initial segments of a linear order are naturally ordered by inclusion. Over $\textbf{ACA}_0$ if a linear order $(L,<_L)$ has countably many initial segments, then we can order by inclusion the initial segments $(I_n)_{n\in\mathbb{N}}$ which are closed in the order topology of $L$. We call this order the Dedekind completion of $L$ and the map $l\mapsto\{j\in L:j\leq_L l\}$ is an embedding from $L$ into its Dedekind completion.
\end{obs}

\begin{thm}\label{Shafer}\textbf{(Shafer \cite{Shafer})} $\textbf{WKL}_0$, and in particular $\textbf{ACA}_0$, proves the order topology of every complete linear order is compact.
\end{thm}

\begin{cor}\label{scattered in comp} $\textbf{ATR}_0$ proves that the order topology of any scattered linear order embeds into a $T_2$ compact CSC space.\\
\\
\begin{proof} By Theorem~\hyperref[Cl]{\ref{Cl}} we have that the Dedekind completion of $X$ is countable, and so $X$ embeds into a complete linear order which by Theorem~\hyperref[Shafer]{\ref{Shafer}}, it will be compact.
\end{proof}
\end{cor}
\noindent
The proof of Theorem~\hyperref[Cl]{\ref{Cl}} requires some results involving the Hausdorff rank, which is not a topological invariant. We also have that non scattered linear orders may have scattered order topology. For example, $\mathbb{Q}\times\mathbb{Z}$ with the lexicographic order is non scattered as a linear order, but its order topology is discrete. So Corollary~\hyperref[scattered in comp]{\ref{scattered in comp}} does not ensure that every $T_3$ scattered CSC space will embed into a $T_2$ compact CSC space. We will need to use the Cantor-Bendixson rank instead. This will allow us to show that every $T_3$ scattered CSC space is homeomorphic to a scattered linear order with its order topology. The Cantor-Bendixson derivative and rank of CSC spaces in reverse mathematics has already been studied by Mont{\'a}lban and Greenberg \cite{Greenberg}, Friedman \cite{FriedC} in the case of countable metric spaces, and by Friedman and Hirst \cite{First} in the form of characteristic systems.

\begin{defin} Let $(X,(U_i)_{i\in\mathbb{N}},k)$ be a CSC space then we write
\begin{equation*}
D(X)=\{x\in X: x\text{ is not isolated}\}
\end{equation*}
We call $D(X)$ the sets of limit points of $X$. We observe that $D(X)$ is arithmetical relative to $X$, however, in general it will not be recursive.
\end{defin}

\begin{defin}\label{rankdef} Over $\textbf{ACA}_0$ we say that a Cantor-Bendixson rank or rank for a CSC space $X$ is a well order $R$ such that there is a sequence $(X_r)_{r\in R}$ such that
\begin{enumerate}[label={(\arabic*)}]
\item $X_{0}=X$.
\item For all $r\in R$ that is not maximal $X_{r+1}=D(X_r)$.
\item For every $r\in R$ which is limit $X_r=\bigcap_{s<_R r} X_s$.
\item $D(\bigcap_{r\in R} X_r)=\bigcap_{r\in R} X_r$.
\end{enumerate}
We call a well order such that $(1),(2),(3)$ hold, but not necessarily $(4)$, a partial rank. A rank $R$ of $X$ is minimal if none of its initial segments are a rank for $X$. Given a rank $R$ on $X$ we write 
\begin{equation*}\text{rank}_R(x)=\text{rank}(x)=\min\{r\in R:x\notin X_{r+1}\}
\end{equation*}
It is straightforward that for any $x$ of rank $r$, we have $x\in X_r$.
\end{defin}

\begin{obs} If $X$ is a $T_3$ space and $R$ is a rank for $X$ then by definition of rank we have $D(\bigcap_{r\in R} X_r)=\bigcap_{r\in R} X_r$. So $\bigcap_{r\in R} X_r$ is invariant under the $D$ and, therefore, does not have any isolated points. This means that $\bigcap_{r\in R} X_r$ is either homeomorphic to $\mathbb{Q}$ or it is empty. In particular, we have that if $X$ is scattered, then $\bigcap_{r\in R} X_r=\emptyset$.
\end{obs}

\begin{obs} Since over $\textbf{ACA}_0$ we have arithmetic transfinite induction, we have that any rank $R$ for a $T_3$ scattered CSC space $X$ has an initial segment that is a minimal rank. Over $\textbf{ACA}_0$, any two minimal ranks $R_0$ and $R_1$ will be order isomorphic. We may say, over $\textbf{ACA}_0$, that $R$ is the Cantor-Bendixson rank of $X$ if $R$ is a minimal rank.
\end{obs}

\begin{thm}\label{ATR>rank}\textbf{(Friedman \cite[Lemma 15]{FriedC})} $\textbf{ATR}_0$ proves that every $T_3$ scattered space is ranked.
\end{thm}

\begin{prop}\label{rankwo>aca}(See \cite[Proposition 6.36]{Greenberg}) Over $\textbf{RCA}_0$ every well order has a rank implies arithmetic comprehension. \\
\\
\begin{proof}
Let $A\subseteq \mathbb{N}$ be a set. Let $L_e=\{t\in \mathbb{N}:\neg\Phi^A_e(e){\downarrow}_{\leq t}\} +1$ and $L=\sum_{e\in\mathbb{N}} L_e$ which is the sum of well orders and so it is a well order. We can compute $A'$ from the set of non isolated points of $L$, which is computable from the rank of $L$.
\end{proof}
\end{prop}

\noindent
We now will show that over $\textbf{ACA}_0$ a scattered $T_3$ CSC space with a rank embeds into a well order. To show these results, we make use of a construction similar to that done in {Theorem~\hyperref[compwell]{\ref{compwell}}}. 

\begin{C} Working over $\textbf{ACA}_0$, let $(X,(U_i)_{i\in\mathbb{N}})$ be a $T_3$ scattered CSC space, which, without loss of generality, we assume it has an algebra of clopen sets. Let $R$ be a rank for $X$. We now lay out an arithmetic procedure to associate to each point $x$ a unique sequence $\alpha_x$, which we call the address of $x$. We will see later that the Kleene-Brouwer order on the set of addresses is a well order and that the map $x\mapsto \alpha_x$ will be an embedding.\\
\\
Let $F:X\times \mathbb{N}\rightarrow \mathbb{N}$ be a partial function such that for all $x$ in $X$ if $x$ is isolated $U_{F(x,0)}=\{x\}$ otherwise we define $F(x,n)$ to be the least $s\in \mathbb{N}$ such that
\begin{enumerate}[label={(\arabic*)}]
\item $U_s\cap X_{\text{rank}(x)}=\{x\}$, or rather, all of the elements in $U_s$ besides $x$ are of lower rank. Since $x$ is isolated in $X_{\text{rank}(x)}$, there must be a neighborhood of $x$ containing $x$ and points of rank strictly lower than $x$.
\item $\forall m< n\;U_s\subsetneq U_{F(x,m)}$.
\end{enumerate}
We have that $F$ exists by arithmetic comprehension. Informally, $F$ lists out a descending sequence of neighborhoods for every point. For ease of notation, we will write $U(x,n)=U_{F(x,n)}$. We observe that the sets of the form $U(x,n)$ form a basis of clopen sets for $X$.\\
\\
For each $x\in X$ and $n\in \mathbb{N}$ let $A^{(x,n)}= U(x,n)\setminus U(x,n+1)$ and $A^{-1}=X$. For $h=-1$ or $h=(x,n)$ we define inductively two sequences $\sigma^h\in X^{<\mathbb{N}}$ and $\tau^{h}\in \mathbb{N}^{<\mathbb{N}}$. If $A^h\setminus \bigcup_{i<k}U(\sigma^h(i),\tau^h(i))$ is empty we terminate the construction, otherwise let $\sigma^h(k)$ be the first element in $A^h\setminus \bigcup_{i<k}U(\sigma^h(i),\tau^h(i))$. Define
\begin{equation*}
\tau^h(k)=\min\left\{s\in\mathbb{N}:U(\sigma^h(k),s)\subseteq A^h\setminus\bigcup_{i<k}U(\sigma^h(i),\tau^h(i))\right\}
\end{equation*}
If $A^h$ is compact then there exists a $k$ such that $A^h=\bigcup_{i<k}U(\sigma^h(i),\tau^h(i))$ and at such $k$ we halt the construction. Otherwise, the construction might go on forever, and $\sigma^h$ and $\tau^h$ will be functions. We observe that the sets $(U(\sigma^h(i),\tau^h(i)))_{i\in \text{dom}(\sigma)}$ form a partition of $A^h$ into clopen sets. Since $\sigma^h$ and $\tau^h$ are uniformly arithmetically defined with respect to $h$ we have that
\begin{equation*}
\{(\sigma^h,\tau^h):h=-1\vee (h=(x,n)\wedge x\in X\wedge n\in \mathbb{N})\}
\end{equation*}
exists by arithmetic comprehension.
\\
\\
For each point $x\in X$ we define $\alpha_x,\beta_x\in \mathbb{N}^{<\mathbb{N}}$,  and $\gamma_x\in\mathbb{N}$ as:
\begin{enumerate}[label={(\arabic*)}]
\item  $\alpha_x(0)$ is the unique $i$ such that $x\in U(\sigma^{(-1)}(i),\tau^{(-1)}(i))$ and $\beta_x(0)=\sigma^{(-1)}(i)$.
\item $\alpha_x(2n+1)$ is the unique $m$ such that $x\in U(\beta_x(n),m)\setminus U(\beta_x(n),m+1)$.
\item $\alpha_x(2n+2)$ is the unique $i$ such that
\begin{equation*}
x\in U(\sigma^{(\beta_x(n),\alpha_x(2n+1))}(i),\tau^{(\beta_x(n),\alpha_x(2n+1))}(i))
\end{equation*}
and let
\begin{equation*}\beta_x(n+1)=\sigma^{(\beta_x(n),\alpha_x(2n+1))}(i)
\end{equation*}
If $x=\beta_x(n+1)$, then $\gamma_x=0$ if $x$ is isolated and $\gamma_x=\tau^{(\beta_x(n),\alpha_x(2n+1))}(i)$ otherwise and the construction terminates.
\end{enumerate}
Since for all $n$ $\text{rank}(\beta_x(n))>\text{rank}(\beta_x(n+1))$, we have that for all $x$ the construction above must eventually terminate.
The set of all $(\alpha_x,\beta_x,\gamma_x)_{x\in X}$ is arithmetically definable, and so it exists by arithmetic comprehension.\\
\\
For each $x\in X$ we call $\alpha_x$ the address of $x$. Since $X$ is Hausdorff, we have that the map $x\mapsto \alpha_x$ is injective. The collection of addresses $\mathcal{A}$ will not be a tree since every address has odd length. But for any $\alpha_x\in\mathcal{A}$, every odd length initial segment of $\alpha_x$ will be an address.  In fact, for each $n<|\beta_x|$  $\alpha_{\beta_x(n)}$ is the initial segment of $\alpha_x$ of length $2n+1$. We observe that for all $x,y\in X$ that $y\in U(x,\gamma_x)$ if and only if $\alpha_x\sqsubseteq \alpha_y$. By construction, we have that if $\alpha_x\sqsubseteq \alpha_y$ then $\beta_x\sqsubseteq \beta_y$. In particular, since the rank of $\beta_x(n)$ is strictly decreasing, there cannot exist an infinite increasing sequence of addresses.  \\
\\
Let $<_\alpha$ be the order given by
\begin{equation*}
x<_\alpha y\leftrightarrow \alpha_x<_{\text{KB}} \alpha_y
\end{equation*} 
Since the Kleene-Brouwer ordering on a well-founded tree is a well order, we have that $<_\alpha$ defines a well order on $X$. In general, we have that $<_\alpha$ induces a coarser topology on $X$. We show that the map $x\mapsto \alpha_x$ defines an embedding from $X$ to ${\downarrow} \mathcal{A}$. That is, the map $x\mapsto \alpha_x$ is a homeomorphism between $(X,(U_i)_{i\in\mathbb{N}})$ and $\mathcal{A}$ with the subspace topology.\\
\\
If $x$ is isolated in $X$, then we have that $\alpha_x$ does not have any extensions in $\mathcal{A}{\downarrow}$, so it has a predecessor with respect to the  Kleene-Brouwer order on ${\downarrow}\mathcal{A}$. So $x$ is isolated in $\mathcal{A}$ with the subspace topology. If $x$ is not isolated in $X$, then we have that $U(x,\gamma_x)$ will be infinite. In particular, for all $m\geq \gamma_x$, the set $U(x,m)\setminus U(x,m+1)$ is non empty. Given a $\sigma \in {\downarrow}\mathcal{A}$ such that $\sigma<_{\text{KB}} \alpha_x$ then there exists an $m\geq \gamma_x$ such that $\sigma<_{\text{KB}}\alpha_x^\frown (m)$. For  any $y\in U(x,m+1)\setminus U(x,m+2)$ we have $\alpha_x^\frown (m)<_{\text{KB}} \alpha_y<_\text{KB}\alpha_x$. This implies that $\alpha_x$ is not isolated in $\mathcal{A}$ with the subspace topology. So $x$ is isolated in $X$ if and only if $\alpha_x$ is isolated in the subspace topology of $\mathcal{A}$. It suffices to show that $x\mapsto \alpha_x$ is continuous and open on the non isolated points of $X$.\\
\\
Let $x,y\in X$ such that $\alpha_y<_{\text{KB}} \alpha_x$ and $\alpha_x$ is not isolated in $\mathcal{A}$ with the subspace topology. If $y\notin U(x,\gamma_x)$ then $U(x,\gamma_x)\subseteq \left]y,x\right]_{<_\alpha}$. Otherwise there exists a unique $n\geq \gamma_x$ such that $y\in U(x,n)\setminus U(x,n+1)$, so $\alpha_x^\frown(n)\sqsubseteq \alpha_y$. Since
\begin{equation*}\forall z\in U(x,n+1)\;(\alpha_x^\frown(n+1)\sqsubseteq \alpha_z)
\end{equation*} we have that $\forall z\in U(x,n+1)$ $\alpha_y<_{\text{KB}} \alpha_z<_{\text{KB}} \alpha_x$ and so $U(x,n+1)\subseteq \left]y,x\right]_{<_\alpha}$. This proves that the map $x\mapsto \alpha_x$ is continuous.
\\
\\
Let $x\in X$ be a non isolated point and $n\in \mathbb{N}$. Since $x$ is not isolated there exists a $y\in U(x,\max\{n,\gamma_x\})\setminus\{x\}$. We have that $\alpha_x\sqsubseteq \alpha_y$ and so $\left]y,x\right]_{<_\alpha}\subseteq U(x,n)$ by construction. This proves the map $x\mapsto \alpha_x$ is open with its image. So $X$ is homeomorphic to $\mathcal{A}$ with the subspace topology.
\end{C}

\begin{thm}\label{scatter in comp} $\textbf{ACA}_0$ proves that every $T_3$ scattered CSC space with rank is homeomorphic to a subspace of a well order.\\
\\
\begin{proof} Using the construction before we have that $X$ embeds into ${\downarrow}\mathcal{A}$ with the Kleene-Brouwer order and so $X$ embeds into a well order.
\end{proof}
\end{thm}

\begin{thm} $\textbf{ACA}_0$ proves that every $T_2$ l.c.\ CSC space with a rank has a $CCN$.\\
\\
\begin{proof} Fix a point $x\in X$, we have shown that for $n\geq \gamma_x$ that $U(x,n)$ is compact if and only if for all $y\in U(x,n)$ such that $\alpha_x\sqsubseteq \alpha_y$ and all $s\in\mathbb{N}$ then $\alpha_y^\frown(s)$ finitely branches. If $U(x,n)$ is compact and $y\in U(x,n)$ then  $U(y,s)\setminus U(y,s+1)$ will also be compact and so $\alpha_y^\frown(s)$ has finitely many extensions of length $|\alpha_y| +2$. If instead for all $y$ such that $\alpha_x\sqsubseteq \alpha_y$ and all $s$ $\alpha_y^\frown(s)$ is finitely branching, then by the previous result we have that the topology on $U(x,n)$ is the same topology as the topology induced by $<_\alpha$, and so, in particular, $U(x,n)$ will be homeomorphic to a well order with maximal element and so it will be compact. We can, therefore, verify if $U(x,n)$ is compact arithmetically, so by arithmetic comprehension, there exists a $CCN$ for $X$.
\end{proof}
\end{thm}

\begin{thm}\label{fruit} The following are equivalent over $\textbf{RCA}_0$:
\begin{enumerate}[label={(\arabic*)}]
\item Arithmetic transfinite recursion.
\item Every l.c.\ $T_2$ CSC space is well-orderable.
\item Bel'nov-Knaster-Urbanik theorem: every $T_3$ scattered CSC space embeds into a well order.
\item Every $T_3$ scattered CSC space has a rank.
\item Kat{\v e}tov's theorem: for every $T_3$ scattered CSC space $X$ there is an effectively continuous bijection from $X$ to the order topology of a well order. That is to say, every $T_3$ scattered CSC space is the refinement of the order topology of some well order.
\end{enumerate}
\begin{proof} 
$(4\rightarrow 3)$ By {Proposition~\hyperref[rankwo>aca]{\ref{rankwo>aca}}}, we have that every well order is ranked implies arithmetic comprehension; so we may work over $\textbf{ACA}_0$. By {Proposition~\hyperref[scatter in comp]{\ref{scatter in comp}}} every $T_3$ scattered CSC space with rank embeds into a well order.\\
\\
$(3\rightarrow 1)$ By {Theorem~\hyperref[theworst]{\ref{theworst}}} we have that every $T_3$ scattered CSC space embeds into a linear order implies arithmetic comprehension. So we may work over $\textbf{ACA}_0$. By {Lemma~\hyperref[lc are scattered]{\ref{lc are scattered}}} and {Proposition~\hyperref[lc>T3]{\ref{lc>T3}}} $T_2$ l.c.\ CSC spaces are $T_3$ and scattered. Since every well order has a $CCN$, we have that every $T_2$ l.c.\ CSC space embeds into a CSC space with a $CCN$. So by {Proposition~\hyperref[CCNstab]{\ref{CCNstab}}} and \hyperref[in comp>CCN]{\ref{in comp>CCN}} every $T_2$ l.c.\ CSC space has a $CCN$ which by {Proposition~\hyperref[T2CCN>ATR]{\ref{T2CCN>ATR}}} implies arithmetic transfinite recursion. \\
\\
$(4\rightarrow 2)$ Every l.c.\ $T_2$ space with Cantor-Bendixson rank has a $CCN$ and by {Proposition~\hyperref[2 for 1]{\ref{2 for 1}}} every $T_2$ l.c.\ CSC space with a $CCN$ is well-orderable.\\
\\
$(2\rightarrow 1)$ If every $T_2$ l.c.\ CSC space is well-orderable, then the upper limit topology of every well order is effectively homeomorphic to the order topology of a well order which by \hyperref[wellmess]{\ref{wellmess}} implies arithmetic comprehension. Furthermore, if every $T_2$ l.c.\ CSC space is well-orderable then every $T_2$ l.c.\ CSC space has a $CCN$ which by {Theorem~\hyperref[T2CCN>ATR]{\ref{T2CCN>ATR}}} implies arithmetic transfinite recursion.\\
\\
$(1\rightarrow 4)$ is simply {Theorem~\hyperref[ATR>rank]{\ref{ATR>rank}}}.\\
\\
$(3\rightarrow 5)$ follows from the fact that the subspaces of well orders have subspace topology, which is finer or equal to their order topology.\\
\\
$(5\rightarrow 3)$ every scattered $T_3$ CSC space is a refinement of a well order implies arithmetic comprehension by {Theorem~\hyperref[horrible]{\ref{horrible}}}. We may therefore work over $\textbf{ACA}_0$. Assume that $f:X\rightarrow W$ is a continuous bijection with $W$, we show that $X$ can embed into a well order. Let $S\subseteq W$ be the set of limit points $x\in W$ such that $f^{-1}(x)$ is isolated in $X$. Define
\begin{equation*}
\widehat{W}=S\times\{0\}\sqcup W \times  \{1\}
\end{equation*}
With the lexicographic order. We have that $\widehat{W}$ is a well order. Informally, $\widehat{W}$ is the order obtained by adding to $W$ a new point under every element of $S$. Let $g:X\rightarrow \widehat{W}$ be the map $x\mapsto (f(x),1)$. It is straightforward to show that $g$ is an embedding.
\end{proof}
\end{thm}

\begin{thm}\label{New} $\textbf{ACA}_0$ proves every subspace of a well order is homeomorphic to a scattered linear order with countably many cuts.\\
\\
\begin{proof} Let $(W,<_W)$ be a well order and let $S\subseteq W$ be a subspace. Let $L\subseteq W\setminus S$ be the collection of all $w\in W\setminus S$ such that $\mathopen]-\infty,w\mathclose[_{<_W} \cap S$ does not have a $<_W$-maximum element. Since $L\subseteq W$ is well-ordered, every element is maximal or has a successor. For each $l\in \{-\infty\}\cup L$, let $l^+\in L\cup\{+\infty\}$ denote the successor of $l$. We observe that for every $l\in L$, the point $\min\{s\in S: l<_W s\}$ will be isolated with respect to the subspace topology but will be a limit point with respect to the order topology (see {Observation~\hyperref[Path]{\ref{Path}}}).\\
\\
The sequence $(\mathopen]l,l^+\mathclose]_{<_W})_{l\in L\cup\{-\infty\}}$ defines a partition of $W$ into clopen sets. So we have
\begin{equation*}
W\cong\coprod_{l\in L\cup\{-\infty\}} \mathopen]l,l^+\mathclose]_{<_W}\quad \text{ and }\quad S\cong\coprod_{l\in L\cup\{-\infty\}} \mathopen]l,l^+\mathclose]_{<_W}\cap S
\end{equation*}
We define a new order $<_N$ on $W$ by flipping and rearranging the intervals $\mathopen]l,l^+\mathclose]$ so that the subspace topology on $S$ will be the same as the order topology of $<_N$. Let $(l_n)_{n\in\mathbb{N}}$ be an enumeration of the elements of $L\cup\{-\infty\}$ (we note that in general $l_n^+\neq l_{n+1}$). We define on $W$ the order $<_N$ where $a<_N b$ if and only if one of the following occurs:
\begin{enumerate}[label={(\arabic*)}]
\item $a\in \mathopen]l_n,l_n^+\mathclose]_{<_W}\wedge b\in \mathopen]l_m,l_m^+\mathclose]_{<_W}\wedge n<_{\mathbb{N}}m$
\item $\exists n\;( a,b\in \mathopen]l_{2n},l_{2n}^+\mathclose]_{<_W}\wedge a<_W b)$
\item $\exists n\;( a,b\in \mathopen]l_{2n+1},l_{2n+1}^+\mathclose]_{<_W}\wedge b<_W a)$
\end{enumerate}
It is routine to show that $<_N$ has the same order topology as $<_W$. We show that the subspace topology on $S$ is the same as the order topology of $<_N$. In general, the subspace topology is finer or equal to the order topology, so it suffices to show that the order topology is finer or equal to the subspace topology. \\
\\
Let $x\in S$ be a point and $a,b\in W$ such that $a<_N x<_N b$, we show that there are $c,d\in S$ such that $x\in \mathopen] c,d\mathclose[_{<_N}\cap S\subseteq \mathopen] a,b\mathclose[_{<_N}\cap S$. There exists a unique $n$ such that $x\in \mathopen]l_{n},l_{n}^+\mathclose]_{<_W}$. We consider the case in which $n$ is even, the odd case is proved similarly. Since $n$ is even we have by definition that $<_N$ restricted to $\mathopen]l_{n},l_{n}^+\mathclose]_{<_W}$ is equal to $<_W$. By definition of $l_n^+$ we have that $\mathopen]l_n,l_n^+\mathclose]_{<_W}\cap S$ is unbounded. Let $d=\min\{ s\in \mathopen]l_n,l_n^+\mathclose]_{<_W}\cap S: x<_Ns\}$, that is, $d$ is the successor of $x$ in $\mathopen]l_n,l_n^+\mathclose]_{<_W}\cap S$. If $x$ is the $<_W$ minimum element of $\mathopen]l_{n},l_{n}^+\mathclose]_{<_W}\cap S$ then let $c=-\infty$ if $n=0$ otherwise let $c=\min_{<_W} \mathopen]l_{n-1},l_{n-1}^+\mathclose]_{<_W}\cap S$. If $x$ is not the $<_W$ minimum in $\mathopen]l_{n},l_{n}^+\mathclose]\cap S$ then if there is an $s\in \mathopen[a,x[_{<_W}\cap S$ we set $c=s$. Otherwise, the set
\begin{equation*}
\{l\in \mathopen]l_{n},l_{n}^+\mathclose]_{<_W}: l\leq a\wedge \forall j\;(l\leq_W j <_W x\rightarrow j\notin S)\}
\end{equation*}
is non empty since it contains $a$ and has a least element $u$ since $W$ is a well order. By definition, the only element of $L$ in $\mathopen]l_{n},l_{n}^+\mathclose]_{<_W}$ is $l^+_n$, so we have that $u\notin L$. By definition of $L$, $\mathopen]l_{n}, u\mathclose[_{<_W}\cap S$ has a maximum and we set $c=\max (\mathopen]l_{n}, u\mathclose[_{<_W}\cap S)$. It is straightforward to verify that $x\in \mathopen] c,d\mathclose[_{<_N}\cap S\subseteq \mathopen] a,b\mathclose[_{<_N}\cap S$.
\\
\\
Any cut of $(S,<_N)$ will either be of the form $\mathopen]-\infty,l\mathclose]_{<_N}\cap S$, where $l\in W$, or of the form $\bigcup_{i\leq 2n+1}\mathopen]l_{i},l_{i}^+\mathclose]_{<_W}\cap S$. Since $(S,<_N)$ has countably many cuts, it is a scattered linear order.
\end{proof}
\end{thm}

\begin{cor} $\textbf{ACA}_0$ proves every ranked $T_3$ scattered CSC space is effectively homeomorphic to the order topology of a scattered linear order with countably many cuts.
\end{cor}

\begin{cor}\label{scatscat} Over $\textbf{RCA}_0$ the following are equivalent:
\begin{enumerate}[label={(\arabic*)}]
\item Arithmetic transfinite recursion.
\item Every $T_3$ scattered CSC space is effectively homeomorphic to the order topology of a scattered linear order with countably many cuts.
\end{enumerate}
\begin{proof} $\textbf{ATR}_0$ proves that every scattered $T_3$ CSC space embeds into a well order and so by {Theorem~\hyperref[New]{\ref{New}}} $\textbf{ATR}_0$ proves every scattered $T_3$ CSC space is effectively homeomorphic to a scattered linear order with countably many cuts.\\
\\
For the converse, assume every $T_3$ scattered CSC space is effectively homeomorphic to the order topology of a scattered linear order with countably many cuts. In particular, we have that every $T_3$ scattered CSC space embeds into a linear order, which by {Theorem~\hyperref[theworst]{\ref{theworst}}} implies arithmetic comprehension. So we may work over $\textbf{ACA}_0$. We also have that every scattered $T_3$ space embeds into the order topology of a complete linear order. By {Theorem~\hyperref[Shafer]{\ref{Shafer}}}, we have that the order topology of a complete linear order is compact.  So every $T_3$ scattered CSC space, and therefore every $T_2$ l.c.\ CSC space embeds, into a compact CSC space. This implies that every $T_2$ l.c.\ CSC space has a $CCN$ which by {Theorem~\hyperref[]{\ref{T2CCN>ATR}}} implies arithmetic transfinite recursion.
\end{proof}
\end{cor}

\section{Complete metrizability}
We have shown that $\textbf{ATR}_0$ proves that a $T_3$ scattered CSC space can embed into a well order with maximum element, and in particular, they can embed into a countable compact metric space. From general topology, we know that compact metric spaces are complete and the $G_\delta$ subspaces of a complete metric space are completely metrizable. For countable $T_1$ spaces, all of the subspaces are $G_\delta$, so all $T_3$ scattered CSC spaces are completely metrizable. We will show that this proof can be carried out in $\textbf{ATR}_0$. As in the previous part of this work, we will only be working with countable metric spaces and will not be considered as being codes for Polish spaces as it is usually done in reverse math. 

\begin{prop} $\textbf{RCA}_0$ proves every complete countable metric space is $T_3$ and scattered.\\
\\
\begin{proof} Let $(X,d)$ be a complete countable metric space, by {Proposition~\hyperref[metstruc]{\ref{metstruc}}} we have that for some $a\in \mathbb{R}_{>0}$ the balls $(B(x,q\cdot a))_{x\in X,q\in\mathbb{Q}}$ are clopen. Without loss of generality, we may assume $a=1$. By {Theorem~\hyperref[Big4]{\ref{Big4}}} $\textbf{RCA}_0$ proves that $(X,(B(x,q)_{x\in X,q\in \mathbb{Q}},k)$, is uniformly $T_3$. If $X$ were not scattered, then it would contain a non empty dense in itself subset $S$ that will be a uniformly $T_3$ space that does not have isolated points. We note that any neighborhood of a point in $S$ will contain infinitely many points.  Let $(q_n)_{n\in\mathbb{N}}$ enumerate the elements of $\mathbb{Q}_{>0}$. Let $x_0$ be the $<_\mathbb{N}$ first element of $S$ and set $r_0=1$. Assume that for all $i\leq n$ $y_i\in X$ and $r_i\in \mathbb{Q}_{>0}$ are defined and for all $i<n$ we have $x_i\notin B(x_{i+1},r_{i+1})\subseteq B(x_{i},r_i)$. Let $x_{n+1}$ be the $<_\mathbb{N}$ least element of $S\cap  B(x_n,r_n)\setminus\{x_n\}$, which is not empty since $S$ is dense in itself, and let $r_{n+1}$ be the first element in the enumeration of $\mathbb{Q}$ such that 
\begin{equation*}r_{n+1}\leq \min\{ d(x_{n+1},x_n),r_n-d(x_{n+1},x_n), \frac{1}{n+1}\}
\end{equation*}
We have by the triangle inequality that $x_n\notin B(x_{n+1},r_{n+1})\subseteq B(x_n,r_n)$. By construction, the sequence $(x_n)_{n\in\mathbb{N}}$ is Cauchy and, therefore, must converge to some $x\in X$. Since for all $n$ we have that $x_{n+1}$ is the $<_\mathbb{N}$ least element of $B(x_n,r_n)\cap S\setminus\{x_n\}$ and $x_{n+2}\in B(x_n,r_n)\cap S\setminus\{x_n\}$ we have that $x_{n+1}<_\mathbb{N} x_{n+2}$. In particular the sequence $(x_n)_{n\in\mathbb{N}}$ is strictly $<_\mathbb{N}$ increasing. For every $n\in \mathbb{N}$ we have that $(x_n)_{n\in\mathbb{N}}$ is eventually in $B(x_n,r_n)$ and so $x\in \overline{B(x_n,r_n)}=B(x_n,r_n)$, so in particular  $x_{n+1}\leq_\mathbb{N} x$. But this is absurd since the sequence $(x_n)_{n\in\mathbb{N}}$ is strictly increasing and therefore is unbounded. So $X$ is a uniformly $T_3$ scattered CSC space.

\end{proof}
\end{prop}
\begin{defin} A countable metric space $(X,d)$ is said to be totally bounded if for every $n\in \mathbb{N}$ there exists a finite set $F\subseteq X$ such that $X=\bigcup_{x\in F} B(x,\frac{1}{n})$. Being a totally bounded space is arithmetically definable, as noted by Hirst \cite{Hirst2}.
\end{defin}

\begin{lemma} $\textbf{RCA}_0$ proves that if $(X,d)$ is totally bounded then any subset of $X$ is totally bounded.\\
\\
\begin{proof} Let $(X,d)$ be a totally bounded countable metric space and $Y\subseteq X$ be a subspace. By {Proposition~\hyperref[metstruc]{\ref{metstruc}}} there is an $a\in \mathbb{R}_{>0}$ such that the balls $(B(x,q\cdot a))_{x\in X,q\in\mathbb{Q}}$ are clopen. Without loss of generality, we may assume $a=1$. Given $n>0$, we show that we can cover $Y$ with finitely many balls of radius $\frac{1}{n}$. Since $X$ is totally bounded, there exists a finite $F\subseteq X$ such that $X=\bigcup_{x\in F} B(x,\frac{1}{2n})$. Define
\begin{equation*}
G=\left\{y\in Y:\exists x\in F\; y\in B(x,\frac{1}{2n})\wedge \forall z<y\; z\notin Y\cap B(x,\frac{1}{2n})\right\}
\end{equation*}
we have that $G$ exists by $\Delta^0_1$ comprehension and is a finite set since it injects into $F$. For all $x\in F$ if $B(x,\frac{1}{2n})\cap Y\neq \emptyset$ then there is a $y\in G$ such that $B(x,\frac{1}{2n})\subseteq B(y,\frac{1}{n})$ and so $Y\subseteq \bigcup_{y\in G}B(y,\frac{1}{n})$. Thus, $Y$ is totally bounded.
\end{proof}
\end{lemma}

\begin{thm}\label{Geometria 2} $\textbf{ACA}_0$ proves that a countable metric space is compact if and only if it is complete and totally bounded.\\
\\
\begin{proof} Let $(X,d)$ be a compact countable metric space. Since over $\textbf{ACA}_0$ sequential compactness is equivalent to compactness, we have that $X$ is sequentially compact. Therefore, every Cauchy sequence in $X$ will have a convergent subsequence, so $X$ is complete. Since $X$ is compact, for every $n\in\mathbb{N}$, it will have a finite covering of balls of radius $\frac{1}{n}$, so $X$ is totally bounded.\\
\\
For the converse, let $X$ be a complete, non empty, and totally bounded countable metric space. We show that $X$ is sequentially compact, which implies that $X$ is compact. Let $(y_n)_{n\in\mathbb{N}}$ be a sequence in $X$, we define inductively a subsequence $(z_n)_{n\in\mathbb{N}}$ of $(y_n)_{n\in\mathbb{N}}$ and a sequence of points $(x_n)_{n\in\mathbb{N}}$. Let $z_0=y_0$ and $x_0$ be the first $y_n$ such that $B(y_n,1)\cap\{y_k:k\in\mathbb{N}\}$ is infinite. Assume that we have defined $z_i$ and $x_i$ for $i\leq k$ such that for infinitely many $n$ $y_n\in\bigcap_{i\leq k} B(x_i,\frac{1}{2i+1})$. Since $X$ is totally bounded, we have that $\bigcap_{i\leq k} B(x_i,\frac{1}{2i+1})$ is totally bounded, and so there is a least finite set $F$ such that
\begin{equation*}
\bigcap_{i\leq k} B(x_i,\frac{1}{2i+1})\subseteq \bigcup_{x\in F} B(x,\frac{1}{2k+3})
\end{equation*}
By the infinite pigeonhole principle let $x_{k+1}$ be the least $x\in F$ such that for infinitely many $n$ we have $y_n\in B(x,\frac{1}{2k+3})\cap \bigcap_{i\leq k} B_i$. Define $z_{k+1}$ to be the first $y_n\in B(x_n,\frac{1}{2k+3})$.\\
\\
For every $m\in\mathbb{N}$ the sequence $(z_n)_{n\in\mathbb{N}}$ is eventually in $B(x_m,\frac{1}{2m+1})$ which has diameter $\frac{2}{2m+1}<\frac{1}{m}$ and so we have that $(z_n)_{n\in\mathbb{N}}$ is a Cauchy subsequence of $(y_n)_{n\in\mathbb{N}}$. Since $(X,d)$ is complete we have that $(z_n)_{n\in\mathbb{N}}$ converges.
\end{proof}
\end{thm}

\begin{cor}\label{compiscomp} $\textbf{ACA}_0$ proves every effectively $T_2$ effectively compact CSC space is completely metrizable.\\
\\
\begin{proof} By Corollary~\hyperref[acacomp]{\ref{acacomp}} we have that $\textbf{ACA}_0$ proves that every $T_2$ compact CSC space is $T_3$ and by Theorem~\hyperref[stuff>aca]{\ref{stuff>aca}} it is metrizable. Any metric on a compact CSC space must be complete by the previous theorem.
\end{proof}
\end{cor}

\begin{defin}\label{defT+}\textbf{(Simpson \cite[Exercise VI.1.8]{Simp})} Let $T\subseteq \mathbb{N}^{<\mathbb{N}}$ be a tree then we define
\begin{equation*}
T^{+}=\{\tau\in \mathbb{N}^{<\mathbb{N}}: \exists \sigma\in T\;(|\tau|= |\sigma| \wedge \forall n<|\tau|\;(\sigma(n)\leq \tau(n)))\}
\end{equation*}
which is the tree of all sequences $\tau$ that dominate some sequence of $T$ of equal length.
\end{defin}

\begin{prop}\label{T+} Over $\textbf{ACA}_0$ we have that $T$ is well-founded if and only if $\text{KB}(T^{+})$ with the order topology is scattered. \\
\\
\begin{proof} The Kleene-Brouwer order of any tree will have the empty sequence as a maximal element. If $\text{KB}(T^+)$ is not scattered, then it is not compact, so it cannot be well-ordered. This means $T^+$ is not well-founded. Given $f\in [T^+]$ we define the subtree
\begin{equation*}
\{\sigma\in T: \forall n<|\sigma|\,(\sigma(n)\leq f(n))\}
\end{equation*}
which is a finitely branching infinite subtree of $T$. By weak K{\H o}nig's lemma, we have that the subtree, and therefore $T$, also has a branch.\\
\\
Assume instead $T$ is not well-founded and let $f\in[T]$ be a branch. We show that
\begin{equation*}
S=\{\sigma\in T^+:\forall n<|\sigma|\;(\sigma(n)\geq f(n))\}\cong\mathbb{Q}
\end{equation*} 
to do this, it suffices to show that $S$ is dense in itself. Let $\rho\in S$  and $\tau,\sigma \in T^+$ such that $\rho \in \mathopen]\sigma,\tau\mathclose[_{<_{\text{KB}}}$, if $\rho\sqsubseteq \sigma$ then we have
$\sigma<_{\text{KB}}\rho^\frown(\sigma(|\rho|)+1)<_{\text{KB}}\rho$, and so $\rho^\frown(\sigma(|\rho|+1))\in \mathopen]\sigma,\tau\mathclose[_{<_{\text{KB}}}\cap S$. Otherwise, if $\rho\notsqsubseteq \sigma$ then we have $\sigma<_{\text{KB}}\rho^\frown(f(|\rho|)+1)<_{\text{KB}}\rho$. So $S$ with the subspace topology is dense in itself, and so $T^+$ is not scattered.
\end{proof}
\end{prop}

\begin{obs} The previous proposition shows that being a scattered CSC space is a $\Pi^1_1$ universal formula. The proof above can be modified to show that if $T$ is well-founded, then $T^+$ with the topology given by being a subspace of $(\mathbb{N}^{<\mathbb{N}},<_{\text{KB}})$ is compact. Similarly, if $T$ is ill founded then $T^+$ with the topology given by being a subspace of $(\mathbb{N}^{<\mathbb{N}},<_{\text{KB}})$ is scattered.\\
\\
Clote \cite{Clote} used a similar construction to show that being a scattered linear order is a $\Pi^1_1$ universal formula. However, as noted before, the notion of scatteredness for linear orders does not, in general, coincide with the topological notion of scatteredness.
\end{obs}

\begin{prop}\label{completecomplete} Over $\textbf{ACA}_0$ being a complete countable metric space is a $\Pi^1_1$ universal formula.\\
\\
\begin{proof} We can express $(M,d)$ being complete as
\begin{equation*}
\forall Y \subseteq M\; (Y\text{ is a Cauchy sequence}\rightarrow \exists x\in M\;(Y\text{ converges to }x))
\end{equation*}
Which is a $\Pi^1_1$ formula since converging to $x$ and being a Cauchy sequence are arithmetically defined.\\
\\
Fix any arithmetically definable metric $d$ that is compatible with the order topology on $(\mathbb{N}^{<\mathbb{N}},<_\text{KB})$. Over $\textbf{ACA}_0$ being a well-founded tree is $\Pi^1_1$ universal so we have
\begin{equation*}
T\text{ is well-founded}\leftrightarrow  (T^{+},<_{\text{KB}}) \text{ is compact} \rightarrow (T^{+},d)\text{ is complete}
\end{equation*}
and
\begin{equation*}
T\text{ is ill founded}\leftrightarrow (T^{+},<_{\text{KB}}) \text{ is not scattered} \rightarrow (T^{+},d)\text{ is not complete}
\end{equation*}
So over $\textbf{ACA}_0$ we have that for any tree $T$
\begin{equation*}
T \text{ is well-founded}\leftrightarrow  (T^+,d) \text{ is complete}
\end{equation*}
which implies that being a complete countable metric space is a $\Pi^1_1$ universal formula.
\end{proof}
\end{prop}

\begin{cor}\label{scat index} Over $\textbf{RCA}_0$ the following are equivalent:
\begin{enumerate}[label={(\arabic*)}]
\item $\Pi^1_1$ comprehension.
\item For any sequence of CSC spaces $(X^i)_{i\in\mathbb{N}}$ the set $\{i\in\mathbb{N}: X^i \text{ is scattered}\,\}$ exists.
\item For any sequence of countable metric spaces $(X^i,d^i)_{i\in\mathbb{N}}$ the set  $$\{i\in\mathbb{N}: (X^i,d^i) \text{ is complete}\,\}$$ exists.
\end{enumerate}
\begin{proof} Since being a complete countable metric space and being scattered are both $\Pi^1_1$ universal formulas over $\textbf{ACA}_0$, it suffices to show that $(1)$ and $(2)$ both imply arithmetic comprehension. Using Theorem~\hyperref[ACAnot]{\ref{ACAnot}}, we show that $(1)$ and $(2)$ both imply that every set has a Turing jump.\\
\\ 
Assume $(2)$ and let $A\subseteq \mathbb{N}$ be a set. Let $(q_n)_{n\in\mathbb{N}}$ be an enumeration of $\mathbb{Q}$. Define
\begin{equation*}
X^e=\{q_n:  \neg \Phi^A_e(e){\downarrow}_{\leq n}\}\subseteq \mathbb{Q}
\end{equation*}
with the subspace topology. We have that $X^e=\mathbb{Q}$ if and only if $\Phi^A_e(e){\uparrow}$ and $X^e$ is a discrete finite space, and therefore scattered, if and only if $\Phi^A_e(e){\downarrow}$. So $A'\leq_T \{e\in\mathbb{N}: X^e \text{ is scattered}\,\}$ and so $A'$ exists by $\Delta^0_1$ comprehension. \\
\\
Assume $(3)$ and let $A\subseteq \mathbb{N}$ be a set. Define $X^e=\{\frac{1}{n+1}: \neg \Phi^A_e(e){\downarrow}_{\leq n}\}$ with the metric $d^e(x,y)=|x-y|$. We have that $X^e$ is complete if and only if $\Phi^A_e(e){\downarrow}$ and so $A'\leq_T \{e\in\mathbb{N}: (X^e,d^e) \text{ is  complete}\,\}$ exists by $\Delta^0_1$ comprehension.
\end{proof}
\end{cor}

\begin{defin} In general topology, a subset of a topological space $X$ is said to be $G_\delta $ if it is the countable intersection of open sets. If $X$ is countable and $T_1$, that is, all singletons are closed, then all subsets of $X$ are $G_\delta$. A classic result in descriptive set theory is that the $G_\delta$ subspaces of complete separable metric spaces are also completely metrizable. We note that the usual proof can be carried out over $\mathbf{ACA}_0$
\end{defin}

\begin{lemma}\label{Gd} (See \cite[Theorem 3.11]{Kechris}) $\textbf{ACA}_0$ proves every subspace of a complete metric space is completely metrizable.
\end{lemma}

\begin{thm}\label{compmet} Over $\textbf{RCA}_0$ the following are equivalent:
\begin{enumerate}[label={(\arabic*)}]
\item Arithmetic transfinite recursion.
\item Every $T_3$ scattered CSC space can be embedded into a completely metrizable space.
\item Every $T_3$ scattered CSC space is completely metrizable.
\end{enumerate}
\begin{proof}
$(1\rightarrow 2)$ By Lemma~\hyperref[scatter in comp]{\ref{scatter in comp}} every $T_3$ scattered CSC space embeds into a  $T_2$ compact space and by Corollary~\hyperref[compiscomp]{\ref{compiscomp}} every $T_2$ compact space is completely metrizable.\\
\\
$(2)$, $(3)$ both imply arithmetic comprehension since they imply that every $T_3$ scattered CSC space is metrizable and by Theorem~\hyperref[Big4]{\ref{Big4}} and Corollary~\hyperref[scatter metriz>aca]{\ref{scatter metriz>aca}} this implies arithmetic comprehension. So we may prove the remaining implications over $\textbf{ACA}_0$.\\
\\
$(2\rightarrow 3)$ Let $X$ be a $T_3$ scattered CSC space. If $X$ can be embedded into a complete countable metric space, then by Lemma~\hyperref[Gd]{\ref{Gd}} $X$ is completely metrizable.\\
\\
$(3\rightarrow 1)$ By Lemma~\hyperref[lc>T3]{\ref{lc>T3}} and  Proposition~\hyperref[lc are scattered]{\ref{lc are scattered}} we have that every $T_2$ l.c.\ space is $T_3$ and scattered. So $(3)$ implies every $T_2$ l.c.\ CSC space is completely metrizable. Let $X$ be a $T_2$ l.c.\ CSC space and let $d$ be a complete metric on $X$. For each $x\in X$ let $i(x)$ be the least number such that $(\overline{B(x,\frac{1}{2^{i(x)}})})_{i\in\mathbb{N}}$ is totally bounded. Since being totally bounded is arithmetical, the sequence $(i(x))_{x\in X}$ exists by arithmetic comprehension. Since closed subsets of complete countable metric spaces are complete, we have that for all $x$, the set $\overline{B(x,\frac{1}{2^{i(x)}})}$ is compact by Theorem~\hyperref[Geometria 2]{\ref{Geometria 2}}. So the sequence $(\overline{B(x,\frac{1}{2^{i(x)}})})_{i\in\mathbb{N}}$ is a $CCN$. We have that every $T_2$ l.c.\ space has a $CCN$, which is equivalent to arithmetic transfinite recursion.
\end{proof}
\end{thm}
\section{Comparability of locally compact $T_2$ spaces}

A classic result in reverse mathematics is that over $\textbf{RCA}_0$ arithmetic transfinite recursion is equivalent to the comparability of well orders. By comparability of well orders, we mean that for any pair of well orders, there is an order preserving map from one to the other. It turns out there is a similar topological equivalent to arithmetic transfinite recursion. It was shown by Friedman \cite{FriedC} that topological comparability of well orders is equivalent to arithmetic transfinite recursion. We will give an alternative proof to a weaker theorem of Hirst \cite{Hirst2}, namely that arithmetic transfinite recursion is equivalent over $\textbf{ACA}_0$ to the topological comparability of $T_2$ l.c.\ CSC spaces.

\begin{defin} We say two CSC spaces $X$ and $Y$ are topologically comparable if there is an effective embedding from $X$ to $Y$ or from $Y$ to $X$.
\end{defin}

\begin{prop}\label{SUBSCCN} Arithmetic transfinite recursion is equivalent over $\textbf{ACA}_0$ to every $T_2$ l.c.\ CSC space can be embedded into a compact space.\\
\\
\begin{proof} By Theorem~\hyperref[TSP>CCN]{\ref{TSP>CCN}} we have that $\textbf{ATR}_0$ proves that every $T_2$ l.c.\ space has a $CCN$. By Proposition~\hyperref[alexemb]{\ref{alexemb}} every l.c.\ $T_2$ CSC space with a $CCN$ embeds into its one point compactification. On the other hand, every compact space trivially has $CCN$. By Proposition~\hyperref[CCNstab]{\ref{CCNstab}}, if every l.c.\ $T_2$ CSC space can be embedded into a compact CSC space, then it will have a $CCN$, which by {Proposition~\hyperref[T2CCN>ATR]{\ref{T2CCN>ATR}}} implies arithmetic transfinite recursion.
\end{proof}
\end{prop} 

\begin{lemma}\label{weball} $\textbf{ACA}_0$ proves that for any CSC space  $(X,(U_i)_{i\in\mathbb{N}},k)$ if every $T_2$ compact CSC space embeds into $X$ then $X$ is not scattered.\\
\\
\begin{proof} Let $X$ be a CSC space such that every compact $T_2$ CSC space $C$ embeds into $X$. For every well-founded tree $T$, we have that $T^{+}$ is well-founded and therefore, $(T^{+},<_{\text{KB}})$ (see {Definition~\hyperref[defT+]{\ref{defT+}}} ) embeds into $X$ since it is compact. So we have that
\begin{equation*}
\forall T\subseteq \mathbb{N}^{<\mathbb{N}} ( (T,<_{\text{KB}})\text{ is a well order}\rightarrow (T^{+},<_{\text{KB}}) \text{ embeds into } X)
\end{equation*}
But $(T^{+},<_{\text{KB}})$ embeds into $X$ is $\Sigma^1_1$ relative to $X$. Since ``$(T,<_{\text{KB}})$ is a well order" is a universal $\Pi^1_1$ formula we have that there exists a non well-founded $T$ such that $(T^{+},<_{\text{KB}})$ embeds into $X$. Since $T$ is not well-founded by the {Proposition~\hyperref[T+]{\ref{T+}}}, we have that  $\mathbb{Q}$ embeds into $(T^{+},<_{\text{KB}})$ and so $\mathbb{Q}$ embeds into $X$ contradicting our assumption that $X$ was scattered.
\end{proof}
\end{lemma}

\begin{thm} Over $\textbf{ACA}_0$ the following are equivalent:
\begin{enumerate}[label={(\arabic*)}]
\item Arithmetic transfinite recursion.
\item All pairs of l.c.\ $T_2$ CSC spaces are topologically comparable.
\end{enumerate}
\begin{proof} $(1\rightarrow 2)$ follows immediately from the fact that 
$\textbf{ATR}_0$ proves that every $T_2$ l.c.\ CSC space is homeomorphic to some well order and that for every pair of well orders, $X$ and $Y$ either $X$ is isomorphic to an initial segment of $Y$ or $Y$ isomorphic to an initial segment of $X$.\\
\\
$(2\rightarrow 1)$. By {Corollary~\hyperref[SUBSCCN]{\ref{SUBSCCN}}}, arithmetic transfinite recursion is equivalent to every $T_2$ l.c.\ CSC space can be embedded into a $T_2$ compact CSC space. Let $X$ be a $T_2$ l.c.\ CSC space. If $X$ does not embed into a $T_2$ compact CSC space, we have that every compact $T_2$ CSC space is homeomorphic to a subspace of $X$. By the previous lemma, this implies that $X$ is not scattered, which by {Lemma~\hyperref[lc are scattered]{\ref{lc are scattered}}} contradicts our assumption that $X$ is l.c.\
\end{proof}
\end{thm}

\begin{obs} We might wonder if $T_3$ scattered spaces are topologically comparable. This is not the case; we can consider the ordinal $\omega\cdot 2+1$ and a point with infinitely many sequences converging to it. These two spaces are both $T_3$ and scattered.  However, they are not topologically comparable. This observation can be carried out easily in $\textbf{RCA}_0$. However, Friedman proved that over $\mathbf{ATR}_0$ countable metric spaces are comparable with continuous injections.
\end{obs}

\begin{thm}\label{balling} \textbf{(Friedman \cite[Theorem 16]{FriedC})}  $\textbf{ATR}_0$ proves that for any two countable metric space there is a continuous injection from one to the other.
\end{thm}
\noindent
We can use this fact to get an alternative proof that every $T_3$ scattered CSC space embeds into a well order. We provide a proof sketch.
\begin{prop} $\textbf{ATR}_0$ proves that any scattered $T_3$ CSC space embeds into a well order.\\
\\
\begin{proof} Let $X$ be a $T_3$ scattered CSC space. If there is a continuous injection from $X$ to a well order then, by an argument similar to that in {Theorem~\hyperref[fruit]{\ref{fruit}}}, there is an embedding from $X$ into a well order. Assume that there aren't any continuous injections from $X$ into a well order, then for every well order $W$ there exists a continuous injection from $W$ to $X$ by {Theorem~\hyperref[balling]{\ref{balling}}}. In particular, for every well order $W$ with maximum, $W$ will have compact order topology. So any continuous injection from $W$ to $X$ will be an embedding.  By {Corollary~\hyperref[corwell]{\ref{corwell}}}, every $T_2$ compact CSC space is well orderable over $\textbf{ACA}_0$, so every compact $T_2$ CSC space embeds into $X$. By {Lemma~\hyperref[weball]{\ref{weball}}}, we have that $X$ is not scattered.
\end{proof}
\end{prop}
\noindent
For completeness, we will lay out Friedman's proof that topological comparability for well orders is equivalent to arithmetic transfinite induction. To do so, we need to define well order exponentiation.

\begin{defin}\textbf{(Hirst \cite[Definition 2.1]{Hirst})} Given two well orders $L$ and $W$, we define $W^L$ to be the set of sequences that includes the empty sequence and all sequences of the form
\begin{equation*}
((a_0,b_0),\dots (a_n,b_n))
\end{equation*}
such that $\forall i\leq n$ $b_i\in L$ and $a_i\in W\setminus\{\min W \}$ and for all $j<i\leq n$ $b_i<_L b_j$.\\
\\
Given $\sigma,\tau\in W^L$ we define $\sigma<_{W^L}\tau$ if and only if either $\sigma\subsetneq \tau$ or given that
\begin{equation*}j=\max\{i<\min\{|\tau|,|\sigma|\},\sigma(i)\neq\tau(i)\}
\end{equation*}
and 
\begin{equation*}
(a_j,b_j)=\sigma(j)\wedge (c_j,d_j)=\tau(j)
\end{equation*} then
\begin{equation*}(b_j<d_j\vee (d_j=b_j\wedge a_j<c_j)) 
\end{equation*} 
Intuitively, the elements of $W^L$ can be viewed as being ordinals less than $W^L$ in their Cantor normal form in base $W$ and ordered in the standard way. 
\end{defin}

\begin{thm}\textbf{(Hirst \cite[Theorem 2.6]{Hirst})} Over $\textbf{RCA}_0$, arithmetic comprehension is equivalent to ordinal exponentiation being well defined. That is, for any pair of well orders $W$ and $L$, the set $W^L$ is well-ordered.
\end{thm}

\begin{obs}\label{rankobs} $\textbf{ACA}_0$ proves that the isolated points of $\mathbb{N}^L$ are the empty sequence and all sequences of the form
\begin{equation*}
((a_0,b_0),\dots (a_n,0))
\end{equation*}
since they have as a successor
\begin{equation*}
((a_0,b_0),\dots (a_n+1,0))
\end{equation*}
and as a predecessor
\begin{equation*}
\begin{cases}((a_0,b_0),\dots (a_n-1,0))\quad  \text{ if }\quad a_n>1\\
((a_0,b_0),\dots (a_{n-1},b_{n-1})) \quad \text{ if } \quad a_n=1
\end{cases}
\end{equation*}
On the other hand, sequences of the form
\begin{equation*}
((a_0,b_0),\dots (a_n,b_n))
\end{equation*}
where $b_n>0$ will not have a predecessor and are not isolated.
\end{obs}

\begin{lemma}\label{NLrank} $\textbf{ACA}_0$ proves that for every well order $L$ $\mathbb{N}^L+1$ is a well order and has $L$ as a rank.\\
\\
\begin{proof} 
For each $l\in L$ let $X_l$ be the set of all sequences in $\mathbb{N}^L$ of the form
\begin{equation*}
((a_0,b_0),\dots (a_n,b_n))
\end{equation*}
such that for all $i\leq n$ we have $b_i\geq l$. We have that the sequence $(X_l)_{l\in L}$ exists by arithmetical comprehension. By {Observation~\hyperref[rankobs]{\ref{rankobs}}}, we have that $X_{l+1}$ are the limit points of $X_l$. So we have that the sequence $(X_l)_{l\in L}$  witnesses that $L$ is a rank for $\mathbb{N}^L$. 
\end{proof}
\end{lemma}

\begin{lemma}\label{Homorank} Over $\textbf{ACA}_0$, let $X$ and $Y$ be CSC spaces with rank $W$ and $Y$ with respectively. If $X$ embeds topologically into $Y$ then there exists an order embedding of $W$ into $L$.\\
\\
\begin{proof} Let $f:X\rightarrow Y$ be an embedding, then we define $\phi:W\rightarrow L$ where given $w\in W$:
\begin{equation*}
\phi(w)=\min\{l\in L: \exists x\in X (\text{rank}_W(x)=w \wedge \text{rank}_L(f(x))=l\}.
\end{equation*}
We show that $f$ preserves the strict order. Seeking a contradiction, assume that there exists a least $w\in W$ and $w_0<_W w$ such that $\phi(w)\leq_L \phi(w_0)$. Let $x\in X$ be of rank $w$ such that $f(x)$ has rank $\phi(w)$. In $X_{w_0}$, $x$ is a limit point of the elements of rank $w_0$ and so $f(x)$ must be the limit point of the image of the elements of rank $w_0$. We have by definition of $\phi$ that $f(W_{w_0})\subseteq Y_{\phi(w_0)}$, so $f(x)$ is a limit point of $Y_{\phi(w_0)}$. Since $\text{rank}_{L}(f(x))=\phi(w)\leq_L\phi(w_0)$ either $f(x)\notin Y_{\phi(w_0)}$ or $f(x)$ is isolated in $Y_{\phi(w_0)}$. Both cases contradict the fact that $f(x)$ is a limit point of $Y_{\phi(w_0)}$
\end{proof}
\end{lemma}
\noindent
Friedman proved a special case of the previous theorem, where $X$ and $Y$ are of the form $\mathbb{N}^{W_0}$ and $\mathbb{N}^{W_1}$, where $W_0$ and $W_1$ are well orders (see \cite[Lemma 25]{FriedC}).
\begin{thm} Over $\textbf{RCA}_0$ the following are equivalent:
\begin{enumerate}[label={(\arabic*)}]
\item Arithmetic transfinite recursion.
\item Every pair of $T_2$ l.c.\ CSC spaces are topologically comparable. (Hirst \cite[Theorem 4.1]{Hirst2})
\item Every pair of well orders are topologically comparable. (Friedman \cite[Lemma 20]{FriedC})
\end{enumerate}
\begin{proof} We already showed that $(1\rightarrow 2)$, and $(2\rightarrow 3)$ is obvious.\\
\\
The proof that $(3)$ implies arithmetic comprehension is not trivial but standard and can be found in \cite{FriedC}. Working over $\textbf{ACA}_0$ we show $(3)$ implies that for any two well-orders $W$ and $L$, either there is a strictly increasing function from $W$ to $L$ or from $L$ to $W$, which is equivalent to arithmetic transfinite recursion. We have by {Lemma~\hyperref[NLrank]{\ref{NLrank}}} that $\mathbb{N}^L+1$ has $L$ as a rank and $\mathbb{N}^W+1$ has $W$ as a rank. $\textbf{ACA}_0$ proves that $\mathbb{N}^L+1$ and $\mathbb{N}^W+1$ are well-orders and so there exists either an embedding from $\mathbb{N}^L+1$ to $\mathbb{N}^W+1$ or an embedding from $\mathbb{N}^W+1$ to $\mathbb{N}^L+1$. Without loss of generality, assume there is an embedding from $\mathbb{N}^L+1$ to $\mathbb{N}^W+1$. By {Lemma~\hyperref[Homorank]{\ref{Homorank}}}, this means that there is a strictly increasing function from $L$ to $W$.
\end{proof}
\end{thm}
\noindent
Montalb{\'a}n and Greenberg \cite{Greenberg} showed that $\textbf{ATR}_0$ is equivalent to every $T_2$ compact CSC space having a rank. We give an alternative proof.
\begin{thm} Over $\textbf{RCA}_0$ the following are equivalent:
\begin{enumerate}[label={(\arabic*)}]
\item Arithmetic transfinite recursion.
\item Every scattered $T_3$ space has Cantor-Bendixson rank. 
\item Every well-order with the order topology has Cantor-Bendixson rank.
\end{enumerate}
\begin{proof} $(1\rightarrow 2)$ follows from {Theorem~\hyperref[ATR>rank]{\ref{ATR>rank}}} and $(2\rightarrow 3)$ is immediate.\\
\\
By Proposition~\hyperref[rankwo>aca]{\ref{rankwo>aca}} $(3)$ implies arithmetic comprehension, so we may work over $\textbf{ACA}_0$. Let $W$ and $L$ be well-orders. We show that one embeds as an initial segment of the other, which over $\textbf{RCA}_0$ is equivalent to arithmetic transfinite recursion. Let $X=(\mathbb{N}^L+1)$ and $Y=(\mathbb{N}^W+1)$. The space $Z$, which is the disjoint union of $X$ and $Y$, is homeomorphic to the well-order $(\mathbb{N}^W+1)+(\mathbb{N}^L+1)$ and so by assumption $Z$ has a rank $R$. We have that $(\mathbb{N}^W+1)$ has $W$ as rank and $(\mathbb{N}^L+1)$ has $L$ as rank. By transfinite induction on $L$, we have that the map $L\rightarrow R$ given by $\text{rank}_R(x)\mapsto\text{rank}_L(x)$ is well defined, and it defines an order isomorphism from $L$ to an initial segment of $R$. Similarly, we can embed $W$ as an initial segment of $R$. So either $W$ or $L$ is order isomorphic to an initial segment of the other.
\end{proof}
\end{thm}

\section{$\Pi^1_1$ comprehension}
\begin{q}\textbf{(Chan \cite{Chan})} 
Are the following equivalent over
 $\textbf{RCA}_0$?
\begin{enumerate}[label={(\arabic*)}]
\item $\Pi^1_1$ comprehension.
\item For every sequence $((X^j,(U^j_i)_{i\in\mathbb{N}}))_{j\in\mathbb{N}}$ of CSC spaces the set $\{j\in\mathbb{N}: X^j \text{ is connected}\}$ exists.
\end{enumerate}
Where a CSC space $(X,(U_i)_{i\in\mathbb{N}},k)$ is said to be connected for any open set $A$ and $B$ such that $X=A\cup B$ then $A\cap B\neq \emptyset$.
\end{q}
\noindent
Chan showed that $(2)$ implies arithmetic comprehension, so the following result gives a positive answer to the question.
\begin{thm}\label{chansol} Being a connected CSC space is $\Pi^1_1$ universal over $\textbf{ACA}_0$.\\
\\
\begin{proof} Let $L$ be a linear order. Up to adding a new element, we may assume that $L$ has a minimal element. For each $l\in L$ let $V_l$ be the set of all $j\in L$ such that either $j\leq_L l$ or there exists a sequence $\sigma$ such that $\sigma(0)=l$ for all $n<|\sigma|$ we have that $\sigma(n+1)$ is the smallest element greater than $\sigma(n)$ and that $\sigma(|\sigma|-1)=j$. That is, $V_l$ is the smallest set containing $l$ that is closed under successor when it's defined and downwards closed. Define the topology $Top(L)$ on $L$ to be generated by the sets $(V_l)_{l\in L}\cup (\,\mathopen]a,+\infty\mathclose[\,)_{l\in L}$. We observe that sets of the form $V_l\cap\mathopen]j,+\infty\mathclose[$ form a basis for this topology and that all open sets are closed under successor when it's defined.\\
\\
We show that $L$ is well-ordered if and only if it is connected with respect to $Top(L)$. If $L$ has a descending sequence $(a_i)_{i\in\mathbb{N}}$ then the set ${\uparrow} \{a_i:i\in\mathbb{N}\}$ is clopen since its complement is downwards closed and closed under successor. We have that $L\neq {\uparrow}\{a_i:i\in\mathbb{N}\}$ since $L$ is assumed to have a minimal element.\\
\\
If $L$ is a well-order, let $A$ be a non empty clopen set of $L$ which contains the least element of $L$. If $A\neq L$, then let $x$ be the least element of $L\setminus A$. Since $L\setminus A$ is open, we have that $x$ must be in a basic open set of the form $V_j\cap \mathopen]l,+\infty\mathclose[\subseteq L\setminus A$. We have that $l<x$, so by minimality of $x$, we have that $l\in A$. Since $A$ is closed under successor we have $l+1\in A$ but $l+1\in \mathopen]l,+\infty\mathclose[$ and since $V_j$ is downwards closed and $l+1\leq x$ we have that $l+1\in V_j\cap \mathopen]l,+\infty\mathclose[\subseteq L\setminus A$ which is absurd since $l+1\in A$. So $L$ must be connected as it has only trivial clopen sets.\\
\\
So we have that $L$ is a well-order if and only if $Top(L)$ is a connected topology on $L$. Since being a connected space is $\Pi_1^1$ definable, it is a $\Pi^1_1$ universal formula.
\end{proof}
\end{thm}

\begin{thm}\label{basically} Over $\textbf{RCA}_0$ the following are equivalent:
\begin{enumerate}[label={(\arabic*)}]
\item $\Pi^1_1$ comprehension.
\item Every $T_3$ CSC space is the disjoint union of a scattered subspace and a dense in itself subspace (this will not, in general, be a topological disjoint union as the dense in itself subspace may not be open).
\item Every $T_3$ space is ranked.
\end{enumerate}
\begin{proof}
$(1)$ implies arithmetic comprehension since arithmetic formulas are $\Pi^1_1$ and $(3)$ implies arithmetic comprehension by \hyperref[rankwo?aca]{\ref{rankwo>aca}}. We show that $(2)$ implies arithmetic comprehension. Let $A\subseteq\mathbb{N}$ be a set and $(q_n)_{n\in\mathbb{N}}$ be an enumeration of $\mathopen]0,1\mathclose[_{<_\mathbb{Q}}\subseteq \mathbb{Q}$ and consider
\begin{equation*}X=\{(e,q_n):e\in\mathbb{N}\wedge \neg\Phi^A_e(e){\downarrow}_{\leq n}\}\cup\{(e,0):e\in \mathbb{N}\}\subseteq \mathbb{N}\times \mathbb{Q}
\end{equation*}
with the subspace topology. Let $D\subseteq X$ be a maximal dense in itself subspace of $X$ then we have that $(e,0)\notin D\leftrightarrow \Phi^A_e(e){\downarrow}$ and so $A'\leq_T D$. So $(2)$ implies the Turing jump of every set exists, which by \hyperref[ACAnot]{\ref{ACAnot}} implies arithmetic comprehension. Since $(1),(2)$ and $(3)$ all imply arithmetic comprehension, it suffices to show that they are equivalent over $\textbf{ACA}_0$.
\\
\\
$(1\rightarrow 2)$ Let $X$ be a $T_3$ space and define $U\subseteq X$ to be the set of points $x\in X$ such that there exists an infinite dense in itself subspace of $X$ containing $x$. The set $U$ exists by $\Sigma^1_1$ comprehension, which is equivalent to $\Pi^1_1$ comprehension. We have that $U$ is a maximal dense in itself subspace of $X$, and therefore, $X\setminus U$ must be scattered.\\
\\
$(2\rightarrow 1)$ Since being scattered is $\Pi^1_1$ universal, we show that for any sequence $(X_i)_{i\in\mathbb{N}}$ of $T_3$ CSC spaces the set of indices $i$ such that $X_i$ is scattered exists. Let $X=\coprod X_i$ be the disjoint union of all the $X_i$ and let $D$ be the maximal dense in itself set of $X$. We have that $X_i$ is scattered if and only if $X_i\cap D= \emptyset$ and so the set
\begin{equation*}
\{i\in\mathbb{N}: X_i \text{ is scattered }\}
\end{equation*}
exists by arithmetic comprehension.
\\
\\
$(1+2\rightarrow 3)$ Let $X$ be a $T_3$ CSC space, then we have by $(2)$ that $X=P\sqcup Y$ where $D$ is dense in itself and $Y$ is scattered. Since $\textbf{ATR}_0$ proves $T_3$ scattered CSC space is ranked we have that $(1)$ implies that $Y$ has a rank $R$, let $(Y_r)_{r\in R}$ witness that $R$ is a rank for $Y$. We have that $(Y_r\sqcup P)_{r\in R}$ witnesses that $R$ is a rank for $X$.\\
\\
$(3\rightarrow 2)$ Let $X$ be a $T_3$ CSC space. By assumption $X$ has a rank $R$, let $(X_r)_{r\in R}$ witness that $R$ is a rank for $X$. By definition of rank \hyperref[rankdef]{\ref{rankdef}}, we have that $\bigcap_{r\in R} X_r=P$ does not have limit points. So $P$ is either empty or is a maximal dense in itself subspace of $X$. So $X\setminus P$ is scattered and $X=P\sqcup (X\setminus P)$ can be written as the disjoint union between a dense in itself subspace and a scattered subspace.
\end{proof}
\end{thm}
\noindent
\textbf{Note:} The fact that a countable subset of $\mathbb{R}^n$ can be written as the sum of a scattered set and a dense in itself set is due to Sierpinski \cite{Sierp2}. A similar proof is used to show that $\Pi^1_1$ comprehension is equivalent to the Cantor-Bendixson theorem for complete separable metric spaces \cite[Theorem VI.1.3]{Simp}.

\section{Summary Part II}
\begin{thm} Over $\textbf{RCA}_0$ the following are equivalent:
\begin{enumerate}[label={(\arabic*)}]
\item Arithmetic transfinite recursion.
\item Every locally compact CSC space has a $CCN$ (Proposition~\hyperref[CCN>ACA]{\ref{CCN>ACA}} and \hyperref[TSP>CCN]{\ref{TSP>CCN}}).
\item Every $T_2$ locally compact CSC space has a $CCN$ (Proposition~\hyperref[T2CCN>ATR]{\ref{T2CCN>ATR}}).
\item Every $T_2$ locally compact CSC space is the disjoint union of open compact sets (Proposition~\hyperref[2 for 1]{\ref{2 for 1}}).
\item Every $T_2$ locally compact CSC space has an exhaustion of compact sets (Proposition~\hyperref[Chonk]{\ref{Chonk}}).
\item Every $T_2$ locally compact CSC space has an effectively $T_2$ $1$ point compactification (Proposition~\hyperref[awful]{\ref{awful}}).
\item Every $T_2$ locally compact CSC space is homeomorphic to a well-order with the order topology (Proposition~\hyperref[2 for 1]{\ref{2 for 1}}).
\item Every well-order has Cantor-Bendixson rank (Friedman \cite{FriedC}).
\item Every scattered $T_3$ space has Cantor-Bendixson rank (Friedman  \cite{FriedC}).
\item Every pair of $T_2$ locally compact CSC spaces are topologically comparable (Hirst \cite{Hirst2}).
\item Every pair of well-orders with the order topology are topologically comparable (Friedman \cite{FriedC}).
\item Bel'nov-Knaster-Urbanik theorem(\cite{Knaster} and \cite{Belnov}): every $T_3$ scattered CSC space effectively embeds into a well-order (Theorem~\hyperref[scatter in comp]{\ref{scatter in comp}}).
\item Kat{\v e}tov's theorem (\cite{Katetov}): every $T_3$ scattered CSC space is the refinement of the order topology of a well-order (Theorem~\hyperref[fruit]{\ref{fruit}}).
\item Every $T_3$ scattered CSC space is effectively homeomorphic to a scattered linear order with countably many cuts (Theorem~\hyperref[scatscat]{\ref{scatscat}}).
\item Every $T_3$ scattered CSC space is completely metrizable (Theorem~\hyperref[compmet]{\ref{compmet}}).

\end{enumerate}
\end{thm}

\begin{thm} Over $\textbf{RCA}_0$ the following are equivalent:
\begin{enumerate}[label={(\arabic*)}]
\item $\Pi^1_1$ comprehension.
\item For any sequence of CSC spaces $(X_i)_{i\in\mathbb{N}}$ then the set $\{i\in\mathbb{N}:X_i\text{ is connected}\}$ exists (Theorem~\hyperref[chansol]{\ref{chansol}}).
\item For any sequence of CSC spaces $(X_i)_{i\in\mathbb{N}}$ then the set $\{i\in\mathbb{N}:X_i\text{ is compact}\}$ exists (Proposition~\hyperref[univcomp]{\ref{univcomp}}).
\item For any sequence of CSC spaces $(X_i)_{i\in\mathbb{N}}$ then the set $\{i\in\mathbb{N}:X_i\text{ is scattered}\}$ exists (Corollary~\hyperref[scat index]{\ref{scat index}}).
\item For any sequence of countable metric spaces $((X_i,d_i))_{i\in\mathbb{N}}$ then the set $\{i\in\mathbb{N}:(X_i,d_i)\text{ is Complete}\}$ exists (Corollary~\hyperref[scat index]{\ref{scat index}}).
\item For any sequence of countable linear orders $((L_i,<_i))_{i\in\mathbb{N}}$ then the set $\{i\in\mathbb{N}:L_i\text{ is Complete}\}$ exists (Proposition~\hyperref[univcomp]{\ref{univcomp}}, Shafer \cite[Corollary 4.2]{Shafer})
\item Every $T_3$ CSC space is the disjoint union of a scattered space and a space homeomorphic to $\mathbb{Q}$ (Theorem~\hyperref[basically]{\ref{basically}}).
\item Every $T_3$ CSC space is ranked  (Theorem~\hyperref[basically]{\ref{basically}}).
\end{enumerate}
\end{thm}

\begin{q} Over $\textbf{RCA}_0$, any sequentially compact metric space will be complete. Is the same true for effectively compact or even strongly effectively compact metric spaces? If not, what is the strength of every effectively compact metric space is complete?
\end{q}

\begin{q} Does $\mathbf{RCA}_0$ prove that every compact CSC space is scattered?
\end{q}

\begin{q} Is the well-foundedness of $\omega^\omega$ necessary to prove that every strongly effectively compact effectively $T_2$ CSC space is well orderable over $\mathbf{RCA}_0$?
\end{q}

\begin{q} Does arithmetic comprehension follow from every uniformly $T_3$ scattered CSC space having a complete metric?
\end{q}

\begin{q} Does arithmetic comprehension follow from every $T_3$ scattered CSC space having an effective embedding into a $T_2$ compact CSC space?
\end{q}

\begin{q} Over $\textbf{ACA}_0$ we have
\begin{equation*}
\text{Compact} \rightarrow \text{Locally Compact}\rightarrow \text{Scattered}
\end{equation*}
and being scattered and being compact are both $\Pi^1_1$ universal formulas. We can express a CSC space  $(X,(U_i)_{i\in\mathbb{N}},k)$ is locally compact as
$
\forall x\,\exists i\,(x\in U_i\wedge \overline{U_i}\text{ is compact}\,)
$.
So over $\Sigma^1_1$ axiom of choice being locally compact is equivalent to a $\Pi^1_1$ formula. Is being locally compact equivalent to a $\Pi^1_1$ formula over $\textbf{ACA}_0$?
\end{q}

\begin{q} $\textbf{ATR}_0$ proves that any $T_3$ scattered CSC space is homeomorphic to the order topology of a scattered linear order with countably many cuts. Can $\textbf{ACA}_0$ prove that any $T_3$ scattered CSC space is homeomorphic to a scattered linear order with the order topology? 
\end{q}
\begin{q}\textbf{(Shafer)} What are the sets of Turing degrees of metrics on a recursive $T_3$ CSC space? An analysis of the proof in Theorem \hyperref[worst proof]{\ref{worst proof}} shows that there is a recursive $T_3$ space $X$ such that a set $A$ computes a metric for $X$ if and only if $\emptyset'\leq_T A$.
\end{q}
\noindent
\textbf{Acknowledgements:} Many thanks to my supervisor, Paul Shafer, for many helpful suggestions and discussions,  for reviewing my drafts, and for presenting the original problem and motivation for this work. I would also like to thank Carl Mummert, Orazio Nicolosi, and the anonymous reviewer for helpful suggestions and references.\\
\\
\textbf{Funding:}
The author's work was supported by a scholarship from the UK Engineering and Physical Sciences Research Council.

\pagestyle{fancy}
\renewcommand{\headrulewidth}{0pt}
\renewcommand{\footrulewidth}{0pt}
\fancyhf{}
\pagenumbering{gobble}
\printbibliography
\fancyfoot[L]{School of Mathematics, University of Leeds\\
url: \url{ https://giorgioggenovesi.github.io/ } \\
email: \url{mmggg@leeds.ac.uk}  }
\end{document}